\documentclass[12pt,oneside,english]{amsart}
\usepackage[T1]{fontenc}
\usepackage[latin1]{inputenc}
\usepackage[a4paper]{geometry}
\usepackage{textcomp}
\usepackage{amsthm}
\usepackage{amstext}

\makeatletter
\numberwithin{equation}{section}
\numberwithin{figure}{section}
  \theoremstyle{plain}
  \newtheorem*{conjecture*}{Conjecture}
\theoremstyle{plain}
\newtheorem{thm}{Theorem}
  \theoremstyle{plain}
  \newtheorem{cor}[thm]{Corollary}
  \theoremstyle{definition}
  \newtheorem{defn}[thm]{Definition}
  \theoremstyle{plain}
  \newtheorem{lem}[thm]{Lemma}
  \theoremstyle{plain}
  \newtheorem{prop}[thm]{Proposition}
  \theoremstyle{remark}
  \newtheorem{rem}[thm]{Remark}
  \theoremstyle{plain}
  \newtheorem{conjecture}[thm]{Conjecture}

\usepackage{amsthm}

\usepackage{amssymb}

\usepackage{amssymb}

\usepackage{amssymb}

\def\makebbb#1{
    \expandafter\gdef\csname#1\endcsname{
        \ensuremath{\Bbb{#1}}}
}
\makebbb{R}
\makebbb{N}
\makebbb{Z}
\makebbb{C}
\makebbb{G}
\makebbb{E}
\makebbb{H}
\makebbb{P}
\makebbb{B}
\makebbb{K}
\makebbb{E}

\usepackage{babel}

\usepackage{babel}

\usepackage{babel}

\usepackage{babel}

\makeatother

\usepackage{babel}

\begin{document}

\title{analytic torsion, vortices and positive Ricci curvature }

\author{Robert J. Berman}

\curraddr{Mathematical Sciences, Chalmers University of Technology and the
University of Gothenburg, SE-412 96 Göteborg, Sweden}

\email{robertb@chalmers.se}
\begin{abstract}
We characterize the global maximizers of a certain non-local functional
defined on the space of all positively curved metrics on an ample
line bundle $L$ over a Kähler manifold $X.$ This functional is an
adjoint version, introduced by Berndtsson, of Donaldson's $L$-functional
and generalizes the Ding-Tian functional whose critical points are
Kähler-Einstein metrics of positive Ricci curvature. Applications
to (1) analytic torsions on Fano manifolds (2) Chern-Simons-Higgs
vortices on tori and (3) Käher geometry are given. In particular,
proofs of conjectures of (1) Gillet-Soulé and Fang (concerning the
regularized determinant of Dolbeault Laplacians on the two-sphere)
(2) Tarantello and (3) Aubin (concerning Moser-Trudinger type inequalities)
in these three settings are obtained. New proofs of some results in
Kähler geometry are also obtained, including a lower bound on Mabuchi's
$K-$energy and the uniqueness result for Kähler-Einstein metrics
on Fano manifolds of Bando-Mabuchi. This paper is a substantially
extended version of the preprint \cite{berm} which it supersedes.
\end{abstract}
\maketitle
\tableofcontents{}

\section{Introduction}

\subsection{Background }

Consider the two-dimensional sphere $S^{2}$ equipped with its standard
Riemannian metric $g_{0}$ of constant positive curvature, normalized
so that the corresponding volume form $\omega_{0}$ gives unit volume
to $S^{2}.$ A celebrated inequality of Moser-Trudinger-Onofri proved
in its sharp form by Onofri \cite{on}, asserts that \begin{equation}
\log\int_{S^{2}}e^{-u}\omega_{0}\leq-\int_{S^{2}}u\omega_{0}+\frac{1}{4}\int_{S^{2}}du\wedge d^{c}u\label{eq:m-t-o intro}\end{equation}
 for any, say smooth, function $u$ on $S^{2},$ where the last term
is the $L^{2}-$norm of the gradient of $u$ in the conformally invariant
notation of section \ref{sub:Setup} below. 

As is well-known the inequality above has a rich geometric content
and appears in a number of seemingly unrelated contexts ranging from
the problem of prescribing the Gauss curvature in a conformal class
of metrics on $S^{2}$ (the Yamabe and Nirenberg problems \cite{cha})
to sharp critical \emph{Sobolev inequalites} \cite{bec} and lower
bounds on \emph{free energy functionals} in mathematical physics \cite{on,k,clmp}.
There is also a \emph{spectral} intepretation of the Moser-Trudinger-Onofri
inequality \cite{on,o-p-s}: the {}``determinant'' $\det\Delta_{g},$
of the Laplacian on $S^{^{2}}$ seen as a functional on space of all
unit area conformal metrics $g$ on $S^{2},$ is maximized precisely
on the standard round metric - up to conformal automorphisms. More
precisely, $\det\Delta_{g}$ is defined using the zeta function regularization
of the product of all strictly positive eigenvalues and its maximizers
are hence the metrics of the form $e^{-u}g_{0}$ satisfying the the
constant positive curvature equation: \[
\omega_{0}+dd^{c}u=e^{-u}\omega_{0},\]
where $dd^{c}u=\Delta_{g_{0}}u\omega_{0}/4\pi$ (using the notation
in section \ref{sub:Setup}).The bridge between this spectral problem
and the inequality \ref{eq:m-t-o intro} is given by the \emph{Polyakov
anomaly formula} \cite{cha}, which first appeared in Physics in the
path integral (random surface) approach to the quantization of the
bosonic string.

One of the main motivations for the present paper comes from a conjecture
of of Gillet-Soulé concerning the determinant of the Dolbeault Laplacian
acting on the space of $(0,q)-$forms of a vector bundle over a Kähler
manifold. In the one-dimensional case it can be formulated as follows:
\begin{conjecture*}
(Gillet-Soulé). Let $(X,\omega_{0})$ be a complex curve with a fixed
Hermitian metric. The determinant of the Dolbeault Laplacian $\Delta_{\bar{\partial}}$
considered as a functional on the space of all smooth Hermitian metrics
on a holomorphic line bundle $L\rightarrow X$ is bounded from above.
\end{conjecture*}
The conjecture was motivated by Arakelov geometry, in particular the
arithmetic Riemann-Roch theorem and is equivalent to the boundedness
(from below) of certain arithmetic Betti-numbers (\cite{g-s}; see
also \cite{g-sa} p. 526-527). In the case of $S^{2}$ the conjecture
was confirmed by Fang \cite{f}, who by symmetrization reduced the
problem to the case when the metric on $L$ is invariant under rotation
around an axes of $S^{2},$ earlier treated by Gilllet-Soulé \cite{g-s}.
Fang also put forward the following more precise form of the conjecture
above, when $L$ is ample:
\begin{conjecture*}
(Fang). In the case of $(S^{2},\omega_{0})$ equipped with the standard
round metric $\omega_{0}$ the upper bound in the previous conjecture
is achieved precisely for the Fubini-Study metric on $L$ (up to scaling).
\end{conjecture*}
In other words the maximizers are conjectured to be precisely the
metrics on the line bundle $L$ whose curvature form is constant in
the sense that it is given by the invariant metric $\omega_{0}$ (assuming
that $L$ has positive degree). Fang was also motivated by the similarities
between the conjecture above and a classical inequality of Szegö for
Toeplitz operators on the unit-circle. It should also be pointed out
that the conjecture above implies, by well-known arguments in Arakelov
geometry involving Minkowski's theorem for lattices,  an effective
arithmetic Riemann-Roch theorem on the projective line \cite{g-sa}
(for a precise statement see the preprint \cite{berm}, section 3.5).

\subsection{The present paper}

In this paper the positive solution of Fang's conjecture will be deduced
from a general result about the maximizers of a certain non-local
functional $\mathcal{F}_{\omega_{0}}$ defined on the space of all
positively curved Hermitian metrics on an ample line bundle $L$ over
a Kähler manifold $(X,\omega_{0}).$ In fact, a more precise inequality
then the one conjectured by Fang will be obtained (Corollary \ref{cor:moser})
which implies both Fang's conjecture \emph{and} the Moser-Trudinger-Onofri
inequality (and hence the extremal properties of $\det\Delta_{g},$
as well). The inequality obtained for $L=\mathcal{O}(m),$ i.e. the
degree $m$ line bundle on $S^{2}$, is equivalent to the following
upper bound, \begin{equation}
\log(\frac{\det\Delta_{\bar{\partial_{u}}}}{\det\Delta_{\bar{\partial}_{0}}})\leq-\frac{1}{2}(\frac{1}{(m+2)})\int du\wedge d^{c}u(\leq0),\label{eq:sharp ineq for det in intr}\end{equation}
where we have expressed a metric $h$ on $\mathcal{O}(m)$ as $h=e^{-u}h_{0}^{\otimes m}$
in terms of the Fubini-Study metric on , which clearly implies Fang's
conjecture above. The extremals in the first inequality above will
also be characterized. 

Generalizations in two different directions will also be considered.
In the case of a higher dimensional Fano manifold $X$ (i.e $X$ admitts
a metric of positive Ricci curvature) upper bounds on the \emph{analytic
torsions} associated to powers of the anti-canonica line bundle $-K_{X}$
will be given. Compared to the higher dimensional version of the conjecture
of Gillet-Soulé this amounts to bounds on \emph{alternating} producs
of determinants of Dolbeault Laplacians in\emph{ all} degrees $(0,q)$.
Interestingly, it will appear that while it is crucial to assume the
positivity of the curvature of the metrics on $L$ in higher dimensions
this is not so when $X$ is a complex curve. In the case of complex
curves of higher genus we will prove uniqueness for mean field type
equations with one vortex on complex curves, partly confirming a conjecture
of Tarantello \cite{ta0,ta} in the context of self-dual Chern-Simons-Higgs
equations. 

The functional $\mathcal{F}_{\omega_{0}}$ referred to above is an
{}``adjoint version'', introduced by Berndtsson, of Donaldson's
(normalized) $L$-functional and generalizes the Ding-Tian functional
whose critical points are Kähler-Einstien metrics. In particular,
new proofs of some results in Kähler geometry are also obtained, including
a lower bound on Mabuchi's $K-$energy and the uniqueness result for
Kähler-Einstein metrics on Fano manifolds of Bando-Mabuchi (see section
\ref{sub:Further-relations-to} for precise references). 

Before turning to the precise statement of the main result we will
first introduce the general setup.

\subsection{Setup\label{sub:Setup} and statement of the main results}

Let $L\rightarrow X$ be a holomorphic line bundle over a compact
complex manifold $X$ of complex dimension $n.$ Denote by $\mbox{Aut}_{0}(X,L)$
the group of automorphism of $(X,L)$ in the connected component of
the identity, modulo the elements covering the identity on $X.$ The
line bundle $L$ will be assumed \emph{ample,} i.e. there exists a
Kähler form $\omega_{0}$ in the first Chern class $c_{1}(L)$ and
a {}``weight'' $\psi_{0}$ on $L$ such that $\omega_{0}$ is the
normalized curvature $(1,1)-$form of the hermitian metric on $L$
locally represented as $h_{0}=e^{-\psi_{0}}.$ In this notation, the
space of all positively curved smooth hermitian metrics on $L$ may
be identified with the open convex subset \[
\mathcal{H}_{\omega_{0}}:=\{u:\,\omega_{u}:=dd^{c}u+\omega_{0}>0\}\]
 of \emph{$\mathcal{C}^{\infty}(X),$} where $d^{c}:=i(-\partial+\overline{\partial})/4\pi,$
so that $dd^{c}=\frac{i}{2\pi}\partial\overline{\partial}.$ Hence,
$\omega_{u}$ is the normalized curvature form of the metric $e^{-u}h_{0}$
on $L$ representing the first Chern class of $L$ in $H^{2}(X,\Z).$
The natural multiplicative action of $\R^{*}$ on metrics on $L$
hence corresponds to an \emph{additive} action of $\R$ on $\mathcal{H}_{\omega_{0}},$
that we will sometimes refer to as {}``scaling''. Note that the
natural action of $\mbox{Aut}_{0}(X,L)$ on the space of all metrics
on $L$ corresponds to the action $(u,F)\mapsto v:=F^{*}(\psi_{0}+u)-\psi_{0}$
so that, in particular, $\omega_{v}=F^{*}\omega_{u}.$ Occasionally,
we will also work with the closure $\overline{\mathcal{H}}_{\omega_{0}}$
of $\mathcal{H}_{\omega_{0}}$ in $L^{1}(X,\omega_{0}),$ coinciding
with the space of all $\omega_{0}-$plurisubharmonic functions on
$X,$ i.e. the space of all upper semi-continuous functions $u$ which
are absolutely integrable and such that $\omega_{u}\geq0$ as a $(1,1)-$current. 

We equip the $N-$dimensional complex vector space $H^{0}(X,L+K_{X})$
of all holomorphic sections of the adjoint bundle $L+K_{X}$ where
$K_{X}$ is the canonical line bundle on $X,$ with the Hermitian
product induced by $\psi_{0},$ i.e. \[
\left\langle s,s\right\rangle _{\psi_{0}}:=i^{n^{2}}\int_{X}s\wedge\bar{s}e^{-\psi_{0}},\]
identifying $s$ with a holomorphic $n-$form with values in $L.$
We will use additive notation for tensor products of line bundles. 

Next, we will introduce the two functionals on $\mathcal{H}_{\omega_{0}}$
which will play a leading role in the following. First, consider the
following energy functional \begin{equation}
\mathcal{E}_{\omega_{0}}(u):=\frac{1}{(n+1)!V}\sum_{i=1}^{n}\int_{X}u(dd^{c}u+\omega_{0})^{j}\wedge(\omega_{0})^{n-j},\label{eq:def of e intro}\end{equation}
where $V:=\mbox{Vol}(\omega_{0})$ is the volume of $L,$ which seems
to first have appeared in the work of Mabuchi \cite{m2} and Aubin
\cite{au} in Kähler geometry ($\mathcal{E}_{\omega_{0}}=-F_{\omega_{0}}^{0}$
in the notation of \cite{ti}). It also appears in Arithmetic (Arakelov)
geometry as the top degree component of the secondary Bott-Chern class
of $L$ attached to the Chern character. 

The second functional $\mathcal{L}_{\omega_{0}}$ may be geometrically
defined as $\frac{1}{N}$ times the logarithm of the quotient of the
volumes of the unit-balls in $H^{0}(X,L+K_{X})$ defined by the Hermitian
products induced by the metrics $\psi_{0}$ and $\psi_{0}+u$ \cite{b-b}.
Concretely, this means that \begin{equation}
\mathcal{L}_{\omega_{0}}(u):=-\frac{1}{N}\log\mbox{det}(\left\langle s_{i},s_{j}\right\rangle _{\psi_{0}+u}),\label{eq:def of l intro}\end{equation}
where $1\leq i,j\leq N$ and $s_{i}$ is any given base in $H^{0}(X,L+K_{X})$
which is orthogonal wrt $\left\langle \cdot_{i},\cdot\right\rangle _{\psi_{0}}.$
The functional $\mathcal{L}_{\omega_{0}}(u)$ may also be invariantly
expressed as a \emph{Toeplitz determinant:}

\begin{equation}
\mathcal{L}_{\omega_{0}}(u):=-\frac{1}{N}\log\mbox{det}(T[e^{-u}]),\label{eq:l as toeplitz det intro}\end{equation}
 where $T[e^{-u}]$ is the \emph{Toeplitz operator with symbol $e^{-u}$}
(compare formula \ref{eq:toeplitz as k} in the appendix).  If $N=0$
we let $\mathcal{L}_{\omega_{0}}(u):=-\infty.$ The normalizations
are made so that the functional

\[
\mathcal{F}_{\omega_{0}}:=\mathcal{E}_{\omega_{0}}-\mathcal{L}_{\omega_{0}}\]
is invariant under the action of $\R$ and hence descends to a functional
on the space of all Kähler metrics in $c_{1}(L).$ An element $u$
in $\mathcal{H}_{\omega_{0}}$ will be said to be \emph{critical}
(wrt $L+K_{X})$ if it is a critical point of the functional $\mathcal{F}_{\omega_{0}}$
on $\mathcal{H}_{\omega_{0}},$ i.e. if $u$ is a smooth solution
in $\mathcal{H}_{\omega_{0}}$ of the Euler-Lagrange equations $(d\mathcal{F}_{\omega_{0}})_{u}=0.$
These equations may be written as the highly non-linear Monge-Ampère
equation: \begin{equation}
\frac{1}{Vn!}(dd^{c}u+\omega_{0})^{n}=\beta(u),\label{eq:e-l for f}\end{equation}
 where $\beta(u)$ is the Bergman measure associated to $u$ (formula
\ref{eq:bergman meas} below). This latter measure depends on $u$
in a \emph{non-local} manner and is strictly positive precisely when
$L+K_{X}$ is \emph{globally generated,} i.e. when there, given any
point $x$ in $X,$ exists an element $s$ in $H^{0}(X,L+K_{X})$
such that $s(x)\neq0.$ For example, since $L$ is ample, this condition
holds when $L$ is replaced by $kL$ for $k$ sufficiently large. 

By definition, a critical point $u$ is a priori only a \emph{local}
extremum of $\mathcal{\mathcal{F}}_{\omega_{0}}.$ But the next theorem
relates \emph{global} maximizers of $\mathcal{\mathcal{F}}_{\omega_{0}}$
and its critical points:
\begin{thm}
\emph{\label{thm:main}Let $L$ be an ample line bundle over a compact
complex manifold $X.$ Then the absolute maximum of the functional}
\textup{$\mathcal{\mathcal{F}}_{\omega_{0}}$ on $\mathcal{H}_{\omega_{0}},$
defined above, is attained at any critical point $u.$ Moreover if}\emph{
the adjoint line bundle $L+K_{X}$ is globally generated}\textup{
any smooth maximizer of $\mathcal{\mathcal{F}}_{\omega_{0}}$ on}
$\overline{\mathcal{H}}_{\omega_{0}}$ \textup{is unique (up to addition
of constants) modulo the action of} $\mbox{Aut}_{0}(X,L).$ \emph{In
particular, such a maximizer is critical.} 
\end{thm}
In the case when the ample line bundle $L=-K_{X},$ so that $X$ is
a Fano manifold, the space $H^{0}(X,L+K_{X})$ is one-dimensional
and hence $\mathcal{L}_{\omega_{0}}(u)=-\frac{1}{N}\log\int e^{-(u+\psi_{0})}.$
Then it is well-known that any critical point may be identified with
a Kähler-Einstein metric on $X$, i.e. the corresponding Kähler metric
satisfies \[
\mbox{Ric }\omega_{u}=\omega_{u}\]
 It should be emphasized that the \emph{existence} of critical points
of $\mathcal{\mathcal{F}}_{\omega_{0}}$ is a difficult issue. For
example, in the case $L=-K_{X}$ fundamental conjectures of Yau, Tian,
Donaldson \cite{ti,do1,th} relate the existence problem to a notion
of algebro-geometric stability. In section \ref{sub:Conjectures}
some conjectures concerning the case of a general line bundle $L$
are proposed. As explained in remark \ref{rem:existence} there is
a weak solution (in the sense of pluripotential theory) to the critical
point equation if the functional $\mathcal{F}_{\omega_{0}}$ is \emph{proper}
of \emph{coercive} in a suitable sense. In the case of $L=-K_{X}$
this amounts to saying that $X$ is \emph{analytically stable} in
the sense of Tian \cite{ti}. Moreover, by Proposition \ref{pro:non-deg max}
and the subsequent remark the property of admitting a critical point
is an open condition in the moduli space of polarized manifolds $(X,L)$
as long as $\mbox{Aut}_{0}(X,L)$ is trivial.

However, in the presense of large symmetry groups existence of critical
points can be readily established. Hence, we next assume that $(X,L)$
is \emph{$K-$homogenous,} i.e. that $X$ admits a transitive action
by a compact semi-simple Lie group $K,$ whose action on $X$ lifts
to $L.$ We will then take $\omega_{0}$ as the unique Kähler form
in $c_{1}(L)$ which is invariant under the action of $K$ on $X.$ 
\begin{cor}
\label{cor:homeg}Let $L\rightarrow X$ be a $K-$homogenous ample
holomorphic line bundle over a compact complex manifold $X$ and denote
by $\omega_{0}$ be the unique $K-$invariant Kähler metric in $c_{1}(L).$
Then, for any function $u$ in $\mathcal{H}_{\omega_{0}}$ \[
-\mathcal{L}_{\omega_{0}}(u)\leq-\mathcal{E}_{\omega_{0}}(u)\]
 with equality iff the function $u$ is constant, modulo the action
of $\mbox{Aut}_{0}(X,L).$ 
\end{cor}
Surprisingly, specializing to the case when $X$ is a complex curve
(i.e. $n=1)$ allows one to take $u$ as \emph{any} smooth function
(which is \emph{not }true in higher dimensions, as shown in \cite{rub0}
in the case when $L=-K_{X};$ see remark \ref{rem:j-f}). More generally,
we can then take $u$ to be in the Sobolev space $W^{2,1}(X)$ of
all functions $u$ on $X$ such that $u$ and its differential $du$
are square integrable. We will first consider the homegenous case
in one dimension, i.e. $X=\P^{1},$ the complex projective line (i.e.
topologically $X=S^{2},$ the two-sphere) and hence $K=SU(2).$ Identifying
$S^{2}$ with the one-point compactification of $\C_{z}$ (compare
section \ref{sub:Explicit-expression}) we then get the following
corollary 
\begin{cor}
\label{cor:moser}Let $u$ be a function in the Sobolev space $W^{2,1}(S^{2})$
on the two-sphere $S^{2}$ and denote by $\omega_{0}$ the volume
form corresponding to the standard invariant metric on $S^{2}$ of
unit-area. Then \[
\log\mbox{det}(c_{ij}\int_{\C}\frac{z^{i}\bar{z}^{j}}{(1+z\bar{z})^{m}}e^{-u}\omega_{0})\leq-(m+1)\int_{S^{2}}u\omega_{0}+(\frac{m+1}{m+2})\frac{1}{2}\int_{S^{2}}du\wedge du^{c}\]
where $1/c_{ij}=(m+1)\binom{m}{i}\binom{m}{j}$ with equality iff
there exists a Möbius transformation $F$ of $S^{2}$ such that $\omega_{u}=F^{*}\omega_{0}.$ 
\end{cor}
The case when $m=0$ gives a new proof of the Moser-Trudinger-Onofri
inequality \ref{eq:m-t-o intro}. The reduction of the proof of Corollary
\ref{cor:moser} to Corollary \ref{cor:homeg} is based on properties
of a projection operator $P_{\omega_{0}}$ (formula \ref{eq:proj oper intro}).
As pointed out by Fang \cite{f} the determinant in the previous corollary
may also be expressed as the integral of $\sum_{i}e^{-u(x_{i})}$
over the $N$ fold product of $S^{2}$\emph{ (with $N=m+1)$} of the
$SU(2)-$invariant probability measure with density \[
\rho_{N}(x_{1},...,x_{N}):=\Pi_{1\leq i<j\leq N}\left\Vert x_{i}-x_{j}\right\Vert ^{2}/Z_{N}\]
where $1/Z_{N}=N^{N}\binom{N-1}{0}...\binom{N-1}{N-1}/N!,$ written
in terms of the ambient Euclidian norm in $\R^{3}.$ In the physics
litterature this integral appears as the free energy for the Gibbs
ensemble of a \emph{Coulomb gas} of unit-charge particles (i.e one
component plasma) confined to the sphere \cite{ca} in a neutralizing
uniform background. In the case of a general line bundle $\mathcal{L}_{\omega_{0}}$
may be expressed in terms of a determinantal random point process
(see section \ref{sub:Bergman-kernels}) 

In the non-homegenous one dimensional case, i.e. higher genus curves
the assumption in Theorem \ref{thm:main} about global generation
of $L+K_{X}$ turns out to be superflous. A simple application of
the Riemann-Roch theorem shows that $L+K_{X}$ is not globally generated
precisely when $L:=L_{p}$ is an effective degree one line bundle,
i.e. $L$ has a holomorphic section $s$ whose zero divisor is an
irreducible point $p$ in $X.$ 
\begin{thm}
\label{thm:main ine one dim}Let $L$ be a degree one line bundle
over a complex curve $X$ of genus $g\geq1.$ Then the\textup{\emph{
}}the functional\emph{ }\textup{\emph{$\mathcal{\mathcal{F}}_{\omega_{0}}$
on $\mathcal{H}_{\omega_{0}},$ admits a unique (modulo scaling) critical
point $u$ on $\mathcal{C}^{\infty}(X).$ Moreover $u$ is a maximizer
of $\mathcal{\mathcal{F}}_{\omega_{0}}$ on $\mathcal{C}^{\infty}(X).$ }}
\end{thm}
The \emph{existence }of a critical point in the previous theorem follows
from a variational argument, using the general Moser-Trudinger inequality
of Fontana \cite{cher}. This variational approach is closely related
to the one used by Troyanov \cite{tr} in his study of constant curvature
metrics on Riemann surfaces with conical singularities.

\subsection{\label{sub:Application-to-determinants}Applications to twisted analytic
torsion}

Given any line bundle $L$ over a complex curve $(X,\omega_{0})$
equipped with a Kähler metric $\omega_{0}$ of constant curvature
we will denote by $\det\Delta_{\bar{\partial}}$ the functional on
the space of all smooth Hermitian metrics on $L$ which to an Hermitian
metric associates the zeta-regularized product of the strictly positive
eigenvalues of the corresponding $\partial-$Laplacian (see section
\ref{sub:Analytic-torsion}). Combining Corollary \ref{cor:moser}
above with the anomaly formula of Bismut-Gillet-Soulé \cite{b-g-s}
now implies the following positive solution of Fang's conjecture
\begin{cor}
\label{cor:max of det}For any line bundle over $L$ over $\P^{1}$
equipped with its invariant metric the functional $\det\Delta_{\bar{\partial}}$
has a unique maximizer (modulo scaling), namely the Hermitian metric
on $L$ whose curvature form $\omega$ has constant curvature. 
\end{cor}
In fact, the proof of the previous corollary will give the stronger
statement that the inequality \ref{eq:sharp ineq for det in intr}
for $\mbox{\ensuremath{\det}}\Delta_{\bar{\partial}_{u}}$ stated
in the introduction holds and that this latter inequality is equivalent
to Corollary \ref{cor:moser}. Note that a direct consequence of the
previous corollary is the following reponse to a variant of Kac's
classical question {}``Can one hear the shape of a drum?\textquotedbl{}
\cite{ka}: if the $\bar{\partial}-$Laplacian on \emph{some} power
$\mathcal{O}(m)$ induced by a smooth metric on $\mathcal{O}(1)\rightarrow\P^{1}$
has the same spectrum (including multiplicities) as the $\bar{\partial}-$Laplacian
induced by the standard $SU(2)$ invariant metric then the two metrics
coincide up to scaling. 

In the case when $X$ is a genus one curve, we obtain the following
\begin{cor}
\label{cor:max of det on torus}Let $L$ be degree one line bundle
over a complex curve $X$ of genus one and write $L=L_{p}$ for a
unique point $p$ on $X.$ Then the functional $\det\Delta_{\bar{\partial}}$
has a unique maximizer (modulo scaling), namely the Hermitian metric
on $L_{p}$ whose curvature form is the constant curvature metric
$\omega_{p}$ with a conical singularity at $p$ (see below: $\omega_{p}$
satisfies equation $(ii)'$ in \ref{eq:vortex ex curv form} for $t=1$).
\end{cor}
A natural higher dimensional generalization of the regularized determinant
associated to a line bundle $L$ over a Kähler manifold $(X,\omega_{0})$
(we will assume that $\omega_{0}\in c_{1}(L))$ is the \emph{analytic
torsion} (see section \ref{sub:Analytic-torsion}). The following
theorem can be seen as a generalization of a weak form of Fang's result
on $S^{2}$ to certain Fano manifolds. It is formulated in terms of
the follolwing invariant $R(X)$ of a Fano manifold $X:$ \begin{equation}
R(X)=\sup_{\omega_{t}}\{t:\,\,\,\mbox{Ric }\mbox{\ensuremath{\omega}}_{t}>t\omega_{t}\}>0\label{eq:def of r}\end{equation}
with $\omega_{t}$ ranging over all Kähler metrics in $c_{1}(-K_{X})$
(this invariant was studied very recently by Székelyhidi,\cite{sz};
the strict positivity follows from \cite{au}). 
\begin{thm}
\label{thm:fano torsion}Let $X$ be a Fano manifold with $R(X)>1-1/n^{2}.$
Then, for $k$ sufficently large, the analytic torsion associated
to $-kK_{X}$ is bounded from above when considered as a functional
on the space of all positively curved metrics on $-kK_{X}.$ In particular,
the boundeness holds for any Fano surface. Moreover, when $(X,\omega_{0})$
is equal to $\P^{2}$ equipped with the standard invariant metric
$\omega_{0}$ the upper bound holds for any positive $k$ and the
maximum is achieved precisely on the metrics of constant curvature
$\omega_{0}$, i.e. precisely on the Fubini-Study metric on $-kK_{X}$
(up to scaling). 
\end{thm}
The fact that the condition on $R(X)$ above is satisfied on Fano
surfaces (i.e. on the Del Pezzo surfaces) follows from \cite{t-y,li,sz}.
If $X$ admits a Kähler-Einstein metric then $R(X)=1$ (but the converse
is not true) \cite{sz}. It would be interesting to know to what extent
the condition on $R(X)$ above is satisfied in higher dimension? One
motivation for the previous theorem comes from the following observation
(see remark \ref{rem:torsion entropy}): after a suitable scaling
the analytic torsion $\mathcal{T}_{k}(u)$ for a metric $h_{ku}$
on $-kK_{X}$ (expressed in terms of $u\in\mathcal{H}_{\omega_{0}}$
as $h_{ku}=e^{-ku}h_{0}^{\otimes k})$ has the following asymptotics
due to Bismut.Vasserot \cite{b-v}: \[
\lim_{k\rightarrow\infty}\mathcal{T}_{k}(u)=-\mathcal{S}(\mu_{u},\mu_{0})\]
 i.e. minus the \emph{relative entropy} of the probability measure
$\mu_{u}:=\omega_{u}^{n}/Vn!$ with respect to the measure $\mu_{0}.$
It is a well-known fact that $\mathcal{S}(\mu,\nu)\geq0$ with equality
if and only if the two measures coincide. Hence, the previous theorem
can be seen as an analogue of this latter fact for \emph{finite} $k.$
We will also show that, in the case of $(\P^{n},\omega_{0})$ for
$n\leq2$ the analytic torsion functionals $\mathcal{T}_{k}$ are
geodesically\emph{ concave} on $\mathcal{H}_{\omega_{0}}$ equipped
with its symmetric space metric. This should be compared with a recent
result of von Renesse-Sturm \cite{r-s} saying that on a Riemannian
manifold $(X,g_{0})$ of unit-volume the relative entropy $\mathcal{S}(\mu,dVol_{g_{0}})$
is convex with respect to the $L^{2}-$Wasserstein metric on the space
of all probability measures $\mu$ on $X$ precisely when $g_{0}$
has semi-positive Ricci curvature.

\subsection{\label{sub:Application-to-vortices}Application to vortices and metrics
with conical singularities}

Let $(X,\omega_{0})$ be a Riemann surface with a Kähler (i.e. area)
form $\omega_{0}$ of unit area. Given a parameter $\mu\in\R_{+}$
and a fixed point $p$ the corresponding (normalized) \emph{mean field
equation} with a single vortex at $p$ is defined as \begin{equation}
(i)\,\Delta u=\eta(e^{g_{p}-u}-1)\label{eq:intro mean field for u}\end{equation}
 where $\Delta$ is the Laplacian with respect $\omega_{0}$ and $g_{p}$
is the correspondong Green function with a pole at $p$ (note that
any solution satisfies $\int_{X}e^{g_{p}-u}=1/\eta)$ A similar mean
field type equation was also studied very recently in \cite{lw}:
\begin{equation}
(ii)\,\Delta w+\eta e^{w}=\eta\delta_{p}\label{eq:intro mean field for v}\end{equation}
which is equivalent to the previous one when $\eta=4\pi$ (just set
$w=g_{p}-u).$ These two equations appear in the mean field limit
of a statistical mechanics model for vortices in fluids \cite{k,clmp}.
For a very recent account of mean field equations and their relations
to self-dual \emph{Chern-Simons-Higgs equations} see for example \cite{ta}
(see also section \ref{sub:Relation-to-random}). Equivalently, setting
$t:=\eta/4\pi$ the pseudo-metrics $\Omega_{t}:=\omega_{u/t}$ and
$\omega_{t}:=\omega_{g_{p}-w/t}$ satisfy \begin{equation}
(i')\,\,\mbox{Ric\,\ensuremath{\Omega_{t}:=}}t\Omega_{t}-[p]+(1-t)\omega_{0},\,\,\,(ii')\,\,\mbox{Ric\,\ensuremath{\omega_{t}:=}}t(\omega_{t}-[p])\label{eq:vortex ex curv form}\end{equation}
For $t=1$ both equations hence describe a metric $\omega_{t}$ of
unit area with unit positive Gauss curvature and a conical singularity
at $p$ of angle $4\pi.$ The following theorem should be seen in
the light of the well-known fact that conical constant curvature metrics
may be non-unique for angles stricly larger than $2\pi$ (for the
case of the two-sphere see \cite{tr2,er}). 
\begin{thm}
\label{cor:mean field}Let $X$ be a Riemann surface of genus one
with a marked point $p$ and let $\omega_{0}$ be the standard invariant
metric on $X$ (i.e. $(X,p,\omega_{0})$ may be identified with $(\C/(\Z+\tau\Z),0,\frac{i}{2}dz\wedge d\bar{z})$
for $\tau$ in the upper half plane). Then the mean field equations
\ref{eq:intro mean field for u} and \ref{eq:intro mean field for v}
both admit a unique solution when $\eta\in]0,4\pi+\epsilon]$ for
some $\epsilon>0.$ 
\end{thm}
The \emph{existence }of a solution for $\eta\in]0,8\pi[$ is well-known
and follows from Moser-Trudinger inequalities and basic variational
arguments as above. Uniqueness in the case when $\eta\in]0,8\pi[$
was conjectured by Tarantello \cite{ta0} (see also the discussion
on page 158 in \cite{ta}). The uniqueness for $\mu=4\pi$ was established
very recently in \cite{lw} (Theorem 3.2) using function theory on
the elliptic curve $\C/(\Z+\tau\Z),$ where the uniqueness of solutions
for equation \ref{eq:intro mean field for v} was conjectured for
$\eta\in[4\pi,8\pi[$ as well. Under the latter assumption on $\eta$
the authors proved the uniqueness of solutions invariant under the
natural $\Z_{2}$ action $z\mapsto-z$ on $\C/(\Z+\tau\Z).$ Remarkably,
it was also shown in \cite{lw} that uniqueness of equation \ref{eq:intro mean field for v}
fails for $\eta=8\pi$ as soon as there is some solution (which was
shown to happen when $\tau=e^{i2\pi/3}).$ For uniqueness results
in the different setting of plane domains see \cite{ba,clmp}. As
a byproduct of the present proof we will also obtain uniqueness of
the \emph{linearized} equations (in the case of $\eta=4\pi$ this
{}``linear stability'' was proved in \cite{lw} under the assumption
of $\Z_{2}$ symmetry).

It would be interesting to know whether the present methods could
be complemented in order to deal with the uniqueness problem when
$\eta\in]4\pi,8\pi[$ and the case when $X$ is a Riemann surface
of genus at least two. In the latter case the present proof gives
a generalization of the previous corollary to a variant of the mean
field equation where the factor $e^{g_{p}-u}$ is multiplied with
a certain smooth and strictly positive function $f_{u}$ (see equation
\ref{eq:critical point eq in one dim}). For a fixed function $f$
(often independent of $u)$ such generalized mean field equations
have been studied extensively in the litterature (see \cite{ta} and
references therein).

Finally, it should be emphasized that it is the positive sign of the
parameter $\mu$ that is the source of the delicate analytical properties
of the equation \ref{eq:intro mean field for u}. Indeed, when $\eta<0$
it follows readily from the maximum principle that \ref{eq:intro mean field for u}
has a unique solution. In this case the equations are equivalent to
the abelian \emph{Yang-Mills-Higgs equations} (see section \ref{sub:Relation-to-random}).

\subsection{Applications to Kähler geometry}

These applications are formulated and proved in section \ref{sec:Convergence-towards-Mabuchi's}.

\subsection{\label{sub:Further-relations-to}Further relations to previous results}

The relation between the Moser-Trudinger type inequalities and Kähler-Einstein
metrics of positive curvature in higher dimensions seems to first
have been suggested by Aubin \cite{au} who also proposed a conjecture
that we will be proved in section \ref{sub:Proof-of-a conj aub}.
In the case when $L=-K_{X}$ the first statement of Theorem \ref{thm:main}
is a result of Ding-Tian\cite{d-t} and the {}``uniqueness'' of
critical points (i.e. Kähler-Einstein metrics in this case) was proved
earlier by Bando-Mabuchi \cite{b-m}. See \cite{bbgz} for a generalization
of this latter result to functions of {}``finite energy'', in the
case when $\mbox{Aut}_{0}(X,L)$ is trivial (compare remark \ref{rem:finite energy}). 

The extremal property of the critical points in Theorem \ref{thm:main}
can also be seen as an analog of a result of Donaldson (Theorem 2
in \cite{don1}) who furthermore assumed that $\mbox{Aut}_{0}(X,L)$
is discrete. In this latter setting the role of the space $H^{0}(X,L+K_{X})$
is played by $H^{0}(X,L)$ equipped with the scalar products induced
by the weight $\psi_{0}+u$ and the integration measure $(\omega_{u})^{n}/n!$
Note however that in Donaldson's setting the functional $\tilde{\mathcal{F}}$
corresponding to $\mathcal{F}_{\omega_{0}}$ is \emph{minimized} on
its critical points (compare section \ref{sub:Comparison-with-Donaldson's}
and the discussion in section 5 in \cite{bern2}). In fact, when $X$
is an arithmetic variety the functional $\tilde{\mathcal{F}}$ seems
to first have appeared in the work of Bost \cite{bo} and Zhang \cite{zh}.
In the terminilogy of \cite{don1} these latter critical points correpond
to \emph{balanced metrics.} Donaldson used his result, combined with
the deep convergence results in \cite{do1} for balanced metrics,
in the limit when $L$ is replaced by a large tensor power, to prove
a lower bound on Mabuchi's K-energy functional (under the assumption
that $\mbox{Aut}_{0}(X,L)$ be trivial). It will be shown in section
\ref{sec:Convergence-towards-Mabuchi's} how to deduce this latter
result more directly from Theorem \ref{thm:main} above. After the
first version of the present paper (the preprint ref) appeared another
proof by Li ref of the lower bound on Mabuchi's K-energy functional
appeared which does not use that $\mbox{Aut}_{0}(X,L)$ be trivial
either.

It should also be pointed out that the inequality proved by Donaldson
corresponds to a \emph{lower} bound on $\mathcal{\mathcal{F}}_{\omega_{0}}(u)$
in the present setting, which however will depend on $u$ through
its volume form $(\omega_{u})^{n}/n!$ (see the end of section\ref{sub:Comparison-with-Donaldson's}
). 

Note also that Rubinstein \cite{rub,rub0} recently gave a different
complex geometric proof of the Moser-Trudinger-Onofri inequality \ref{eq:m-t-o intro}
using the inverse Ricci operator and its relation to various energy
functionals in Kähler geometry. See also Müller-Wendland \cite{m}
for a proof of the result on extremals of determinants of the scalar
Laplacian using the Ricci flow. It should also be pointed out that
extremals related to twisted analytic torsion on Kähler manifolds
were also studied in \cite{m}. However the torsion was twisted by
the Quillen metric on a certain virtual vector bundle earlier introduced
by Donaldson in the seminal work \cite{doo}.

\subsection{Outiline of the proofs of the main results.}

\subsubsection*{Theorem \textup{\ref{thm:main}}}

The starting point is the recent work \cite{bern2} of Berndtsson
which shows that $\mathcal{F}_{\omega_{0}}$ is concave along geodesics
with respect to the Riemann metric on the space $\mathcal{H}_{\omega_{0}}$
introduced by Mabuchi \cite{m}. Moreover, the concavity is \emph{strict}
modulo the action of $\mbox{Aut \ensuremath{(X,L)_{0}}}.$ Hence,
if $\mathcal{H}_{\omega_{0}}$ were geodesically convex (i.e. if any
two points in $\mathcal{H}_{\omega_{0}}$ could be connected by a
smooth geodesic) Theorem \ref{thm:main} would follow immediately.
However, since there are no existence results for such geodesic segments
(apart from the toric setting) one of the the key points of the present
paper is to show how to systematically use the global pluripotential
theory developed in \cite{b-b,b-d} to get around this difficulty.
This is achieved by working with continuous {}``geodesics'' in the
$L^{1}-$closure $\overline{\mathcal{H}}_{\omega_{0}}$ of $\mathcal{H}_{\omega_{0}}.$
A priori such a geodesic may hence touch the ''boundary'' $\overline{\mathcal{H}}_{\omega_{0}}-\mathcal{H}_{\omega_{0}},$
but we will exlude this scenario by showing that $\mathcal{\mathcal{F}}_{\omega_{0}}$
has no maximizers on the {}``boundary'' of $\mathcal{H}_{\omega_{0}}$
(if $K_{X}+L$ is globally generated). Following \cite{b-b,bbgz}
this crucial latter fact is shown by extending $\mathcal{\mathcal{F}}_{\omega_{0}}$
to a (Gâteaux) differentiable function on all of $C^{0}(X),$ by replacing
$\mathcal{E}_{\omega_{0}}$ with the composed map $\mathcal{E}_{\omega_{0}}\circ P_{\omega_{0}},$
where $P_{\omega_{0}}$ is the following (non-linear) projection operator
from $C^{0}(X)$ onto $\mathcal{C}^{0}(X)\cap\overline{\mathcal{H}}_{\omega_{0}}:$
\begin{equation}
P_{\omega_{0}}[u](x)=\sup\left\{ v(x):\, v\in\mathcal{H}_{\omega_{0}},\,\, v\leq u\right\} \label{eq:proj oper intro}\end{equation}
Somewhat remarkbly, this latter projection operator is also used to
show that that any critical point of $\mathcal{F}_{\omega_{0}}$ is
automatically a global maximizer on all of $C^{\infty}(X)$ when $n=1.$
Moreover, it is shown that the assumption that $K_{X}+L$ be globally
generated is not needed when $n=1$ by refining the arguments in \cite{bern2}.

\subsubsection*{The applications to analytic torsion }

These applications are based on the fact that the analytic torsion
functional may be expressed as $\mathcal{F}_{\omega_{0}}$ plus an
explicit curvature term (which vanishes when $X$ is a torus). Moreover,
the main contribution in this term (when $L$ is replaced by a large
tensor power) comes from the Ricci curvature for $\omega_{0}$ and
this leading term appears with a useful negative sign when the Ricci
curvature is positive.

\subsubsection*{The applications to vortices }

These are obtained by studying a deformation of $\mathcal{\mathcal{F}}_{\omega_{0}},$
depending on $\mu,$ which preserves the concavity when $\mu\in]0,4\pi].$
The main new difficulty comes from the fact that one has to work with
forms $\omega_{u}$ degenerating at the fixed point $p.$

\subsection*{Acknowledgement}

It is a pleasure to thank Bo Berntdsson for illuminating discussions
on the topic of the present paper, in particular in connection to
\cite{bern2}. It is an equal pleasure to thank Sébastien Boucksom,
Vincent Guedj and Ahmed Zeriahi for discussions and stimulation coming
from the colaboration \cite{bbgz}. The author is also grateful to
Yanir Rubinstein for helpful comments on a preliminary version of
this paper, as well as Nikolas Akerblom and Gunther Cornelissen for
drawing my attention to the references \cite{ta} and \cite{lw}.

\subsection*{Organization }

In section 2 preliminaires for the proofs of the main results appearing
in section 3 are given. The proof of the uniqueness statement in the
main theorem relies on higher order regularity for {}``geodesics''
defined by inhomogenous Monge-Ampère equations. An alternative proof
based on considerably more elementary regularity results is given
in section \ref{sub:Alternative-proof-of}, as well extensions to
degenerate boundary data in the one-dimensional case. Finally, in
section 4 the limit when the line bundle $L$ is replaced by a large
tensor power is studied and a new proof of the lower bound on Mabuchi's
$K-$energy for a polarized projective manifold. The proof of Theorem
\ref{thm:fano torsion} concerning upper bounds on twisted analytic
torsions on Fano manifolds and a proof of a conjecture of Aubin are
also given and relations to Donaldson's work are discussed. Some conjectures
are also proposed. In the appendix some formulas involving Bergman
kernels are recalled and a {}``Bergman kernel proof'' of Theorem
\ref{thm:(Berndtsson)-Let-} is given, as well as generalizations
to a degenerate setting in the one-dimensional case.

\section{Preliminaries: Geodesics and functionals }

\subsection{\label{sub:Geodesics-and-psh}Geodesics and psh paths}

The infinite dimensional space $\mathcal{H}_{\omega}$ inherits an
\emph{affine} Riemannian structure from its natural imbedding as on
open set in $\mathcal{C}^{\infty}(X).$ Mabuchi, Semmes and Donaldson
(see \cite{ch} and references therein) introduced another Riemannian
structure on $\mathcal{H}_{\omega}$ (modolo the constants) defined
in the following way. Identifying the tangent space of $\mathcal{H}_{\omega}$
at the point $u$ with $\mathcal{C}^{\infty}(X)$ the squared norm
of a tangent vector $v$ at the point $u$ is defined as \[
\int_{X}v^{2}(\omega_{u})^{n}/n!.\]
In fact, this metric is induced by a natural symmetric space structure
on $\mathcal{H}_{\omega}.$ However, the \emph{existence} of a geodesic
$u_{t}$ in $\mathcal{H}_{\omega}$ connecting any given points $u_{0}$
and $u_{1}$ is an open and perhaps even dubious problem. There are
two problems: it is not known if $i)$ $u_{t}$ smooth, $ii)$ $\omega_{u_{t}}$
is strictly positive, as a current. As is well-known such a geodesic
may, if it exists, be obtained as the solution of a homogenous Monge-Ampère
equation (see below). In the following we will simply take this characterization
as the \emph{definition} of a geodesic. It will also be important
to consider the larger space $\overline{\mathcal{H}}_{\omega_{0}}\cap C^{0}(X),$
since a priori the path $u_{t}$ may leave $\mathcal{H}_{\omega}.$
\begin{defn}
A continuous path in $\overline{\mathcal{H}}_{\omega_{0}}\cap C^{0}(X)$
$u_{t}$ will be called a \emph{$\mathcal{C}^{0}-$geodesic} connecting
$u_{0}$ and $u_{1}$ if $U(w,x):=u_{t}(x),$ where $t=\log\left|w\right|,$
is continuous on \[
M:=\{1\leq\left|w\right|\leq e\}\times X:=A\times X\]
with $dd^{c}U+\pi_{X}^{*}\omega_{0}\geq0$ and \begin{equation}
(dd^{c}U+\pi_{X}^{*}\omega_{0})^{n+1}=0\label{eq:dirichlet pr}\end{equation}
 in the interiour of $M$ in the sense of pluripotential theory \cite{g-z,de3},
where $\pi_{X}$ denotes the projection from $M$ to $X.$ 
\end{defn}
As shown in \cite{b-d,bbgz} $U(w,x)$ exists and is uniquely defined
as the extension from $\partial M$ obtained as the upper envelope
\begin{equation}
U(w,x)=\sup\left\{ V(w,x):\, V\in\mathcal{H}_{\pi_{X}^{*}\omega_{0}}(M),\,\, V\leq U\,\mbox{on}\,\partial M\right\} ,\label{eq:envelop}\end{equation}
where $\mathcal{H}_{\pi_{X}^{*}\omega_{0}}(M)$ denotes the set of
all smooth functions $V$ on $M$ such that $dd^{c}U+\pi_{X}^{*}\omega_{0}>0.$
If $u_{t}$ is such that $dd^{c}U+\pi_{X}^{*}\omega_{0}\geq0$ then
$u_{t}$ will be called a \emph{psh path} (or a \emph{subgeodesic}).
In local computations we will often make the identification $u_{t}(x)=U(w,x)$
extending $t$ to a complex variable in a strip in the complex plane..
Then $u_{t}(x)$ is independent of the imaginary part of $t$ and
is hence \emph{convex} wrt real $t.$ In this notation a binomial
expansion in equation \ref{eq:dirichlet pr} shows the equation for
a smooth geodesic $u_{t}$ may be written as \[
\partial_{\bar{t}}\partial_{t}u_{t}\omega_{u_{t}}{}^{n}/n!-\partial_{X}\partial_{\bar{t}}u\wedge\bar{\partial}_{X}\partial_{t}u\wedge\omega_{u_{t}}{}^{n-1}/(n-1)!=0,\]
where the subscript $X$ indicate derivatives along $X$ and $\partial_{t}$
and $\partial_{\bar{t}}$ denote the holomorphic and anti-holomorphic
partial derivatives wrt $t$, i.e. \begin{equation}
c(u_{t}):=\partial_{\bar{t}}\partial_{t}u_{t}-|(\partial(\partial_{\bar{t}}u)|_{\omega_{u_{t}}}^{2}=0\label{eq:geod eq concrete}\end{equation}
In the proof of the uniqueness part of Theorem \ref{thm:main} we
will have great use for the following regularity result for geodesics
in $\overline{\mathcal{H}}_{\omega_{0}},$ shown by Chen \cite{ch}.
See also \cite{bl} for a detailed analysis of the proof and some
refinements. The proof uses the method of continuity combined with
very precise a prioiri estimates on the perturbed Monge-Ampère equations.
\begin{thm}
\label{thm:(Chen)-Assume-that}(Chen) Assume that the boundary data
in the Dirichlet problem \ref{eq:dirichlet pr} for the Monge-Ampère
operator on $M$ is smooth on $\partial M.$ Then $U\in\mathcal{C}_{\C}^{1,1}(M).$
More precisely, the mixed second order complex derivatives of $U$
are uniformly bounded, i.e. there is a positive constant $C$ such
that \[
0\leq(dd^{c}U+\pi_{X}^{*}\omega_{0})\leq C(\pi_{X}^{*}\omega_{0}+\pi_{A}^{*}\omega_{A})\]
 where $\omega_{A}$ is the Eucledian metric on $A.$
\end{thm}
In the statement above we have used the (non-standard) notation $\mathcal{C}_{\C}^{1,1}(M)$
for the set of all functions $U$ such that, locally, the current
$dd^{c}U$ has coefficents in $L^{\infty}.$ Such a $U$ is called
\emph{almost $\mathcal{C}^{1,1}$} in \cite{bl}. Note that if $U\in\mathcal{H}_{\pi_{X}^{*}\omega_{0}}(M)$
then this is equivalent to $U$ having a bounded Laplacian $\Delta_{M}U,$
where $\Delta_{M}$ is the Laplacian on $M$ wrt the Kähler metric
$\pi_{X}^{*}\omega_{0}+\pi_{A}^{*}\omega_{A}$ on $M.$ As will be
explained in section \ref{sub:Alternative-proof-of} the proof of
the uniqueness statement in Theorem \ref{thm:main} may actually be
obtained by only using the bounds on the derivatives of $u_{t}$ on
$X$ for $t$ fixed. As shown very recently in \cite{b-d} such bounds
may be obtained by working directly with the envelope \ref{eq:envelop}. 
\begin{thm}
\label{thm:berman-dem}Assume that the boundary data in the Dirichlet
problem \ref{eq:dirichlet pr} for the Monge-Ampère operator on $M$
is in $\mathcal{C}^{1,1}(\partial M).$ Then $u_{t}\in\mathcal{C}_{\C}^{1,1}(X).$
More precisely, the mixed second order complex derivatives of $u_{t}$
on $X$ are uniformly bounded, i.e. there is a positive constant $C$
such that \[
0\leq(dd^{c}u_{t}+\omega_{0})\leq C\omega_{0}\]
 on $X.$
\end{thm}
One of the virtues of this latter approach is that the proof is remarklably
simple when $X$ is homogenous.

\subsection{The functional $\mathcal{L}_{\omega_{0}}.$ }

First note that the functional $\mathcal{L}_{\omega_{0}}(u)$ defined
by formula \ref{eq:l as toeplitz det intro} is increasing on $\mathcal{C}^{0}(X),$
wrt the usual order relation. This is an immediate consequence of
the basic geometric interpretation in \cite{b-b} of $\mathcal{L}_{\omega_{0}}(u)$
as proportional to the logarithmic volume of the unit-ball in the
Hilbert space $H^{0}(X,L+K_{X})$ equipped with the Hermitian product
induced by the weight $\psi_{0}+u.$ Alternatively, this fact follows
from formula \ref{eq:deriv of l} below which shows that the differential
of the functional $\mathcal{L}_{\omega_{0}}$ on $\mathcal{C}^{0}(X)$
may be represented by the positive measure $\beta_{u}.$ Integrating
$\beta_{u}$ along a line segment in $\mathcal{C}^{0}(X)$ equipped
with its affine structure then shows that $\mathcal{L}_{\omega_{0}}(u)$
is increasing. 

The \emph{differential} of the functional $\mathcal{L}_{\omega_{0}}$
on $\mathcal{C}^{0}(X)$ is given by 

\begin{equation}
(d\mathcal{L}_{\omega_{0}})_{u}=\beta_{u},\label{eq:deriv of l}\end{equation}
in the sense that given any smooth function $v$ we have that \[
d(\mathcal{L}_{\omega_{0}}(u+tv))/dt_{t=0}=\int_{X}\beta_{u}v,\]
 where $\beta_{u}$ is the \emph{Bergman measure associated to $u.$}
This latter measure is the positive measure on $X$ defined as \begin{equation}
\beta_{u}=(i^{n^{2}}\frac{1}{N}\sum_{i=1}^{N}s_{i}\wedge\bar{s_{i}}e^{-\psi_{0}})e^{-u}\label{eq:bergman meas}\end{equation}
in terms of any given orthonormal base $(s_{i})$ in the Hilbert space
$H^{0}(X,L+K_{X})$ equipped with the Hermitian product induced by
the weight $\psi_{0}+u$ (compare section \ref{sub:Bergman-kernels}).
In particular this means that $\beta_{u}$ may be represented as $e^{-u}$
times a strictly positive smooth measure on $X$ if $L+K_{X}$ is
globally generated. The proof of formula \ref{eq:deriv of l} follows
more or less directly from the definition (see \cite{bern2} for a
geometric argument).

The following theorem, which is direct consequence of a result of
Berndtsson about the curvature of direct image bundles \cite{bern2},
considers the \emph{second} derivatives of $\mathcal{L}_{\omega_{0}}$
along a psh path. As a courtesy to the reader a new proof of the theorem,
using Bergman kernels, is given in the appendix.
\begin{thm}
\label{thm:(Berndtsson)-Let-}(Berndtsson) Let $u_{t}$ be a continuous
psh path in $\mathcal{H}_{\omega_{0}}.$ Then the function $t\mapsto\mathcal{L}_{\omega_{0}}(u_{t})$
is convex. Moreover, if $\mathcal{L}_{\omega_{0}}(u_{t})$ is affine
and $u_{t}$ is a \emph{smooth} psh path with $\omega_{u_{t}}>0$
on $X$ for all $t,$ then there is an automorphism $S_{1}$ of $(X,L)$,
homotopic to the identity, such that $u_{1}-u_{0}=S_{1}^{*}\psi_{0}-\psi_{0}.$ 
\end{thm}
The convexity statement in \cite{bern2} assumed in fact that $u_{t}$
be\emph{ smooth.} However, by uniform approximation the convexity
statement above in fact holds for any \emph{continuous} psh path in
$C^{0}(X).$ Indeed, if $u_{t}$ is such a path, then there exists,
for example by Richbergs's approximation theorem \cite{d1}, a sequence
$U^{j}$ converging uniformly towards $U$ on $M$ such that $dd^{c}U^{j}+\pi^{*}\omega_{0}>0.$
Applying the theorem above to each $U^{j}$ and letting $j$ tend
to infinity then gives that $f(t):=\mathcal{L}_{\omega_{0}}(u_{t})$
is a uniform limit of convex functions and hence convex, proving the
claim. 

However, for the statement concerning the affine properties of $\mathcal{L}_{\omega_{0}}(u_{t})$
the argument in \cite{bern2} requires that $\omega_{u_{t}}$ be reasonably
smooth in $(t,x)$ (essentially $\mathcal{C}^{3}-$smooth). Moreover,
the assumption that $\omega_{t}>0$ is crucial to be able to define
the vector fields $V_{t}$ that integrate to the automorphism $S_{1}$
(see formula \ref{eq:int multip} below).

\subsubsection{The vector field $V_{t}$ and the strict convexity of $\mathcal{L}_{\omega_{0}}$}

Given a psh path $u_{t}$ such that $u_{t}$ and $\partial_{t}u$
(the partial holomorphic derivative wrt $t)$ are smooth on $X$ with
$\omega_{u_{t}}>0$ we let $V_{t}$ be the family of vector fields
on $X$ of type $(1,0)$ defined by the equation \begin{equation}
\omega_{u_{t}}(V_{t},\cdot)=\bar{\partial}_{X}(\partial_{t}u),\label{eq:int multip}\end{equation}
 where $\bar{\partial}_{X}$ is the $\bar{\partial}-$operator on
$X.$ As shown in \cite{bern2} if the functional $\mathcal{L}(u_{t})$
is affine wrt $t$ then the vector field $V_{t}$ has to be holomorphic
on $X$ for each $t.$ Moreover, it then follows that $V_{t}$ is
holomorphic wrt $t$ as well (this uses that $u_{t}$ is automatically
a geodesic; a slight variant of this argument is recalled in section
\ref{sub:Alternative-proof-of}). Integrating $V_{t}$ will hence
give rise to the automorphism $S_{1}$ in Theorem \ref{thm:(Berndtsson)-Let-}.
As explained above a major difficulty is that, for a given geodesic
$u_{t},$ it is not known whether $\omega_{u_{t}}>0,$ nor if it is
smooth on $X\times]0,1[$ and the main part of the proof of the uniqueness
statement in Theorem \ref{thm:main} will consist in showing that
these conditions \emph{are} satisfied for a geodesic connecting two
critical points of $\mathcal{F}_{\omega_{0}}.$ In the proof of Theorem
\ref{thm:main ine one dim} we will also have to handle the case when
$\omega_{u_{t}}$ vanishes at a finite number of points. To this end
we will prove a generalization of Berndtsson's theorem in the one-dimensional
case in section \ref{sub:Conditions-for-equality} in the appendix.

It was not explicetly pointed out in \cite{bern2} that $V_{t}$ defined
above lifts to the total space of $L,$ but this fact follows from
Lemma 12 in \cite{do1}. For completeness and future reference we
next give the proof of this lifting property emphasizing that it holds
in a non-compact setting as well:
\begin{lem}
\label{lem:lifting}Let $Y$ be a (possibly non-compact) complex manifold
and $\pi:\, L\rightarrow X$ a holomorphic line bundle equipped with
a metric $h_{0}$ with strictly positive normalized curvature form
$\omega_{0}.$ If $V$ is a holomorphic vector field of type $(1,0)$
on $X$ such that \[
\omega_{0}(V,\cdot)=\bar{\partial}_{X}v\]
 for some complex valued smooth function $f$ on $Y,$ then $V$ lifts
to a holomorphic vector field on the total space of $L.$
\begin{proof}
Following \cite{don1} the lifted vector field $\tilde{V}$ may be
explictely defined by \[
\tilde{V}:=V_{hor}-vW,\]
 where $V_{hor}$ is the horisontell lift to $L$ of V (wrt the Chern
connection induced by $h_{0}$) and $W$ is the generator of the natural
$\C^{*}$ action along the fibers of $L.$ To verify that $\tilde{V}$
is indeed holomorphic fix an arbitrary point $p$ in $Y$ and a trivialization
$s$ of $L$ over a neighbourhood of $p$ such that $|s|_{h_{0}}^{2}=e^{-\phi(z)}.$
Take holomorphic coordinates $z$ centered at $p.$ Then the trivialization
of $L$ induce local holomorphic coordinates $(z,w)$ on $L$ such
that the various objects above may be expressed as follows \cite{d1}:
\[
W=w\frac{\partial}{\partial w},\,\,\, V=\phi^{i\bar{j}}v_{\bar{i}}\frac{\partial}{\partial z_{j}}\,\,\, A=\frac{dw}{w}-\phi_{j}dz_{j}\]
where $A$ is the one-form $A$ on $L$ induced by the Chern connection
(using the usual index notation and summation convention for partial
derivatives wrt $z$ and where $\phi^{i\bar{j}}$ is the matrix representing
the inverse of the matrix $\phi_{i\bar{j}}).$ Now by definition $d\pi_{*}(V_{hor})=V$
and $A[V_{hor}]=0,$ which in view of the above relations forces \[
V_{hor}=V+W\phi^{i\bar{j}}v_{\bar{i}}\phi_{j}=V+W\left\langle \partial\phi,V\right\rangle \]
Since by assumption $V$ is holomorphic $\overline{\partial}V_{hor}=W(\overline{\partial}(\phi^{i\bar{j}}v_{\bar{i}}\phi_{j}))$
and hence, using $\overline{\partial}V=0$ again gives \[
(\omega^{ij}v_{\bar{i}}\phi_{j})_{\bar{k}}=(\phi^{i\bar{j}}v_{\bar{i}})\phi_{j\bar{k}}=\phi^{i\bar{j}}v_{\bar{i}}\phi_{j\bar{k}}=v_{\bar{k}},\]
 i.e. $\overline{\partial}V_{hor}=W(\overline{\partial}v).$ All in
all this means that $\overline{\partial}\tilde{V}=0.$ 
\end{proof}
\end{lem}

\subsection{Energy functionals }

The functional $\mathcal{E}_{\omega_{0}}$ defined in the introduction
by formula \ref{eq:def of e intro} may be alternatively defined by
the following variational property: the differential of the functional
$\mathcal{E}_{\omega_{0}}$ on $\mathcal{H}_{\omega_{0}}$ is given
by the following formula:

\begin{equation}
(d\mathcal{E}_{\omega_{0}})_{u}=\omega_{u}^{n}/n!\label{eq:deriv of e}\end{equation}
$\mathcal{E}_{\omega_{0}}$ may also be expressed in terms of a generalized
Dirichlet type energy $J_{\omega_{0}}:$ \begin{equation}
-\mathcal{E}_{\omega_{0}}(u)=\frac{1}{V}(J_{\omega_{0}}(u)-\int u\omega_{0})\label{eq:relation e and j}\end{equation}
 where $J_{\omega_{0}}$ is Aubin's energy functional \begin{equation}
J_{\omega_{0}}(u):=\sum_{i=1}^{n-1}\frac{i+1}{(n+1)!}\int du\wedge du^{c}\wedge(\omega_{0})^{i}\wedge(\omega_{u})^{n-1-i}\label{eq:def of j}\end{equation}
which is clearly non-negative on $\mathcal{H}_{\omega_{0}}$ (the
definition differs from that in \cite{ti} p. 58 by a factor of $Vn!$).
Note that if $n=1$ then \[
J_{\omega_{0}}:=\frac{1}{2}\int du\wedge du^{c}\]
 is the classical Dirichlet energy on a Riemann surface, which is
independent of $\omega_{0}$ and non-negative for \emph{any} $u.$
Occasionaly, we will also use the functional \begin{equation}
I_{\omega_{0}}(u)=\int_{X}u(\omega_{0}^{n}-\omega_{u}^{n})/n!,\label{eq:def of i}\end{equation}
which is also non-negative on $\mathcal{H}_{\omega_{0}}.$ In fact,
the following well-known inequality holds \begin{equation}
I-J\geq J/n,\label{eq:ineq for i-j}\end{equation}
which is proved by a simple integration by parts argument. 

The following generalization from \cite{b-b} of formula \ref{eq:deriv of e}
to the functional $\mathcal{E}_{\omega_{0}}\circ P_{\omega_{0}},$
where $P_{\omega_{0}}$ is the non-linear projection \ref{eq:proj oper intro},
will be crucial for the proof of Theorem \ref{thm:main}:
\begin{thm}
\label{thm:deriv of composed }The functional $\mathcal{E}_{\omega_{0}}\circ P_{\omega_{0}}$
is Gâteaux differentiable on \textup{$\mathcal{C}^{0}(X).$} Its differential
at the \textup{point $u$ is represented by the measure $\omega_{P_{\omega_{0}}u}^{n}/n!,$
i.e.} given $u,v\in\mathcal{C}^{0}(X)$ the function $\mathcal{E}_{\omega_{0}}P_{\omega_{0}}(u+tv)$
is differentiable on $\R_{t}$ and \begin{equation}
d\mathcal{E}_{\omega_{0}}P_{\omega_{0}}(u+tv)/dt_{t=0}=\int_{X}v\omega_{P_{\omega_{0}}u}^{n}/n!\label{eq:deriv of comp}\end{equation}

\end{thm}
Note that for a merely continuous function $u$ the functional $\mathcal{E}_{\omega_{0}}\circ P_{\omega_{0}}$
may be defined in the sense of pluripotential theory, but in the present
case $P_{\omega_{0}}u$ is actually $C^{1,1}-$smooth \cite{be1}
and hence $\mathcal{E}_{\omega_{0}}\circ P_{\omega_{0}}$ is defined
as an integral of a density in $L_{loc}^{\infty}.$ As for the second
derivatives of $\mathcal{E}_{\omega_{0}}$ we have the following Proposition
which is well-known (at least in the smooth case):
\begin{prop}
\label{pro:e is affine}The following properties of $\mathcal{E}_{\omega_{0}}$
and $J_{\omega_{0}}$ hold:\end{prop}
\begin{itemize}
\item The functional $\mathcal{E}_{\omega_{0}}$ on $\overline{\mathcal{H}}_{\omega_{0}}\cap\mathcal{C}^{0}(X)$
is concave wrt the affine structure on $\mathcal{C}^{0}(X).$
\item Let $u_{t}$ be a $\mathcal{C}^{0}$- geodesic in $\overline{\mathcal{H}}_{\omega_{0}}$
connecting $u_{0}$ and $u_{1}.$ Then the functional $t\mapsto\mathcal{E}_{\omega_{0}}(u_{t})$
is affine and continuous on $[0,1],$ while $J_{\omega_{0}}(u_{t})$
is convex. \end{itemize}
\begin{proof}
(A proof also appears in \cite{bbgz}). Recall the following well-known
formula (see for example \cite{b-b}): \begin{equation}
d_{t}d_{t}^{c}\mathcal{E}_{\omega_{0}}(u_{t})=t_{*}(dd^{c}U+\pi^{*}\omega_{0})^{n+1}/(n+1)!,\label{eq:genaral second deriv}\end{equation}
 where $t_{*}$ denotes the natural push-forward map from $M$ to
$\C_{t}.$ In particular, setting $u_{t}=u_{0}+tu$ gives for real
$t$ $d^{2}\mathcal{E}_{\omega_{0}}(u_{t})/d^{2}t=-\int_{X}\left|\partial u_{t}\right|^{2}\omega_{0}^{n}\leq0$
(compare formula \ref{eq:geod eq concrete}) which proves the first
point of the proposition when $u$ is smooth. To handle the general
case one takes $u_{j}$ in $\mathcal{H}_{\omega_{0}}$ converging
uniformly to $u$ and uses that, according to Bedford-Taylor's classical
results, $\mathcal{E}_{\omega_{0}}$ is continuous under uniform limits
in $\overline{\mathcal{H}}_{\omega_{0}}\cap\mathcal{C}^{0}(X)$ (see
also \cite{b-b}). This shows that $\mathcal{E}_{\omega_{0}}(u_{t})$
is the limit of concave functions and hence concave. To prove the
last point take a sequence $U^{j}$ converging uniformly to $U$ on
$M$ and such that $dd^{c}U^{j}+\pi^{*}\omega_{0}>0$ (compare the
discussion below Theorem \ref{thm:(Berndtsson)-Let-}). By Bedford-Taylor
$(dd^{c}U^{j}+\pi^{*}\omega_{0})^{n+1}$ tends weakly to $(dd^{c}U^{j}+\pi^{*}\omega_{0})^{n+1}$
in the interiour of $M.$ Hence, formula \ref{eq:genaral second deriv}
shows that the second real derivatives of $\mathcal{E}_{j}(t):=\mathcal{E}(u_{t}^{j})$
tend weakly to zero in the sense of distributions for $t\in]0,1[.$
But since the sequence $\mathcal{E}_{j}(t)$ of smooth concave functions
tends to $\mathcal{E}(t)$ it follows that $\mathcal{E}(t)$ is affine
on $]0,1[.$ Finally, to see that continuity up to the boundary holds
one proceeds as follows: since $U$ is continuous on the compact set
$M$ the family $u_{t}$ tends to $u_{0}$ and $u_{1}$ uniformly
when $t\rightarrow0$ and $t\rightarrow1,$ respectively. Finally,
since $\mathcal{E}$ is continuous under uniform limits in $\overline{\mathcal{H}}_{\omega_{0}}\cap\mathcal{C}^{0}(X)$
this proves that $\mathcal{E}$ is continuous up to the boundary on
$[0,1].$ Finally, the convexity of $J_{\omega_{0}}(u_{t})$ now follows
immediately from formula \ref{eq:relation e and j} and the fact that
$u_{t}$ is convex in $t.$\end{proof}
\begin{rem}
\label{rem:j-f} The natural condition to obtain non-negativity of
the functional $J_{\omega_{0}}$ when $n>1$ is that $\omega_{u}\geq0.$
On the other hand as shown in \cite{rub0} (Lemma 2.1), there are
examples of general smooth $u$ with $J_{\omega_{0}}<0$ for any manifold
$X$ of dimension $n>1.$ For a completeness we give a simple proof
of this fact when $n=2$ which is valid for the functional $I_{\omega_{0}}-J_{\omega_{0}}$
as well: both these functionals may be decomposed as \[
A\int(-u)(ddu^{c})^{2}+B\int(-u)(ddu^{c})\wedge\omega_{0},\]
 where $A$ and $B$ are strictly positive constants. Now take $u$
such that the first term above is non-zero. Replacing $u$ with $-u$
if necesseary we may even assume that the first term is strictly\emph{
negative.} Finally, replacing $u$ with a sufficently large positive
multiple $tu$ we see that the whole expression is negative (more
precisely the first term is of the order $t^{3},$ while the second
term is of the order $t^{2})$ proving the fact alluded to above.
As a direct consequence it was shown in \cite{rub0}, in the case
$L=-K_{X},$ that any such fucntion $u$ violates the inequality in
Theorem \ref{thm:main}. A similar argument applies to any homogeneous
line bundle $L$ as in Corollary \ref{cor:homeg}. Indeed, without
effecting the value of $J_{\omega_{0}}(u)$ we may assume that $\int_{X}u\omega_{0}^{n}=0$
so that $-\mathcal{E}_{\omega_{0}}(u)=J_{\omega_{0}}(u)<0.$ Now,
using the notation explained in the beginning of of section \ref{sub:Bergman-kernels}
in the appendix \[
-\mathcal{L}_{\omega_{0}}(u)=\log\E_{\psi_{0}}(e^{-(u(x_{1})+...+u(x_{n})})\geq-\E_{\psi_{0}}(u(x_{1})+...+u(x_{n})),\]
using Jensen's inequality in the last step. Moreover, by formula \ref{eq:beta as one-pt}
$\E_{N}(u(x_{1})+...+u(x_{n}))=\int_{X}u\beta_{u}).$ Since $\beta_{u}=\omega_{0}^{n}/V$
in the homogenous case (compare the proof of Corollary \ref{cor:homeg}),
this means that $-\mathcal{L}_{\omega_{0}}(u)\geq0.$ Hence, $u$
violates the inequality referred to above.
\end{rem}
We also recall the following basic cocycle property of the functional
$\mathcal{F}_{\omega_{0}}:=\mathcal{E}_{\omega_{0}}-\mathcal{L}_{\omega_{0}}:$
\begin{equation}
\mathcal{F}_{\omega_{u_{2}}}(u_{1})+\mathcal{F}_{\omega_{u_{3}}}(u_{2})=\mathcal{F}_{\omega_{u_{3}}}(u_{1}),\label{eq:cocy}\end{equation}
 which is a direct consequences of the corresponding cocyle properties
of $\mathcal{E}_{\omega_{0}}$ and $\mathcal{L}_{\omega_{0}}.$ These
latter properties in turn are immediately obtained by integrating
the corresponding differentials along line segments in $\mathcal{H}_{\omega_{0}}$(compare
\cite{ti}).

\subsection{\label{sub:Analytic-torsion}Analytic torsion}

Let $F$ be a holomorphic vector bundle over the Kähler manifold $(X,\omega_{0}).$
Given an Hermitian metric $h$ on $F$ the pair $(h,\omega_{0})$
induces, in the standard way, Hermitian products on the space $\Omega^{0,q}(X,F)$
of smooth $(0,q)-$forms with values in $F.$ The corresponding\emph{
analytic torsion, }first introduced by Ray-Singer,\emph{ }is defined
as \[
\tau_{F}(h,\omega_{0}):=-\sum_{q=0}^{n}(-1)^{q}q\log\det\Delta_{\bar{\partial}}^{(q)}\]
(in our sign convention), where $\Delta_{\bar{\partial}}^{(q)}:=\bar{\partial}\bar{\partial}^{*}+\bar{\partial}^{*}\bar{\partial}$
is the Dolbeault Laplacian acting on $\Omega^{0,q}(X,F)$ and $\log\det\Delta_{\bar{\partial}}^{(q)}:=-\frac{\partial\zeta^{(q)}}{\partial s}_{s=0},$
where $\zeta^{(q)}(s)=\sum_{i}(\lambda_{i}^{(q)})^{-s}$ is the meromorphic
continuation to $\C_{s}$ of the zeta function for the \emph{strictly}
postive eigenvalues $\{\lambda_{i}^{(q)}\}$ of $\Delta_{\bar{\partial}}^{(q)}$
(see \cite{b-g-s,g-s} and references therein). We will be concerned
mainly with the case when $n=1,$ where \[
\tau_{F}(h,\omega_{0})=\log\det\Delta_{\bar{\partial}}^{(1)}=\log\det\Delta_{\bar{\partial}}^{(0)}\]
 (using that the non-zero eigenvalues of $\Delta_{\bar{\partial}}^{(0)}$
and $\Delta_{\bar{\partial}}^{(1)}$ coincide) and study the dependence
of $\tau_{F}(h,\omega_{0})$ on $h.$ The following proposition is
a special case of the general anomaly formula of Bismut-Gillet-Soulé
\cite{b-g-s}:
\begin{prop}
\label{pro:f as anal tors on curve}Let $F$ be a line bundle over
a complex curve $X.$ Let $\omega_{0}$ be an Kähler metric on $X$
with constant Gauss curvature (as a metric on $TX).$ Assume that
$L:=F-K_{X}$ is ample and that $\omega_{0}$ is normalized so that
$\omega_{0}\in c_{1}(L).$ Then \[
\tau_{F}(e^{-u}h_{0},\omega_{0})-\tau_{F}(h_{0},\omega_{0})=\frac{\mbox{deg}(K_{X})}{\mbox{2deg}(L)}J(u)+N\mathcal{F}_{\omega_{0}}(u),\]
 where $h_{0}$ is a metric on $F$ whose normalized curvature form
equals $\omega_{0}$ and $N=\dim H^{0}(X,F).$ \end{prop}
\begin{proof}
It will convenient to use weight notation, i.e. $h_{0}=e^{-\psi_{0}}$
etc. The Hermitian product on $H^{0}(X,L+K_{X})$ induced by a given
weight $\psi=\psi_{0}+u$ on $L$ may be expressed as \[
\left\langle s,t\right\rangle _{\psi}:=i\int s\wedge\bar{t}e^{-\psi}=\int s\bar{t}e^{-(\psi+\psi_{\omega_{0}})}\omega_{0},\]
 where $\psi_{\omega_{0}}:=\log\omega_{0}$ defines a weight on $K_{X}$
such that its curvature $dd^{c}\psi_{\omega_{0}}:=-\mbox{Ric}\omega_{0}.$
In other words, the Hilbert space $H^{0}(X,L+K_{X})$ (with the Hermitian
product determined by $\psi)$ is unitary equivalent to the Hilbert
space $H^{0}(X,F)$ determined by the weight $\phi:=\psi+\psi_{\omega_{0}}$
on $F$ and the metric $\omega_{0}$ on the tangent bundle $TX.$
Next, we will use the anomaly formula of Bismut-Gillet-Soulé \cite{b-g-s}
for the Quillen metric on the determinant line $\bigwedge^{N}H^{0}(X,F)$
(note that $H^{1}(X,F)$ is trivial since $L$ is ample). It reads
as follows in our notation

\begin{equation}
\tau_{F}(e^{-u}h_{0},\omega_{0})-\tau_{F}(h_{0},\omega_{0})=\int_{X}\mbox{Td}(X,\omega_{0})\wedge\tilde{ch}(\phi_{0}+u,\phi_{0})-N\mathcal{L}_{\omega_{0}}(u),\label{eq:pf of prop anal tors}\end{equation}
 where \[
\mbox{Td}(X,\omega_{0})=(1+\frac{1}{2}\mbox{Ric}\omega_{0})=(1-\frac{\mbox{deg}(K_{X})}{\mbox{2deg}(L)}\omega_{0})\]
 is the \emph{Todd class }of $TX$ represented by the constant curvature
metric $\omega_{0}\in c_{1}(L)$ and \[
\tilde{ch}(\phi_{1},\phi_{0})=(\phi_{1}-\phi_{0})(1+\frac{dd^{c}\phi_{1}+dd^{c}\phi_{0}}{2})\]
 is the \emph{Bott-Chern class} of the pair of weights $(\phi_{1},\phi_{0})$
on $F$ associated to the \emph{Chern character} of $F.$ Since by
assumption $dd^{c}\psi_{\omega_{0}}=-\mbox{Ric}\omega_{0}=\frac{\mbox{deg}(K_{X})}{\mbox{deg}(L)}\omega_{0}$
this means that \[
\tilde{ch}(\phi_{0}+u,\phi_{0})=u(1+\frac{\omega_{u}+\omega_{0}(1+\frac{\mbox{\mbox{deg}}(K_{X})}{\mbox{\mbox{deg}}(L)})}{2}).\]
Expanding the integrand in \ref{eq:pf of prop anal tors} hence gives
\[
\tau_{F}(e^{-u}h_{0},\omega_{0})-\tau_{F}(h_{0},\omega_{0})=\frac{\mbox{deg}(K_{X})}{2\mbox{deg}(L)}\int_{X}u\omega_{0}\mbox{+\mbox{deg}}(L)\mathcal{E}_{\omega_{0}}(u)-N\mathcal{L}_{\omega_{0}}(u)).\]
Now by the Riemann-Roch theorem $N=\mbox{deg}(F)-\mbox{deg}(K_{X})/2=\mbox{deg}(L)+\mbox{deg}(K_{X})/2,$
which finally proves the proposition using the relation \ref{eq:relation e and j}. \end{proof}
\begin{rem}
In the general anomaly formula in \cite{b-g-s} the metric $\omega_{0}$
is allowed to vary as well. In particular, when $L=\mathcal{O}(0)$
is the trivial holomorphic line bundle over $S^{2},$ the metric $h=1$
is kept constant, but the conformal metric $g_{u}=e^{-u}g_{0}$ on
$TS^{2}$ varies with $u,$ the anomaly formula in \cite{b-g-s} is
equivalent to Polyakov's formula and then $\log(\frac{\det\Delta_{g_{u}}}{\det\Delta_{g_{0}}})$
coincides with the functional $\mathcal{F}_{0}$ (up to a a multiplicative
constant) \cite{ch}. 

Next, we follow standard notation in Kähler geometry and denote by
$h_{\omega_{0}}$ any Ricci potential of $\omega_{0},$ i.e. \begin{equation}
\mbox{Ric \ensuremath{\omega_{0}-\mu\omega_{0}=dd^{c}h_{0}.}}\label{eq:ricci pot eq}\end{equation}
\end{rem}
\begin{prop}
\label{pro:form for torsion in general n}Fix $\mu\in\{0,-1,1\}$
and assume that $L:=\mu K_{X}$ is ample line bundle or trivial (
if $\mu=0).$ Then $(\tau_{F}(e^{-ku}h_{0},\omega_{0})-\tau_{F}(h_{0},\omega_{0}))/N_{k}=\mathcal{F}_{k\omega_{0}}(u)+$\[
+\frac{1}{2\mbox{Vol}(L)}\left(\int h_{\omega_{0}}(\frac{\omega_{0}^{n}}{n!}-\frac{\omega_{u}^{n}}{n!})-\mu(I_{\omega_{0}}(u)-J_{\omega_{0}}(u))+O(\frac{1}{k}))J_{\omega_{0}}(u)+O(\frac{1}{k})\right)\]
where the error term $O(\frac{1}{k})$ means that $\left|O(\frac{1}{k})\right|\leq C/k,$
where $C$ only depends on $\omega_{0}\in c_{1}(L).$ \end{prop}
\begin{proof}
By scaling invariance we may assume that $\sup_{X}u=0.$ Applying
the anomaly formula as above gives the formula in the proposition
but with a sum of error terms of the form $O(\frac{1}{k})$ times
\[
\int_{X}(-u(\omega_{u})^{i}\wedge(\omega_{0})^{j}\wedge\beta_{i,j}^{n-(i+j)})\]
 where $i+j<n$ and $\beta_{i,j}$ is an $(n-(i+j),n-(i+j))$ form
depending on $\omega_{0}$ only. Since there is an alternative way
to get the explicit form of the leading term we skip the calculation
(compare the remark below). Now, by a well-known inequality for weakly
positive forms \cite{h-k} there is a constant $C$ such that \[
-C(\omega_{0})^{n-(i+j)}\leq\beta_{i,j}\leq C(\omega_{0})^{n-(i+j)}\]
 Hence, we can conclude the proof of the proposition by observing
that, since $\sup_{X}u=0,$ \[
|\int_{X}(u(\omega_{u})^{i}\wedge(\omega_{0})^{n-i}|\leq C(J_{\omega_{0}}(u)-\int u\omega_{0}^{n}))\]
 using formula \ref{eq:relation e and j}. Finally, we just use the
basic fact that $u\mapsto\int u\omega_{0}^{n}$ is a bounded functional
on the subspace of $\mathcal{H}_{\omega_{0}}$ where $\sup_{X}u=0$
(see \cite{g-z}). \end{proof}
\begin{rem}
\label{rem:torsion entropy}Letting $k$ tend to infinity in the previous
formula gives \[
\lim_{k\rightarrow\infty}(\tau_{F}(e^{-ku}h_{0},\omega_{0})-\tau_{F}(h_{0},\omega_{0}))/k^{n}=\]
\[
=-\mathcal{S}(\frac{\omega_{u}^{n}}{\mbox{Vol}(L)n!},\frac{\omega_{0}^{n}}{\mbox{Vol}(L)n!}):=-\frac{1}{2\mbox{Vol}(L)}\int_{X}\log(\frac{\omega_{u}^{n}}{\omega_{0}^{n}})\frac{\omega_{u}^{n}}{n!}\]
 which gives a new proof of a result in \cite{b-v} (which is valid
for any ample line bundle). To see this one first uses that $\mathcal{F}_{k\omega_{0}}(u)$
converges to $-\mathcal{M}_{\omega_{0}}$ where $\mathcal{M}_{\omega_{0}}$
denotes Mabuchi's $K-$energy with our normalizations (formula \ref{eq:asympt of f}).
Finally, by the explicit formula of Tian (\cite{ti}, page 95) and
Chen \cite{ch2} for $\mathcal{M}_{\omega_{0}}$ associated to an
ample line bundle $L$ this then proves the asymptotics above. It
is interesting to see that the limiting funtional (i.e. minus the
relative entropy of the corresponding probability measures) is always
bounded from above and maximized precisely for $u$ such that $\omega_{u}^{n}=\omega_{0}^{n}$
(by Jensen's inequality), which in turn means that $u$ is constant
by well-known uniqueness results for the Monge-Ampère equation. 
\end{rem}

\section{Proofs of the main results}

\subsection{\label{sub:Proof-of-Theorem-main}Proof of Theorem \ref{thm:main}}

By the cocycle property of $\mathcal{F}_{\omega_{0}}$ (see \cite{bbgz,b-b})
we may without loss of generality assume that $u=0$ is critical.
Take a \emph{continuous} element $u_{1}$ in $\mathcal{H}_{\omega_{0}}$
and the corresponding $\mathcal{C}^{0}-$geodesic $u_{t}$ connecting
$u_{0}=0$ and $u_{1}.$ Since $u_{t}$ is a continuous path, combining
Theorem \ref{thm:(Berndtsson)-Let-} and Proposition \ref{pro:e is affine}
gives that $\mathcal{F}_{\omega_{0}}(t):=\mathcal{F}_{\omega_{0}}(u_{t})$
is a continous concave function on $[0,1].$ Hence, the inequality
in Theorem \ref{thm:main} will follow once we have shown that \begin{equation}
\frac{d}{dt}_{t=0+}\mathcal{F}(u_{t})\leq0.\label{eq:ineq for right der}\end{equation}
Of course, if $u_{t}$ were known to be a \emph{smooth} path then
this would be an immediate consequence of the assumption that $u_{0}$
is critical combined with the chain rule (which would even yield equality
above). To prove \ref{eq:ineq for right der} first observe that by
the concavity in Prop \ref{pro:e is affine}\[
(\mathcal{E}_{\omega_{0}}(u_{t})-\mathcal{E}_{\omega_{0}}(u_{0}))/t\leq\frac{1}{t}\int_{X}(u_{t}-u_{0})(\omega_{u_{0}})^{n}/n!\]
Hence, the monotone convergence theorem applied to the sequence $(u_{t}-u_{0})/t$
which decreases to the right derivative $v_{0}$ of $u_{t}$ at $t=0$
(using that $u_{t}$ is convex in $t)$ gives \begin{equation}
\frac{d}{dt}_{t=0+}\mathcal{E}_{\omega_{0}}(u_{t})\leq\int_{X}v_{0}(\omega_{u_{0}})^{n}/n!\label{eq:upper bound on deriv of energy}\end{equation}
Hence, \[
\frac{d}{dt}_{t=0+}\mathcal{F}(u_{t})\leq\int_{X}((\omega_{u_{0}})^{n}/n!-\beta_{u_{u}})v_{0}=0,\]
 where we have also used the dominated convergence theorem to differentiate
$\mathcal{L}_{\omega_{0}}(u_{t})$ (compare \cite{bbgz,b-b}). This
finishes the proof of \ref{eq:ineq for right der}and hence the first
statement in the theorem follows.

\emph{Uniqueness:} 

\emph{Step 1 ($u_{t}$ is critical for all $t)$:} Assume now that
$u_{1}$ is a smooth maximizer of $\mathcal{F}_{\omega_{0}}$ on $\overline{\mathcal{H}}_{\omega_{0}}$
i.e. that $\mathcal{F}_{\omega_{0}}(u_{1})=\mathcal{F}_{\omega_{0}}(u_{0})$
by the previous step. Since $\mathcal{F}_{\omega_{0}}(t):=\mathcal{F}_{\omega_{0}}(u_{t})$
is continuous and concave it follows that $u_{t}$ maximizes $\mathcal{F}_{\omega_{0}}$
on $\overline{\mathcal{H}}_{\omega_{0}}\cap\mathcal{C_{\C}}^{1,1}(X)$
for all $t.$ Next, we will show that $u_{t}$ satisfies the Euler-Lagrange
equation \ref{eq:e-l for f} for any fixed $t$ (see \cite{bbgz}
for similar arguments). To this end fix $t=t_{0}$ and set $u_{t_{0}}:=u.$
Given a smooth function $v$ on $X$ consider the function $f(t):=\mathcal{E}_{\omega_{0}}(P_{\omega_{0}}(u+tv))-\mathcal{L}_{\omega_{0}}(u+tv)$
on $\R_{t}.$ Since, the functional $\mathcal{L}_{\omega_{0}}$ is
increasing on $\mathcal{C}^{0}(X)$ we have $f(t)\leq\mathcal{F}_{\omega_{0}}(P_{\omega_{0}}(u+tv)).$
By assumption this means that the maximal value of the function $f(t)$
is attained for $t=0$ (also using that $P_{\omega_{0}}u=u).$ In
particular, since by Theorem \ref{thm:deriv of composed } $f(t)$
is differentiable $df/dt=0$ at $t=0$ and Theorem \ref{thm:deriv of composed }
and formula \ref{eq:deriv of l} hence show that the Euler-Lagrange
equation \ref{eq:e-l for f} holds (since it holds when tested on
any smooth function $v).$ 

\emph{Step 2 ($U\in\mathcal{C}^{\infty}(\dot{M})):$ }(where $\dot{M}$
denotes the interiour of $M).$ By Theorem \ref{thm:(Chen)-Assume-that}
$U$ is in $\mathcal{C}_{\C}^{1,1}(M).$ Moreover, by the homogenous
Monge-Ampère equation \ref{eq:dirichlet pr} and the Euler-Lagrange
equation \ref{eq:e-l for f} we have \[
(dd^{c}(U+\left|w\right|^{2})+\pi_{X}^{*}\omega_{0})^{n+1}=i\beta_{u}\wedge dw\wedge d\bar{w}\]
 Hence, the following equation holds locally on $\C^{n+1}$ (where
we for simplicity have kept the notation $U$ for the function obtained
after subtracting a smooth and hence harmless function from $U):$
\begin{equation}
\det(\partial_{\zeta_{i}}\partial_{\bar{\zeta}_{j}}U)=e^{-U}\rho,\label{eq:local exp ma}\end{equation}
 where $\rho$ is a positive smooth function, depending on $U$ (compare
the discussion below formula \ref{eq:bergman meas}). In particular,
$\det(\partial_{\zeta_{i}}\partial_{\bar{\zeta}_{j}}U)$ is locally
in $\mathcal{C}_{\C}^{1,1}.$ But then Theorem 2.5 in \cite{bl0},
which is a complex analog of a result of Trudinger for fully non-linear
elliptic operators (compare Evans-Krylov theory), gives that $U$
is locally in the Hölder space $\mathcal{C}^{2,\alpha}$ for some
$\alpha>0.$ Now the equation \ref{eq:local exp ma} shows that $\det(\partial_{\zeta_{i}}\partial_{\bar{\zeta}_{j}}U)$
is also in $\mathcal{C}^{2,\alpha}.$ Finally, since we have hence
shown that $U\in C^{2},$ standard theory of uniformly elliptic operators
then allows us to boot strap using \ref{eq:local exp ma} and deduce
that $U\in\mathcal{C}^{\infty}$ locally (see Theorem 2.2 in \cite{bl}).

\emph{Step 3: the inequalities \begin{equation}
1/C'\omega_{0}\leq\omega_{u_{t}}\leq C'\omega_{0}\label{eq:bound on omega t}\end{equation}
hold. }To see this first note that by the Euler-Lagrange equation
\ref{eq:e-l for f} we have a uniform lower bound $\omega_{u_{t}}^{n}>\delta\omega_{0}^{n}$
(also using the lower bound in formula \ref{eq:bounds on beta} in
the appendix). Combining the previous lower bound with the upper bound
$\omega_{u_{t}}\leq C\omega_{0}$ from Theorem \ref{thm:(Chen)-Assume-that}
then shows that there is a positive constant $C',$independent of
$t,$ such that \ref{eq:bound on omega t} holds.

\emph{Step 4: Application of {}``strict convexity'':}

By the above arguments $\mathcal{F}_{\omega_{0}}(u_{t})$ and $\mathcal{E}_{\omega_{0}}(u_{t})$
are both affine (and even constant) and hence it follows that $\mathcal{L}_{\omega_{0}}(u_{t})$
is affine. In case $U$ were smooth \emph{up to the boundary} of $M$
applying Theorem \ref{thm:(Berndtsson)-Let-} would hence prove the
uniqueness statement in Theorem \ref{thm:main}. To prove the general
case we may without loss of generality assume that $u_{t}(x)$ is
smooth on $[0,1[\times X$ (otherwise we just apply the same argument
on $[1/2,1[$ and $]0,1/2]).$ For any $\epsilon>0$ Theorem \ref{thm:(Berndtsson)-Let-}
(see also Theorem 2.6 in \cite{bern2}) furnishes a 1-parameter holomorphic
family $S_{t}$ in $\mbox{Aut}_{0}(X,L)$ with $t\in[0,1-\epsilon]$
defined by the ordinary differential equation \begin{equation}
\frac{dS_{t}(x(t))}{dt}=d_{X}(S(x(t))[V_{t}]_{x(t)}\label{eq:ode}\end{equation}
with the iniatial data $S_{0}=I$ (the identity), where $V_{t}$ is
the vector field on $X$ of type $(1,0)$ defined by the equation
\ref{eq:int multip}. As shown in \cite{bern2} $V_{t}$ defines a
holomorphic vector field on the interiour of $A\times$ $X.$ Furthermore,
as shown in \cite{bern2} \begin{equation}
\psi_{t}-S_{t}^{*}\psi_{0}=C_{t}\label{eq:action on weights}\end{equation}
 where $\psi_{t}=\psi_{0}+u_{t}$ and $C_{t}$ is a constant for each
$t,$ i.e. \begin{equation}
\omega_{u_{t}}=S_{t}^{*}\omega_{0}.\label{eq:action on curvature}\end{equation}
Now, by the bound \ref{eq:bound on omega t} on $\omega_{u_{t}}$
the point-wise norm of the vector field $V_{t}$ wrt the metric $\omega_{0}$
is uniformly bounded in $t$ on all of $X.$ Hence, the equation \ref{eq:ode}
and a basic normal families argument applied to the family $S_{t}$
yields a subsequence $S_{t_{j}}$ and a holomorphic map $S_{1}$ on
$X$ such that $S_{t_{j}}(x)\rightarrow S_{1}(x)$ uniformly on $X$
(wrt the distance defined by the metric $\omega_{0})$ where $S_{1}$
is a biholomorphism according to the relation \ref{eq:action on curvature}.
Finally, letting $t_{j}\rightarrow1$ in the relation \ref{eq:action on weights}
and using that $u_{t}$ is continuous on $[0,1]\times X$ finishes
the proof of the uniqueness statement in the theorem.

\subsection{Proof of Corollary \ref{cor:homeg}}

First observe that we may assume that $H^{0}(X,L+K_{X})$ has a non-zero
element (otherwise the corrollary is trivally true). But since $(X,L)$
is homogenuous it then follows immediately that $L+K_{X}$ is globally
generated. Hence, the conditions in Theorem \ref{thm:main} are satisfied.

Assume now that $\omega_{0}$ is invariant under the holomorphic and
transitive action of $K$ on $X.$ Then it follows that $0$ is a
critical point. Indeed, the volume form $\omega_{0}^{n}/n!$ is invariant
under the action of $K$ on $X$ and so is the Bergman measure $\beta(0)$
(since it is defined in terms of the $K-$invariant weight $\psi_{0}).$
Since the action of $K$ is transitive and both measures are normalized
it follows that the function $(\omega_{0}^{n}/n!)/\beta(0)$ on $X$
is constant and hence equal to one. In other words, $0$ is a critial
point and by Theorem \ref{thm:main} the inequality in the statement
of Corollary \ref{cor:homeg} then holds. Finally, the last statement
of the corollary is a direct consequence of the uniqueness part of
Theorem \ref{thm:main}.

\subsection{Proof of Corollary \ref{cor:moser}}

Let us first prove the first statement of the corollary. Since $\mathcal{C}^{\infty}(X)$
is dense in $W^{1,2}(X)$ we may assume that $u$ is smooth. First
observe that \begin{equation}
\mathcal{F}_{\omega_{0}}(u)\leq\mathcal{F}_{\omega_{0}}(P_{\omega_{0}}u).\label{eq:pf of cor mos}\end{equation}
 To see this note that, since, by definition, $P_{\omega_{0}}u\leq u$
the fact that $\mathcal{L}_{\omega_{0}}$ is increasing immediately
implies $\mathcal{L}_{\omega_{0}}(u)\geq\mathcal{L}_{\omega_{0}}(P_{\omega_{0}}u).$
Next, observe that by the cocycle property of $\mathcal{F}_{\omega_{0}}(u)$
\[
\mathcal{E}_{\omega_{0}}(u)=\mathcal{E}_{\omega_{0}}(P_{\omega_{0}}u)+\int_{X}(u-P_{\omega_{0}}u)(\omega_{u}+\omega_{P_{\omega_{0}}u})/2\]
 But, since, as is well-known the measure $\omega_{P_{\omega_{0}}u}$
is supported on the open set $\{u>P_{\omega_{0}}u\}$ (cf. Prop. 1.10
in \cite{b-b} for a generalization) we have that the last term above
is equal to \[
\int_{X}(u-P_{\omega_{0}}u)(\omega_{u}-\omega_{P_{\omega_{0}}u})/2=\int_{X}(u-P_{\omega_{0}}u)(dd^{c}(u-P_{\omega_{0}}u)=\]
\[
=-\int_{X}d(u-P_{\omega_{0}}u)\wedge d^{c}(u-P_{\omega_{0}}u)\leq0,\]
 where we have integrated by parts in the last equality, which is
justified since, for example, by Theorem \cite{be1} $P_{\omega_{0}}u$
is in $\mathcal{C}^{1,1}(X)$ (but using that $P_{\omega_{0}}u$ is
in $\mathcal{C}^{0}(X)$ is certainly enough by classical potential
theory). Hence, $\mathcal{E}_{\omega_{0}}(u)\leq\mathcal{E}_{\omega_{0}}(P_{\omega_{0}}u)$
which finishes the proof of \ref{eq:pf of cor mos}. Since, $\omega_{P_{\omega_{0}}u}\geq0$
uniform approximation let's us apply Corollary \ref{cor:homeg} to
deduce \[
\mathcal{F}_{\omega_{0}}(u)\leq\mathcal{F}_{\omega_{0}}(P_{\omega_{0}}u)\leq0\]
 which proves the first statement of the corollary. 

Finally, the uniqueness will follow from Corollary \ref{cor:homeg}
once we know that a maximizer $u$ of $\mathcal{F}_{\omega_{0}}$
on $W^{1,2}(S^{2})$ is smooth with $\omega_{u}>0.$ By the previous
step we may assume that $\omega_{u}\geq0.$ But since $W^{1,2}(S^{2})$
is a linear space containing $\mathcal{C}^{\infty}(X)$ the Euler-Lagrange
equations $\omega_{u_{0}}+dd^{c}u=\beta(u)$ hold for the maximizer
$u.$ Since $\beta(u)=e^{-u}\rho>0$ with $\rho$ smooth, local elliptic
estimates for the Laplacian then show that $u$ is in fact smoth with
$\omega_{t}>0.$ All in all we have proved that \begin{equation}
-\mathcal{L}_{k\omega_{0}}(u)\leq-\mathcal{E}_{k\omega_{0}}(u)\label{eq:abs ineq in pf cor fang}\end{equation}
 for $L=k\mathcal{O}(1)$ with conditions for equality.

\subsubsection{\label{sub:Explicit-expression}Explicit expression}

To make the previous inequality more explicit note that, by definition,
\[
\mathcal{E}_{k\omega_{0}}(u):=\frac{1}{2\int k\omega_{0}}\int(udd^{c}u+u2k\omega_{0})=\frac{1}{2k}\int udd^{c}u+\int u\omega_{0})\]
 Moreover, since for $X=\P^{1}$ we have $K_{X}=-\mathcal{O}(2)$
it follows that $L+K_{X}=\mathcal{O}(k-2)=:\mathcal{O}(m).$ Under
this identification the scalar product on $H^{0}(X,L+K_{X})$ may
be written as \[
\left\langle s,t\right\rangle _{k\psi_{0}+u}=c\int s\bar{t}e^{-(u+m\psi_{0})}\omega_{0}\]
 using that $\omega_{0}$ is a Kähler-Einstein metric, i.e. $\omega_{0}(z):=dd^{c}\psi_{0}=ce^{-2\psi_{0}}idz\wedge d\bar{z}$
for some numerical constant $c.$ Since the funtional $\mathcal{L}$
is invariant under on overall scaling in definition of the scalar
product $\left\langle \cdot,\cdot\right\rangle _{\psi}$ we may as
well assume that $c=1.$ Hence, since $N_{m}=m+1,$ we have\begin{equation}
\mathcal{L}_{m}(u):=(m+1)\mathcal{L}_{\omega_{0,k}}(u)=-\log\mbox{det}(c_{i}c_{j}\int_{\C}\frac{z^{i}\bar{z}^{j}}{(1+z\bar{z})^{m}}e^{-u}\omega_{0}),\label{eq:def of l m in pf fang}\end{equation}
 where $c_{i}=(\int\frac{\left|z^{i}\right|^{2}}{(1+z\bar{z})^{m}}\omega_{0})^{-1/2}.$
Hence, the inequality \ref{eq:abs ineq in pf cor fang} may be expressed
as \begin{equation}
\log\mbox{det}(c_{i}c_{j}\int_{\C}\frac{z^{i}\bar{z}^{j}}{(1+z\bar{z})^{m}}e^{-u}\omega_{0})\leq-\frac{m+1}{(m+2)}\frac{1}{2}\int(udd^{c}u)-(m+1)u\omega_{0}.\label{eq:ineq for l m in pf fang}\end{equation}
 In particular, when $m=0$ the inequality above reads \[
\log(\int_{S^{2}}e^{-u}\omega_{0})\leq\frac{1}{4}\int(udd^{c}u)+\int u\omega_{0}).\]
Finally, to compare with the notation of Onofri \cite{on}, note that,
by definition, $dd^{c}u=\frac{i}{2\pi}\partial\bar{\partial u}$ and
hence, integration by parts gives, \[
-\int udd^{c}u=\frac{1}{\pi}\frac{i}{2}\int\partial u\wedge\bar{\partial u}.\]
Moreover, in terms of a given local holomorphic coordinate $z=x+iy,$
we have $\frac{i}{2}\partial u\wedge\bar{\partial u}=\frac{1}{4}\left|\nabla u\right|^{2}dx\wedge dy$,
where $\nabla=(\partial_{x},\partial_{y})$ is the gradient wrt the
local Euclidian metric. By conformal invariance we hence obtain $-\int udd^{c}u=\frac{1}{4\pi}\int\left|\nabla u\right|^{2}d\mbox{Vol}_{g}$
for\emph{ any} Riemannian metric $g$ on $S^{2}$ conformally equivalent
to $g_{0}.$ In particular, taking $g$ as the usual round metric
on $S^{2}$ induced by its embedding as the unit-sphere in Eucledian
$\R^{3}$ finally gives \[
\log(\int_{S^{2}}e^{-u}d\mbox{Vol}{}_{g}/4\pi)\leq\frac{1}{4}\int\left|\nabla u\right|^{2}-u)d\mbox{Vol}{}_{g}/4\pi),\]
 using that $\omega_{0}=d\mbox{Vol}{}_{g}/4\pi.$ This is precisely
the inequality proved by Onofri \cite{on}.

\subsection{\label{sub:Proof-of-Theorem main one dim}Proof of Theorem \ref{thm:main ine one dim}}

\subsubsection*{Uniqueness }

First note that, by the basic short exact sequence for restriction
to a divisor (here the point $p$) \cite{d1}, $K_{X}+L$ is not globally
generated iff there is a point $p\in X$ such that $H^{1}(L-L_{p}+K_{X})\neq\{0\},$
where $L_{p}$ is the holomorphic line bundle which has a holomophic
section $s_{p}$ vanishing to order one precisely at $p.$ By Serre
duality means that $H^{0}(L_{p}-L)\neq\{0\}.$ But since $\deg L\geq1$
this happens precisely when $\deg L=1$ and $L=L_{p}.$ Moreover,
it then follows from the Riemann-Roch theorem that $\dim H^{0}(X,L+K_{X})=$
\[
=\deg L+\frac{1}{2}\deg K_{X}=1+\frac{2g-2}{2}=\dim H^{0}(X,K_{X})\]
 and hence that \begin{equation}
H^{0}(X,L+K_{X})=H^{0}(X,K_{X})\otimes s_{p}.\label{eq:decomp in pf uniq deg one}\end{equation}
 As a consequence the Bergman measure can be factorized as (after
replacing $s_{p}$ with $s_{p}/\left\Vert s_{p}\right\Vert _{\psi_{0}+u})$
\[
\beta_{u}=i|s_{p}|^{2}e^{-(\psi_{0}+u)}\omega_{0}f_{u},\]
 where \[
f_{u}=\frac{1}{g}\sum_{i}\alpha_{i}\wedge\bar{\alpha}_{i}/\omega_{0}>0\]
 and where $\alpha_{i}$ is a base of holomorphic $(1,0)-$forms in
$H^{1,0}(X)=H^{0}(K_{X})$ which is orthonormal with respect to the
Hermitian product on $K_{X}$ obtained by twisting the canonical Hermitian
product by $|s_{p}|_{\psi_{0}}^{2}e^{-u}$ i.e. \[
\int_{X}\alpha_{i}\wedge\bar{\alpha}_{j}|s_{p}|_{\psi_{0}}^{2}e^{-u}=\delta_{ij}\]
In particular, when $(X,\omega_{0})$ is a torus with its invariant
metric we have $f_{u}=1.$ The equation for a critical point, normalized
so that $\mathcal{L}_{\omega_{0}}(u)=0$ hence reads\begin{equation}
\omega_{0}+dd^{c}u=|s|_{\psi_{0}}^{2}e^{-u}f_{u}\omega_{0}\label{eq:critical point eq in one dim}\end{equation}
Now assume that $u_{0}$ and $u_{1}$ both satisfy the previous equation
(by basic linear elliptic theory we may assume that $u_{i}$ are smooth).
As will be shown in section \ref{sub:The-one-dimensional-case} the
geodesic $U:=u_{t}$ connecting $u_{0}$ and $u_{1}$ has suitable
regularity properties. In particular $u_{t}$ is smooth on $X$. Let
$V_{t}$ be the vector field of type $(1,0)$ defined by the relation
\ref{eq:int multip} on $X-\{p\},$ where we recall that $p=\{s_{p}=0\}.$
It will be enough to prove that $V_{t}=0$ for all $t.$ Indeed, it
then follows immediately from the definition of $V_{t}$ that $\bar{\partial}(\partial_{t}u)=0$
on $X$ and hence $\partial_{t}u$ is constant on $X$ for all $t.$
Integrating over $t\in[0,1]$ then show that $u_{1}-u_{0}$ is constant
on $X,$ proving the uniqueness.

\subsubsection{Proof of $V_{t}=0$ \label{sub:Proof-of-v vanishes}}

Just as in section \ref{sub:Proof-of-Theorem-main} we have that $\mathcal{L}(u_{t})$
is affine wrt $t.$ Hence Corollary \ref{cor:strict convex in append}
in the appendix yields a holomorphic section $h^{1,0}$ of $L+K_{X}$
such that $h^{1,0}\wedge\bar{\partial}(\partial_{t}\psi_{t})/\partial_{X}\bar{\partial_{X}}\psi_{t}$
defines a holomorphic section of $L$ over $X.$ In particular, this
means that $h^{1,0}(V_{t})$ defines a holomorphic section of $L$
over $X-\{p\}$ and hence $V_{t}$ is holomorphic on $X-\{p\}-\{h^{1,0}=0\}.$
But since, by definition, $V_{t}$ is locally bounded on $X-\{p\}$
it then follows, by Riemann's extension theorem, that it is\emph{
}holomorphic on all of $X-\{p\}.$ Next, by the decomposition \ref{eq:decomp in pf uniq deg one}
$h^{1,0}=\gamma\otimes s_{p}$ close to $p,$ were $\gamma\in H^{1,0}(X)$
is non-vanishing close to $p.$ Hence, $V_{t}\otimes s_{p}$ defines
an element of $H^{0}(X,TX+L)=H^{0}(X,-K_{X}+L).$ In the case when
$g\geq2$ this latter space is trivial (since $-K_{X}+L$ then has
negative degree) and hence it follows that $V_{t}=0$ on $X-\{p\}$
in this case. Finally, in the case when $g=1$ we have that $TX$
is trivial and hence $H^{0}(X,TX+L)=H^{0}(X,L)=\C s,$ forcing $V=C\frac{\partial}{\partial z},$
where $\frac{\partial}{\partial z}$ denotes an invariant holomorphic
vector field on the torus $X$ and $C$ is constant on $X$. To see
that $C=0$ we use Lemma \ref{lem:lifting} which says that $V$ lifts
to a holomorphic vector field $\tilde{V}$ on $L$ over $X-\{p\}.$
The explicit formula for $\tilde{V}$ in the proof of Lemma \ref{lem:lifting}
now shows that the coefficents of $\tilde{V}$ are locally bounded
(since $V$ and $u$ are smooth on $X)$ and hence it follows from
Riemann's extension theorem that $\tilde{V}$ extends holomorphically
over the fiber of $L$ over $p$ (which is of codimension one in the
total space of $L$). Hence, $\tilde{V}$ is the generator of an automorphism
in $\mbox{Aut}_{0}(X,L).$ But this latter group is trivial over the
torus (as the corresponding automorphism in $\mbox{Aut}_{0}(X)$ must
leave the point $p,$ determined by $L,$ invariant). It follows that
$\tilde{V}$ is tangent to the fibers of $L$ over $X,$ i.e. its
projection $V$ to $X$ vanishes. This hence finishes the proof of
uniqueness. Alternatively, the vanishing of $C$ above can also be
seen directly from the defining equation for $V$ which gives \[
\bar{\partial}(\partial_{t}u)=C|s|_{\psi_{0}}^{2}d\bar{z}\]
 on $X-\{p\}$ and hence on all of $X$ since $\partial_{t}u\in\mathcal{C}^{\infty}(X).$
But integrating over $X$ and using Stokes theorem then forces $C=0$
and hence $\partial_{t}u$ is constant on $X$ proving uniqueness
as before.

\subsubsection{Existence}

First note that the following \emph{coercivity estimate} holds for
some $\epsilon>0$ and some constant $C>0:$ \begin{equation}
\mathcal{F}_{\omega_{0}}(u)\leq-\epsilon J(u),\,\,\, J(u):=\frac{1}{2}\int du\wedge d^{c}u\label{eq:coercivity est}\end{equation}
 for all smooth $u$ on $X.$ To see this it is convenient to use
the formula \ref{eq:l as expect} in the appendix. Estimating the
density of the corresponding probabilitu measure for $\psi_{0}$ from
above and using Fubini's theorem gives a constant $A$ such that \[
-\mathcal{L}_{\omega}(u)=\log\E_{\psi_{0}}(e^{-(u(x_{1})+...+u(x_{n})})\leq\log\int_{X}e^{-u}\omega+A\]
(a variant of this argument appears in \cite{g-s}). Next, we may
assume that $u$ $\int_{X}u\omega_{0}=0.$ By Fontana's generalization
\cite{cher} of the Moser-Trudinger inequality on $S^{2}$ we have
that the left hand side above is bounded from above by $J(u)/2+B$
(i.e. just as in the case when $X=S^{2}$ ). Hence, \[
-\mathcal{L}_{\omega}(u)\leq J(u)-\epsilon J(u)+A+B\]
for $\epsilon=1/2.$ Since, by assumption $Vol(L)=1$ this proves
the coercivity estimate \ref{eq:coercivity est}. Finally, the existence
of a a critical point now follows from basic variational arguments
(just like in \cite{tr}). 
\begin{rem}
\label{rem:existence}For a general complex manifold $X$ the coercivity
estimate (in terms of Aubin's $J-$functional) on the subspace $\mathcal{H}_{\omega_{0}}$
of $\mathcal{C}^{\infty}(X)$ implies the existence of a solution
$u$ of the critical point equation of $\mathcal{F}_{\omega}$ of
finite energy; see remark below for the definition of finite energy
($u$ is smooth when $n=1,$ by elliptic regularity). In fact, it
is enough to assume that $\mathcal{F}_{\omega_{0}}$ is $J_{\omega_{0}}-$proper
\cite{bbgz}. This existence result is proved by a slight modification
of the proof in \cite{bbgz} concerning the case when $L=-K_{X}$
(the coercivity then means that the Fano manifold $X$ is analytically
stable in the sense of Tian). The point is that the functional $\mathcal{L}_{\omega_{0}}(u)$
is continuous on the subspace of all normalized $\omega-$psh functions
such that $J_{\omega_{0}}(u)\leq C$ (for the same reason as in the
case $L=-K_{X}),$ which yields the existence of an extremizer $u$
with finite energy. The fact that $u$ satisfies the critical point
equation then follows from the differentiability theorem \ref{thm:deriv of composed }.
Essentially the same argument was used in the proof of the uniqueness
in Theorem \ref{thm:main} (note that since we are optimizing $\mathcal{F}_{\omega}$
over a convex set with boundary it is a prioiri far from clear that
the extremizer satisfy the critical point equation).
\end{rem}
Finally, a remark about finite energy:
\begin{rem}
\label{rem:finite energy}Consider the setting of Theorem \ref{thm:main}
and assume that there exists a (smooth) critical point, which we may
assume is given by $0.$ Then the inequality furnised by the theorem,
i.e. \[
\mathcal{F}_{\omega_{0}}(u):=\mathcal{E}_{\omega_{0}}(u)-\mathcal{L}_{\omega_{0}}(u)\leq0\]
actually holds for all $u$ in $\mathcal{E}^{1}(X,\omega_{0}),$ i.e.
for al $u$ in the convex set of all $u$ in $\overline{\mathcal{H}}_{\omega_{0}}$
with \emph{finite energy;} $\mathcal{E}(u)>-\infty,$ where \[
\mathcal{E}(u):=\inf_{u'\geq u}\mathcal{E}(u')\]
 when $u'$ ranges over all elements in $\mathcal{H}_{\omega_{0}}$
such that $u'\geq u.$ Equivalently, $\int_{X}(\omega_{u})^{n}=\mbox{Vol}(L)$
and $-\int_{X}u(\omega_{u})^{n}<\infty$ in terms of \emph{non-pluripolar
products} (see \cite{bbgz} and references therein). The inequality
on all of $\mathcal{E}^{1}(X,\omega_{0})$ is simply obtained by writing
$u$ as a decreasing limit of elements in $\mathcal{H}_{\omega_{0}}$
and using the continuity of $\mathcal{E}$ and $\mathcal{L}_{\omega_{0}}$
under such limits \cite{bbgz} (note that $e^{-u}$ is integrable
if $\mathcal{E}(u)>-\infty$ \cite{bbgz}). Moreover, in the case
when $\mbox{Aut}_{0}(X,L)$ is discrete it can be shown that any maximizer
of $\mathcal{F}_{\omega_{0}}$ on $\mathcal{E}^{1}(X,\omega_{0}),$
is in fact equal to a constant. The proof is a simple adaptation of
the argument in \cite{bbgz} concerning the case $L=-K_{X}.$ It would
be interesting to know if the general uniqueness statement in Theorem
\ref{thm:main} also remains true in the larger class $\mathcal{E}^{1}(X,\omega_{0})?$ 
\end{rem}

\subsection{Proofs of Corollary \ref{cor:max of det} and Corollary \ref{cor:max of det on torus}}

First we consider the first corollary where $X$ has genus zero. We
will use the notation from section \ref{sub:Explicit-expression},
i.e. $F=\mathcal{O}(m),$ $K_{X}=\mathcal{O}(-2)$ and $L=\mathcal{O}(m+2).$
Applying the formula in Proposition \ref{pro:f as anal tors on curve}
gives, since$\mbox{Vol }(L)=m+2,$$\log(\frac{\det\Delta_{u}}{\det\Delta_{0}})=$
\[
=\frac{\mbox{deg}(K_{X})}{\mbox{2deg}(L)}J_{\omega_{0}}(u)+N\mathcal{F}_{\omega_{0}}(u)\leq-\frac{1}{(m+2)}\int du\wedge d^{c}u+0\leq0\]
In particular, the left hand side above vanishes precisely when the
gradient of $u$ does, i.e. when $u$ is a constant. This hence finishes
the proof of Corollary \ref{cor:max of det}. 

Finally, when $X$ has genus one Proposition \ref{pro:f as anal tors on curve}
says that $\log(\frac{\det\Delta_{u}}{\det\Delta_{0}})$ is proportional
to $\mathcal{F}_{\omega_{0}}(u)$ and hence Corollary \ref{cor:max of det on torus}
follows from Theorem \ref{thm:main ine one dim}.

\subsection{Proof of Theorem \ref{cor:mean field}}

\emph{The equation $(i):$} Consider the functional \[
\mathcal{F}_{\omega_{0}}^{\delta}(u):=\mathcal{F}_{\omega_{0}}(u)-\delta J(u):=(\mathcal{E}_{\omega_{0}}(u)+\log\int(e^{g_{p}-u})\omega_{0})-\delta J(u)\]
whose critical point equation is precisely equation \ref{eq:intro mean field for u}
with $\mu=4\pi(1+\delta)^{-1}$ under the normalization $\int e^{g_{p}-u}=1.$
Indeed, with our normalizations $dd^{c}u/\omega_{0}=\Delta u/4\pi.$
As shown in the proof of the uniqueness in Theorem \ref{thm:main ine one dim}
(when $L=L_{p})$ the first term $\mathcal{F}_{\omega_{0}}(u)$ above
is strictly concave along geodesics connecting critical points. Moreover,
by the last point in Proposition \ref{pro:e is affine} $J_{\omega_{0}}$
is convex along geodesics and hence $\mathcal{F}_{\omega_{0}}^{\delta}(u)$
is strictly concave for $\delta\geq0,$ i.e. $\mu\leq4\pi$ proving
uniqueness in this case just as before. 

Finally, to prove uniqueness for $\mu\leq4\pi+\epsilon$ it is, by
a standard application of the implicit function theorem, enough to
prove that the linearization of equation \ref{eq:intro mean field for u}
for $\mu=4\pi$ has zero as its unique solution. This follows immediately
from Proposition \ref{pro:non-deg max} below (recall that in the
present case $n=1,\, V=1$ and $\beta_{u}=e^{g_{p}-u}/\int_{X}e^{g_{p}-u})).$

\emph{The equation $(ii):$ }applying the implicit function theorem
as above gives existence for $\mu\in]4\pi-\epsilon,4\pi+\epsilon[.$
Next assume that $\mu\in]0,4\pi[,$ i.e. $0<q<1$ for $q:=\mu/4\pi.$
It is equivalent to prove uniqueness of solutions to the equation
\[
\omega_{0}+dd^{c}u=(|s|_{\psi_{0}}^{2}e^{-u})^{q}\omega_{0}\]
which describe the critical points of \[
\mathcal{F}_{q}(u):=\mathcal{E}_{\omega_{0}}(u)+\frac{1}{q}\log\int(|s|_{\psi_{0}}^{2}e^{-u})^{q}\omega_{0}\]
which are normalized in the sense that $\int(|s|_{\psi_{0}}^{2}e^{-u})^{q}=1,$
which will henceforth be assumed. 

To this end we let (for $q$ fixed) $w:=qu+(1-q)g_{p}$ so that $(|s|_{\psi_{0}}^{2}e^{-u})^{q}=|s|_{\psi_{0}}^{2}e^{-w}.$
Note that $w\in\overline{\mathcal{H}_{\omega_{0}}},$ i.e. $\psi_{0}+w$
is the weight of a singular metric on $L$ whose curvature form $\omega_{w}$
is positive in the sense of currents. By a standard regurization argument
we may write $g_{p}$ as a decreasing limit of $g^{(\epsilon)}\in\mathcal{H}_{\omega_{0}}$
and let $w^{(\epsilon)}=qu+(1-q)g^{(\epsilon)}.$ Hence, \[
\log\int(|s|_{\psi_{0}}^{2}e^{-u})^{q}\omega_{0}=\mathcal{L}_{\omega_{0}}(w)=\lim_{\epsilon\rightarrow0}\mathcal{L}_{\omega_{0}}(w_{\epsilon})\]
and then it follows as before that $\mathcal{F}_{q}(u_{t})$ is concave
along any $\mathcal{C}^{0}-$geodesic $u_{t}$ (since $\mathcal{E}_{\omega_{0}}(u_{t})$
is affine and $\log\int(|s|_{\psi_{0}}^{2}e^{-u_{t}})^{q}\omega_{0}$
is a limit of concave functions). To prove uniqueness it will, as
before, be enough to prove that $sV_{u}$ defines an element of $H^{0}(X,L)$
where $V_{u_{t}}$ is the vector field associated to a geodesic of
critical points $u_{t}$ as before (since $TX$ is trivial we will
identify $V_{u}$ with a function as before). Moreover, since we have
assumed that $q<1$ it will be enough to prove the following \begin{equation}
\mbox{Claim:\,\,\,}\bar{\partial}V_{u}=0\,\mbox{on\ensuremath{\, X-\{p\}}}\label{eq:claim mean d}\end{equation}
 Indeed, \[
|V_{u}s|_{\psi_{0}}^{2}=|\frac{\bar{\partial}(\partial_{t}u)}{(|s|_{\psi_{0}}^{2}e^{-u})^{q}}|^{2}|s|_{\psi_{0}}^{2}\leq C\frac{1}{(|s|_{\psi_{0}}^{2q-1})^{2}}\]
which is clearly in $L^{2}(X-\{p\})$ when $0<q<1$ and hence, if
the claim above holds, $V_{u}s$ extends over $\{p\}$ to a global
holomorphic section of $L,$ by basic local properties of holomorphic
functions. To prove the previous claim first observe that explicit
differentiation gives (this is essentially a simple special case of
the calculation appearing in the appendix taking advantage of the
fact that $H^{0}(X,L)=\C s):$ \[
\partial_{t}\partial_{\bar{t}}\mathcal{L}_{\omega_{0}}(w_{t})=\int\beta(w_{t})((\partial_{t}\partial_{t}qu_{t})-\int|s\partial_{t}qu_{t}-C_{t}|_{\psi_{0}}^{2}e^{-w_{t}}\omega_{0})\]
 where $\beta(w_{t})=(|s|_{\psi_{0}}^{2}e^{-u_{t}})^{q}$ and the
constant $C_{t}$ is uniquely determined by requiring that \[
\left\langle s(\partial_{t}qu_{t}-C_{t}),s\right\rangle _{\psi_{0}+w_{t}}=0\]
 i.e. that \[
\alpha_{t}:=s(q\partial_{t}u_{t}-C_{t})\]
 be orthogonal to $H^{0}(X,L).$ In other words, $\alpha_{t}$ is
the $L^{2}(e^{-(\psi_{0}+w_{t}})$ minimal solution of the $\bar{\partial}$-equation
\ref{eq:dbar equation appendix}, with $s=s_{i},$ in the appendix.
Now by a standard regularization argument (see for example below)
\begin{equation}
\int_{X}|\alpha|_{\psi_{0}}^{2}e^{-w}\omega_{0}\leq\int_{X}|\bar{\partial}(\partial_{t}u_{t})|_{\omega_{qu}}^{2}|s|_{\psi_{0}}^{2}e^{-w}\omega_{0}<\infty,\label{eq:hormand estimate sing case}\end{equation}
(note that $\omega_{qu}$ is the part of the current $\omega_{w}$
which is absolute continuous wrt $\omega_{0}).$ All in all this means
that \[
\partial_{t}\partial_{\bar{t}}\mathcal{L}_{\omega_{0}}(w_{t})_{t=0}\geq\int\beta(w)c(qu)=0\]
 where $c(u)$ is the non-negative function appearing in equation
\ref{eq:geod eq concrete} (which in fact vanishes since $u_{t}$
is a $\mathcal{C}^{0}-$geodesic). Moreover, equality holds above
precisely when equality holds in the estimate \ref{eq:hormand estimate sing case}.
As in the previous case of smooth weights we will prove that the latter
equality implies the claim above. We will proceed by regularization:
we replace $w$ with $w^{(\epsilon)}$ and accordingly for all objects
defined in terms of $w^{(\epsilon)}.$ Then from the explicit definition
of $\alpha^{(\epsilon)}$ it follows directly that \[
\left\Vert \alpha^{(\epsilon)}\right\Vert _{w^{(\epsilon)}}\rightarrow\left\Vert \alpha\right\Vert _{w},\,\,\,\liminf_{\epsilon}\left\Vert \bar{\partial}\alpha^{(\epsilon)}\right\Vert _{w^{(\epsilon)},\omega_{w^{(\epsilon)}}}\leq\left\Vert \bar{\partial}\alpha\right\Vert _{w,\omega_{qu}}\]
where a subscript of the form $w,\omega$ indicates that we are taking
the point-wise norm along the fibers of $L$ and $T^{0,1}X$ wrt the
metrics $e^{-(\psi_{0}+w)}$ and $\omega,$ respectively. But since
\[
\left\Vert \alpha^{(\epsilon)}\right\Vert _{w^{(\epsilon)}}\leq\left\Vert \bar{\partial}\alpha^{(\epsilon)}\right\Vert _{w^{(\epsilon)},\omega_{w^{(\epsilon)}}},\,\,\,\left\Vert \alpha\right\Vert _{w}=\left\Vert \bar{\partial}\alpha\right\Vert _{w,\omega_{qu}}\]
 it then follows that \[
\left\Vert \alpha^{(\epsilon)}\right\Vert _{w^{(\epsilon)}}=(1+o(1))\left\Vert \bar{\partial}\alpha^{(\epsilon)}\right\Vert _{w^{(\epsilon)},\omega_{w^{(\epsilon)}}},\]
 where $o(1)$ denotes a sequence tending to zero when $\epsilon\rightarrow0.$
Repeating the arguments in the end of the proof of Proposition \ref{pro:horm append}
in the appendix it then follows readily that \[
\left\Vert \bar{\partial}(sV_{w^{(\epsilon)}})\right\Vert _{w^{(\epsilon)},\omega_{w^{(\epsilon)}}}\rightarrow0\]
and hence $V_{u}$ is holomorphic on $X-\{p\}$ (since we may assume
that $w^{(\epsilon)}\rightarrow w$ uniformly with all derivatives
on compacts of $X-\{p\}$ and $V_{u}=V_{w}$ there). This proves the
claim \ref{eq:claim mean d} and hence finishes the proof of the proposition.

The following proposition was used in the previous proof:
\begin{prop}
\label{pro:non-deg max}Assume that either the assumptions in Theorem
\ref{thm:main} hold and that $\mbox{Aut\ensuremath{_{0}(X)}}$is
trivial or that $X$ is a complex curve of genus at least one. Then
the null-space of the linearization $A_{0}$ of the operator \[
u\mapsto(d\mathcal{F}_{\omega_{0}})_{u}:=\frac{1}{Vn!}\omega_{u}^{n}-\beta_{u}\]
at the unique critial point $u_{0}$ of $\mathcal{F}_{\omega_{0}}$
consists of the\emph{ constant }functions.\end{prop}
\begin{proof}
At least formally, (i.e. if $\mathcal{H}_{\omega_{0}}$ were finite
dimensional) this would follow from the strict convexity of $\mathcal{F}_{\omega_{0}},$
modulo constants, along smooth geodesics in $\mathcal{H}_{\omega_{0}}.$
Next, we will provide a rigorous proof which avoids the use of geodesics.
To this end fix $v\in\mathcal{C}^{\infty}(X)$ and observe that \[
\int_{X}A_{0}[v]v=\frac{\partial^{2}\mathcal{F}_{\omega_{0}}(u_{t})}{\partial^{2}t},\,\,\, u_{t}=u_{0}+tv\]
 at $t=0,$ only using that $t=0$ is a critical point of $\mathcal{F}_{\omega_{0}}(u_{t})$.
Assume that $A_{0}v=0.$ Recall that $\mathcal{F}_{\omega_{0}}=\mathcal{E}_{\omega_{0}}-\mathcal{L}_{\omega_{0}}$
and note that, by formula \ref{eq:genaral second deriv} $\frac{\partial^{2}\mathcal{E}_{\omega_{0}}(u_{t})}{\partial^{2}t}$
coincides with the right hand side in formula \ref{eq:pf of convexity of l}
at $t=0$ since $u_{0}$ satisfies the critical point equation. As
a consequence the argument leading up to the inequality \ref{eq:pf of convexity of l}
shows that equality holds in Hörmanders's $L^{2}-$estimate. But then
it follows from the previous arguments (for $n=1$ see section \ref{sub:Proof-of-v vanishes})
that the corresponding vector field $V_{0}$ associated to $u_{t}$
at $t=0$ vanishes on $X,$ which by definition is equivalent to $v$
beeing constant. This finishes the proof of the uniqueness result.\end{proof}
\begin{rem}
\label{rem:family}The argument using the implicit function theorem
above also shows the following fact of independent interest: if $\mathcal{L}$
is a relatively ample line bundle over a holomorphic submersion $\mathcal{X}\rightarrow S$
such that the fibers have no holomorphic vector fields and the functional
$\mathcal{F}_{\omega_{s_{0}}}$ defined over a given fiber $X:=\mathcal{X}_{s_{0}}$
has a critical point, then there also exist critical points for $s$
close to $s_{0}.$ In the case when $\mathcal{L}$ is the relative
anti-canonical line bundle, i.e. the critical points are Kähler-Einstein
metrics, this is well-known (a similar argument appears in Aubin's
continuity method \cite{au,ti} on Fano manifolds). 
\end{rem}

\subsection{\label{sub:Alternative-proof-of}Alternative proof of regularity
and uniqueness}

In this section we will show how to prove the {}``uniqueness'' in
Theorem \ref{thm:main} only using the regularity of the geodesics
furnished by Theorem \ref{thm:berman-dem} and the theory of fully
non-linear elliptic operators in $n$ complex dimensions (applied
to the Monge-Ampère operator on $X$ as in \cite{bl0}). In particular,
this latter theory amounts to the basic \emph{linear} elliptic estimates
for the Laplacian when $n=1$ which will allow us to handle geodesics
connecting degenerate metrics in section \ref{sub:The-one-dimensional-case}.

We will denote by $W^{r,p}(X)$ the Sobolev space of all distributions
$f$ on $X$ such that $f$ and the local derivatives of total order
$r$ are in $W^{0,p}(X):=L^{p}(X)$ (equvalently, all local derivatives
of total order $\leq r$ are in $L^{p}(X));$ see \cite{au2}. If
$f$ is function on $M=[0,1]\times X$ we will write \emph{$f_{t}\in W^{r,p}(X)$
uniformly wrt $t$} if the corresponding Sobolev norms on $(X,\omega_{0})$
of $f_{t}$ are uniformly bounded in $t.$ We will also use the following
basic facts repeatedly:
\begin{itemize}
\item If $f$ is a function on $M$ such that \emph{$f_{t}\in W^{r,p}(X)$
uniformly wrt $t,$} then the distribution $f$ is in $W^{r,p}(\dot{M})$
and the corresponding Sobolev norms on $\dot{M}$ are bounded.
\item Partial derivatives of distributions commute
\item If $f,g\in W^{1,p}(X)$ for any $p>1.$ Then $fg\in W^{1,p}(X)$ for
any $p>1$ and Leibniz product rule holds for the distributional derivatives.
\end{itemize}
Note that as in section \ref{sub:Proof-of-Theorem-main} it will be
enough to prove that the geodesic $u_{t}$ is smooth wrt $(t,x)$
in the \emph{interiour} of $M.$ However, the arguments below will
even give uniform estimates on the local Sobolev norms up to the boundary
of $M.$

Assume now that the boundy data $u_{0}$ and $u_{1},$ defining the
geodesic $u_{t}$ are in $\mathcal{C}^{1,1}(X).$ Since $u_{t}$ is
convex in $t$ the right derivative (or tangent vector) $v_{t}(x):=\frac{d}{dt}_{+}u_{t}$
exists for all $(t,x).$
\begin{lem}
\emph{\label{lem:v is bded}The right tangent vector} \emph{$v_{t}$
of $u_{t}$ at $t$ is uniformly bounded on $M.$} \end{lem}
\begin{proof}
First observe that by the convexity in $t$\[
u_{t}-u_{0}\leq t(u_{1}-u_{0})\leq C_{1}t,\]
 using that $u_{0}$ and $u_{1}$ are continuos and hence uniformly
bouned on $X$ in the last step. Hence, $v_{t}\leq C.$ To get a lower
bound first observe that there is a {}``psh extension'' $\tilde{u}_{t}$
which is uniformly Lipshitz. Indeed, just take $\tilde{u}_{t}:=(1-t)u_{0}+tu_{1}+Ae^{t}$
for $A>>1.$ Using that $0\leq dd^{c}u_{0},\, dd^{c}u_{1}\leq C\omega_{0}$
it is straight-forward to check that $dd^{c}\tilde{U}+\pi^{*}\omega_{0}\geq0$
on $M$ for $A$ sufficently large. Since $U$ is defined by the upper
envelope \ref{eq:envelop} it follows that $\tilde{u}_{t}\leq u_{t}$
and hence \[
u_{t}-u_{0}\geq\tilde{u}_{t}-\tilde{u}_{0}\geq C_{2}t.\]
giving $v_{0}\geq C_{2}.$ Finally, by convexity we get $C_{2}\leq v_{0}\leq v_{t}\leq C_{1}$
which proves the lemma.\end{proof}
\begin{prop}
\label{pro:unique alt}Let $u_{0}$ be a critical point of $\mathcal{F}_{\omega_{0}}$
on $\overline{\mathcal{H}}_{\omega_{0}}\cap\mathcal{C}^{1,1}(X),$
$u_{1}$ an arbitrary element in $\overline{\mathcal{H}}_{\omega_{0}}\cap\mathcal{C}^{1,1}(X)$
and $u_{t}$ the geodesic connecting $u_{0}$ and $u_{1}.$ If $\mathcal{L}_{\omega_{0}}(u_{t})$
is affine, then there is an automorphism $S_{1}$ of $(X,L)$, homotopic
to the identity, such that $u_{1}-u_{0}=S_{1}^{*}\psi_{0}-\psi_{0}.$\end{prop}
\begin{proof}
\emph{Step 1: $u_{t}\in\mathcal{C}^{\infty}(X).$} First note that
by Theorem \ref{thm:berman-dem} $u_{t}\in\mathcal{C}_{\C}^{1,1}(X).$
Moreover, as shown in the beginning of section \ref{sub:Proof-of-Theorem-main}
it follows under the assumptions above that, for any $t,$ the function
$u_{t}$ satisfies the Euler-Lagrange equations \ref{eq:e-l for f}
on $X.$ Hence, just as in section \ref{sub:Proof-of-Theorem-main}
Blocki's complex version of the regularity result of Trudinger, now
applied to local patches of $\{t\}\times X$ immediately gives that
$u_{t}\in\mathcal{C}^{\infty}(X)$ (when $n=1$ this follows from
basic linear elliptic theory).

\emph{Step 2: $\Delta_{X}v_{t}\in L^{\infty}(X)$ uniformly wrt $t.$}
Differentiating the Euler-Lagrange equation wrt $t$ from the right
gives \begin{equation}
ndd^{c}v_{t}\wedge(\omega_{t})^{n-1}=\frac{d\beta_{u_{t}}(x)}{dt}_{+}=:R[v_{t}],\label{eq:linearized euler-l}\end{equation}
in the sense of currents. Of course, this would follow immediately
from the chain rule if $u_{t}$ were smooth in $(t,x)$ In the present
case the right hand side is handled as in Lemma \ref{lem:deriv of beta}
in the appendix. As for the left hand side it is (a part from the
trivial case $n=1)$ handled precisely as in Lemma in the proof of
formula 6.11 in \cite{bbgz}. Moreover, lemma \ref{lem:deriv of beta}
in the appendix implies the bound \begin{equation}
\left\Vert R[v_{t}]/(\omega_{0})^{n}\right\Vert _{L^{\infty}(X)}\leq C\left\Vert v_{t}\right\Vert _{L^{\infty}(X)}\label{eq:ineq for r}\end{equation}
To see this, just note that \[
R[v]\leq2\left\Vert v\right\Vert _{L^{\infty}(X)}\int_{X}\left|K(x,y)\right|^{2}e^{-(\psi(x)+\psi(y))}=2\left\Vert v\right\Vert _{L^{\infty}(X)}\beta_{u},\]
 using the well-known {}``reproducing property'' of the Bergman
kernel (formula \ref{eq:integr out} in the appendix). By formula
\ref{eq:bounds on beta} in the appendix this proves the inequality
\ref{eq:ineq for r}.

Now, since $\omega_{t}>\delta\omega_{0},$ formula \ref{eq:linearized euler-l}
gives that the distribution $\Delta_{\omega_{t}}v_{t},$ where $\Delta_{\omega_{t}}$
is the Laplacian on $X$ wrt the metric $\omega_{t}:=\omega_{u_{t}},$
is in $L^{\infty}(X)$ uniformly wrt $t$ and \[
\left\Vert \Delta_{\omega_{t}}v_{t}\right\Vert _{L^{\infty}(X)}\leq C\left\Vert v_{t}\right\Vert _{L^{\infty}(X)}\leq C',\]
 by lemma \ref{lem:v is bded}. 

\emph{Step 3: $\Delta_{M}u\in W^{1,p}(M)$ for any $p\geq1.$} First
observe that by step 1 \begin{equation}
\partial_{z}(\partial_{z_{i}}\partial_{\bar{z}_{j}}u)\in L^{\infty}(X),\label{eq:zzz deriv}\end{equation}
 uniformly wrt $t.$ Also note that \begin{equation}
\partial_{t}(\partial_{z_{i}}\partial_{\bar{z}_{j}}u)\in L^{\infty}(X),\label{eq:tzz}\end{equation}
uniformly wrt $t.$ Indeed, $\partial_{t}(\partial_{z_{i}}\partial_{\bar{z}_{j}}u)=\partial_{z_{i}}(\partial_{\bar{z_{j}}}\partial_{t}u)=(\partial_{z_{i}}\partial_{\bar{z}_{j}})v_{t}\in L^{p}(X),$
uniformly wrt $t$\textbf{,} for any $p>1,$ by step 2 and local elliptic
estimates for $\Delta_{X}.$ Next, we will use that the following
identity proved in lemma \ref{lem:normal-normal} below:\begin{equation}
\partial_{t}\partial_{\bar{t}}u=\left|V_{t}\right|_{\omega_{t}}^{2}=\left|\partial_{\bar{z}}v_{t}\right|_{\omega_{t}}^{2},\label{eq:ma eq in prop unique}\end{equation}
 where $\left|V_{t}\right|_{\omega_{t}}^{2}$ denotes the point-wise
norm of $V_{t}$ wrt the metric $\omega_{t}$ (where we have used
that $\omega_{t}>0)$. First we have \begin{equation}
\partial_{z}(\partial_{t}\partial_{\bar{t}}u)=\partial_{z}\left|\partial_{\bar{z}}v_{t}\right|_{\omega_{t}}^{2}\in L^{p}(X),\label{eq:ztt}\end{equation}
 uniformly wrt $t,$ for any $p>1$ using Step 1 and Step 2 combined
with local elliptic estimates on $X$ for $\Delta_{X}.$ Next, \begin{equation}
\partial_{t}(\partial_{t}\partial_{\bar{t}}u)\in L^{p}(X),\label{eq:trippel deriv wrt t}\end{equation}
 uniformly wrt $t.$ Indeed, $\partial_{t}(\partial_{t}\partial_{\bar{t}}u)=\partial_{t}\left|\partial_{\bar{z}}\partial_{t}u\right|_{\omega_{t}}^{2}$
and since locally $\partial_{t}\omega_{t}=\partial_{t}(\partial_{z}\partial_{\bar{z}}u)$
\ref{eq:trippel deriv wrt t} follows from \ref{eq:tzz} and \ref{eq:ztt}
combined with Leibniz product rule. All in all this proves Step 3. 

Now by Step 3 and elliptic estimates for the Laplacian we have $u\in W^{3,p}(M).$
In particular, $u$ is locally in $\mathcal{C}^{2}(M).$ As a consequence
the proof of Theorem 2.6 in \cite{bern2} immediately gives that $V_{t}$
is a holomorphic vector field on $X$ for any $t.$ Finally, we will
recall a slight variant of the argument in \cite{bern2} which shows
that $\partial_{\bar{t}}V_{t}=0$ for $V_{t}$ seen as a distribution
on the interiour of $M.$ To simplify the notation we assume that
$n=1,$ but modulo the change to matrix notation the case $n>1$ is
the same. First we write \ref{eq:int multip} in the form \begin{equation}
\omega V=\partial_{\bar{z}}\partial_{t}u,\label{eq:int mult local}\end{equation}
where we have identified $V$ and $\omega$ with elements in $L^{p}(M)$
for $p>>1.$ By Leibniz rule \[
\partial_{\bar{t}}(V\omega)=(\partial_{\bar{t}}V)\omega+V(\partial_{\bar{t}}\omega)\]
Next, observe that \[
\partial_{\bar{t}}\omega=\partial_{\bar{t}}(\partial_{z}\partial_{\bar{z}}u)=\partial_{\bar{z}}(\partial_{\bar{t}}\partial_{z}u)=\partial_{\bar{z}}(\omega\bar{V}),\]
 using \ref{eq:int mult local} in the last step. Hence, since, as
shown above, $\partial_{\bar{z}}V=0,$ the two previous equations
together give \[
\partial_{\bar{t}}(V\omega)=(\partial_{\bar{t}}V)\omega+\partial_{\bar{z}}(V\omega\bar{V})=\partial_{\bar{z}}(\partial_{\bar{t}}\partial_{t}u),\]
 also using \ref{eq:int mult local} in the last step and commuting
$\partial_{\bar{z}}$ and $\partial_{\bar{t}}.$ Since, $V\omega\bar{V}=\left|V\right|_{\omega}^{2}$
it follows by \ref{eq:ma eq in prop unique} that $(\partial_{\bar{t}}V)\omega=0.$
But since, $\omega>0$ and $(\partial_{\bar{t}}V)$ is in $L^{p}(M)$
for all $p>1$ this forces $(\partial_{\bar{t}}V)=0$ a.e. on $M.$
In particular, $(\partial_{\bar{t}}V)=0$ as a distribution on $M.$
Hence, it follows that the distribution $V_{t}$ is in the null-space
of the $\bar{\partial}-$operator on $M.$ By local elliptic theory
it follows that $V_{t}$ is smooth and hence holomorphic in the interiour
of $M.$ Finally, the automorphism $S_{1}$ is obtained precisely
as in the end of section \ref{sub:Proof-of-Theorem-main}.
\end{proof}
In the previous proof we used the following
\begin{lem}
\emph{\label{lem:normal-normal}Under the assumptions in the previous
proposition the following holds: $\partial_{t}\partial_{\bar{t}}u\in L^{\infty}(X)$}
\textup{uniformly in $t$} \emph{and \[
\partial_{t}\partial_{\bar{t}}u=\left|\bar{\partial}_{X}\partial_{t}u\right|_{\omega_{u_{t}}}^{2}.\]
}\end{lem}
\begin{proof}
By assumption the Monge-Ampère measure $(dd^{c}U+\pi_{X}^{*}\omega_{0})^{n+1}$
vanishes on $M.$ Moreover, by Step 1 in the proposition above $\Delta_{X}u_{t}\in C^{\infty}(X)$
for any $t$ with bounds on the Sobolev norms which are uniform wrt
$t.$ Combining this latter fact with lemma \ref{lem:v is bded} gives
that $U$ is Lipschitz on $M.$ Finally, as shown in Step 2 in the
proof of proposition above $\Delta_{X}\partial_{t}u_{t}\in L^{\infty}(X)$
uniformly wrt $t.$ We will next show that these properties are enough
to prove the lemma. As the statement is local we may as well consider
the restriction of $u:=U$ to an open set biholomorphic to a domain
in $\C^{n+1}=\C_{t}\times\C_{z}^{n}.$ Denote by $u^{\epsilon}$ the
local smooth function obtain as the convolution of $u$ with a fixed
local compactly supported smooth family of approximations of the identity.
Expanding gives \begin{equation}
(dd^{c}U+\pi_{X}^{*}\omega_{0})^{n+1}=(\partial_{t}\partial_{\bar{t}}u^{\epsilon}-\left|\partial_{\bar{z}}\partial_{t}u^{\epsilon}\right|_{\omega_{u^{\epsilon}}}^{2})(\omega_{u^{\epsilon}})^{n}\wedge dt\wedge d\bar{t}.\label{eq:expanding monge as c}\end{equation}
 Now since, by assumption, $\left|\partial_{\bar{z}}\partial_{t}u^{\epsilon}\right|_{\omega_{u^{\epsilon}}}^{2}\leq C$
the second term tends to $\left|\partial_{\bar{z}}\partial_{t}u\right|_{\omega_{u}}^{2})(\omega_{u})^{n}\wedge dt\wedge d\bar{t}$
weakly when $\epsilon\rightarrow0.$ Moreover, by assumption $u^{\epsilon}\rightarrow u$
uniformly locally and since the Monge-Ampère operator is continous,
as a measure, under uniform limits of psh functions \cite{de3} it
will now be enough to prove that \begin{equation}
(\partial_{t}\partial_{\bar{t}}u^{\epsilon})(\omega_{u^{\epsilon}})^{n}\wedge dt\wedge d\bar{t}\rightarrow(\partial_{t}\partial_{\bar{t}}u)(\omega_{u})^{n}\wedge dt\wedge d\bar{t}\label{eq:weak conv in pf lemma n-n}\end{equation}
 weakly, where the right hand sice is well-defined since $\partial_{t}\partial_{\bar{t}}u_{t}$
defines a positive measure on $\C^{n+1}$ and $(\omega_{u_{t}})^{n}/\omega_{0}^{n}$
is continous on $\C^{n+1}.$ To this end fix a test funtion $f$ i.e.
a smooth and compactly supported function on $\C^{n+1}.$ Then, with
$\int$ denoting the integral over $\C^{n+1},$ \[
\int f(\omega_{u^{\epsilon}})^{n}(\partial_{t}\partial_{\bar{t}}u^{\epsilon})\wedge dt\wedge d\bar{t}=:\int g_{\epsilon}(\partial_{t}\partial_{\bar{t}}u^{\epsilon})=-\int(\partial_{t}g_{\epsilon})(\partial_{\bar{t}}u^{\epsilon})\]
By assumption $(\partial_{t}g_{\epsilon})$ and $(\partial_{\bar{t}}u^{\epsilon})$
tend to $(\partial_{\bar{t}}u)$ and $(\partial_{\bar{t}}u),$ respectively
in $L^{p}(X)$ for any $p>1,$ uniformly wrt $t$ (more precisily
by the assumption on $\Delta_{X}u_{t}$ and the fact that $u$ is
Lipschitz). Hence, by Hölders's inequality \[
\int g_{\epsilon}(\partial_{t}\partial_{\bar{t}}u^{\epsilon})\rightarrow-\int(\partial_{t}g)(\partial_{\bar{t}}u).\]
 Finally, since $(\partial_{t}g)\in L^{\infty}(X)$ uniformly wrt
$t$ (by the assumption on $\Delta_{X}u_{t}$) and since $\partial_{t}\partial_{\bar{t}}u$
defines a positive measure, Leibniz rule combined with the dominated
convergence theorem gives (by a simple argument using a regularization
of $g)$ \[
-\int(\partial_{t}g)(\partial_{\bar{t}}u)=\int g(\partial_{t}\partial_{\bar{t}}u)\]
This proves \ref{eq:weak conv in pf lemma n-n} and hence finishes
the proof of the lemma. 
\end{proof}

\subsection{\label{sub:The-one-dimensional-case}The one-dimensional case with
degenerate boundary data}

In this section we will show that in the one-dimensional case the
argument in the previous section goes through without assuming strict
positivity of the curvature forms of the boundary data. This fact
was used in the proof of Theorem \ref{thm:main ine one dim}. 
\begin{prop}
Assume that $\dim X$=1 and let $u_{0}$ and $u_{1}$ be two smooth
critical points of the functional $\mathcal{F}_{\omega_{0}}.$ Then
the geodesic $u_{t}$ connectiong $u_{0}$ and $u_{1}$ satisfies
the following regularity properties: 
\begin{itemize}
\item $u_{t}$ and in particular $\omega_{t}$ are smooth on $X$ for each
fixed $t$
\item $\partial_{X}u_{t}$ is smooth on $X$ 
\item and $(\partial_{t}\partial_{\bar{t}}u)(\frac{\omega_{t}}{\omega_{0}})$
is locally bounded in the interiour of $M$ and $U(t,x):=u_{t}(x)$
is locally in $\mathcal{C}_{\C}^{1,1}$ in the interiour of $\{(t,x)\in M:\,\omega_{t}>0\}$
\end{itemize}
\end{prop}
\begin{proof}
The key point is that Step1 and Step 2 in the proof in the previous
section are still valid when $n=1$ (even if $L+K_{X}$ is not globally
generated). Indeed, as shown in the beginning of section \ref{sub:Proof-of-Theorem main one dim}
$u_{t}$ is continuous and satisfies the following non-linear elliptic
equation on $X$ \[
\Delta_{X}u_{t}+1=fe^{-u_{t}}\]
where $f$ is a non-negative smooth function and \textbf{$\Delta_{X}$
}denotes the (normalized) Laplacian on $X$ wrt $\omega_{0}.$ In
particular, $\Delta_{X}u_{t}\in L^{p}(M)$ for all $p>1.$ Repeating
the argument and using the standard elliptic estimates for $\Delta_{X}$
finally proves that $u_{t}\in\mathcal{C}^{\infty}(X)$ (with uniform
bounds wrt $t$ on the derivatives). This proves the first point in
the proposition.

As for Step 2, differentiating the equation above wrt $t$ also shows
that $v_{t}\in\mathcal{C}^{\infty}(X)$ (uniformly wrt $t)$, the
point beeing that the linearized equation does not involve $\omega_{t}$
when $n=1.$ This proves the second point in the proposition. 

Finally, to prove the last point we may apply a slight variant of
Lemma \ref{lem:normal-normal} obtained by multiplying formula formula
\ref{eq:ma eq in prop unique} with the function $(\frac{\omega_{t}}{\omega_{0}})$
and using the previous points.\end{proof}
\begin{rem}
When $n=1$ any critical point with finite Dirichlet energy is in
fact smooth. Indeed, by the Moser-Trudinger inequality $e^{-u_{t}}$
in in $L^{p}(X)$ for any $p$ and hence the previous boot strapping
argument still applies.
\end{rem}

\section{\label{sec:Convergence-towards-Mabuchi's}The large tensor power
limit: analytic torsions and Mabuchi's K-energy energy}

In this section we will consider the asympotic situation when the
ample line bundle $L$ is replaced by a multiple $kL$ for a large
positive integer $k$. Building on \cite{bern2} Berndtsson we will
relate the large $k$ asymptotics of $\mathcal{F}_{k\omega_{0}}$
to \emph{Mabuchi's K-energy.} The work \cite{bern2} was in turn inspired
by the seminal work of Donaldson \cite{don1} where a functional closely
related to $\mathcal{F}_{k\omega_{0}}$ was introduced (see section
\ref{sub:Comparison-with-Donaldson's} below). We will give a new
simple proof of the lower bound on Mabuchi's K-energy (Theorem \ref{thm:mabuchi}
below) which only uses the $C^{0}-$regularity of the geodesic connecting
two given smooth points in $\mathcal{H}_{\omega_{0}}.$ See \cite{ch}
for a proof which uses the $\mathcal{C}_{\C}^{1.,1}-$regularity (Theorem
\cite{ch}) in the case when the first Chern class of $X$ is assumed
non-positive. We will also consider the analytic torsion for large
tensor powers over a Fano manifold and give the proof of Theorem \ref{thm:fano torsion}.
Finally, some conjectures will be proposed.

Fixing $\omega_{0}\in c_{1}(L)$ we willl take $k\omega_{0}$ as the
reference Kähler metric in $c_{1}(kL).$ Throughout the section $u$
will denote an element in $\mathcal{H}_{\omega_{0}}.$We will make
use of the following essentially well-known asymptotics for the differential
of the functional $\mathcal{L}_{k\omega_{0}}:$ 

\begin{equation}
(d\mathcal{L}_{k\omega_{0}})_{ku}=\frac{k}{\mbox{Vol}(L)}(\frac{\omega_{u}^{n}}{n!}-\frac{1}{k}(\frac{1}{2}\mbox{(Ric}\omega_{u}-S\omega_{u})\wedge\omega_{u}^{n-1}/(n-1)!)+O(1/k))/n!\label{eq:asym of dl}\end{equation}
 where $S$ is certain invariant associated to $(X,L).$ Using formula
\ref{eq:deriv of l} the proof of the previous formula is reduced
to the well-known asymptotics of the Bergman measure on $kL+E,$ where
$E$ is a given line bundle on $E,$ due to Tian-Catlin-Zelditch.
Here we take $E=K_{X}$ (see \cite{bern2}). In particular, we obtain
\begin{equation}
(d\mathcal{F}_{k})_{ku}:=\frac{1}{\mbox{Vol}(L)}\frac{1}{2}(\mbox{Ric}\omega_{u}-\omega_{u})\wedge\frac{\omega_{u}^{n-1}}{(n-1)!}+o(1)\label{eq:asympt of df}\end{equation}
Following Mabuchi \cite{m2,ti} the \emph{K-energy} (also called the
\emph{Mabuchi functional}) is defined, up to an additive constant,
by letting $-\mathcal{M}$ be the primitive on $\mathcal{H}_{\omega_{0}}$
of the exact one-form defined as the leading term above, i.e. as $\omega_{u}^{n}$
times the scalar curvature minus its average. Hence, $u$ is a critical
point of $\mathcal{M}$ on $\mathcal{H}_{\omega_{0}}$ iff the Kähler
metric $\omega_{u}$ has \emph{constant scalar curvature}. We will
denote by $\mathcal{M}_{\omega_{0}}$ the K-energy normalized so that
$\mathcal{M}_{\omega_{0}}(0)=0.$ Integrating along line segments
in $\mathcal{H}_{\omega_{0}}$ and using \ref{eq:asympt of df} immediately
gives the asympotics \begin{equation}
\mathcal{F}_{k}(u)=-\mathcal{M}_{\omega_{0}}(u)+o(1).\label{eq:asympt of f}\end{equation}
 For the most general version of the following theorem see \cite{c-t}. 
\begin{thm}
\label{thm:mabuchi}Assume that the Kähler metric $\omega_{u_{0}}$
has constant scalar curvature. Then $u_{0}$ minimizes Mabuchi's K-energy
$\mathcal{M}_{\omega_{0}}$ on $\mathcal{H}_{\omega_{0}}.$\end{thm}
\begin{proof}
By the cocycle property of $\mathcal{M}_{\omega_{0}}$ we may as well
assume that $u=0$ in the statement above. Now fix an arbitrary $u$
in $\mathcal{H}_{\omega_{0}}$ and take the $C^{0}-$geodesics $u_{t}$
connecting $0$ and $u.$ Given a positive integer $k$ the fact that
$\mathcal{F}_{k}$ is concave along $u_{t}$ (compare the proof of
Theorem \ref{thm:main}) immediately gives \[
\mathcal{F}_{k}(u)\leq\mathcal{F}_{k}(0)+\frac{d}{dt}_{t=0+}\mathcal{F}_{k}(u_{t})\]
Combining formulas \ref{eq:asympt of f}, \ref{eq:asympt of df} then
gives \begin{equation}
\mathcal{F}_{k}(u)\leq\mathcal{F}_{k}(0)+\frac{1}{\mbox{Vol}(L)}\frac{1}{2}\int(\omega_{0}-\mbox{Ric}\omega_{0})\wedge\frac{\omega_{0}^{n-1}}{(n-1)!})+O(k^{-1}))\omega_{0}^{n})(-v_{0}),\label{eq:fk in terms of deriv-1}\end{equation}
 where $v_{0}=\frac{du}{dt}_{t=0+}$. By the assumption on $\omega_{0}$
the integral above vanishes. Moreover, by lemma \ref{lem:v is bded}
we have that $v_{0}$ is uniformly bounded (in fact it is enough to
know that its $L^{1}-$norm is uniformly bounded, which can be proved
as in the proof of Theorem \ref{thm:fano torsion} below. Letting
$k$ tend to infinity the assumption on $u$ hence gives, also using
\ref{eq:asympt of f}, \[
-\mathcal{M}_{\omega_{0}}(u)\leq-\mathcal{M}_{\omega_{0}}(0),\]
 which hence finishes the proof of the theorem.
\end{proof}
In particular, the proof above shows that, $\mathcal{M}_{\omega_{0}}$
is {}``convex along a geodesic'', in the sense that it is the point-wise
limit of the \emph{convex} functionals $\mathcal{F}_{k}$ along a
geodesic connecting two points in $\mathcal{H}_{\omega_{0}},$ only
using the $C^{0}-$regularity of the corresponding geodesic. Note
however that the definition of $\mathcal{M}_{\omega_{0}}$ as given
above does not even make sense unless $u_{t}$ is in $\mathcal{C}^{4}(X)$,
for $t$ fixed and $\omega_{t}>0$ (the smoothness assumption may
be relaxed to $u_{t}\in\mathcal{C}_{\C}^{1,1}(X)$ using the alternative
formula for $\mathcal{M}_{\omega_{0}}$ from \cite{ti2,ch2}). In
the case when the geodesic $u_{t}$ is assumed \emph{smooth} and $\omega_{t}>0$
the argument in the proof of the theorem above is essentially contained
in \cite{bern2}. In this latter case the convexity statement seems
to first have appeared in \cite{m} (see also \cite{d00} ). In \cite{don1}
the previous theorem was proved using the deep results in \cite{do1}
and the {}``finite dimensional geodesics'' in approximations of
$\mathcal{H}_{\omega_{0}}$ as briefly explained below.

\subsection{\label{sub:Proof-of-Theorem fano tors}Proof of Theorem \ref{thm:fano torsion}}

Set $L=-K_{X}$ so that $S=1$ and let $\mathcal{T}_{k}(u)=(\tau_{F}(e^{-ku}h_{0},\omega_{0})/N_{k}$
in terms of the notation in section \ref{sub:Analytic-torsion}.

\emph{Step 1: }For any given $\epsilon>0$ there is a metric $\omega_{0}\in c_{1}(-K_{X})$
such that \begin{equation}
\mathcal{F}_{k}(u)\leq\frac{n}{2}(1-R(X)+\epsilon)(J_{\omega_{0}}(u)+C_{\epsilon})\label{eq:pf of torsion step one}\end{equation}
 for $k\geq k_{\epsilon}.$ 

To prove this we continue with the notation in the proof of Theorem
\ref{thm:mabuchi} above. By scaling invariance we may assume that
$\sup u=0.$ Now, by the convexity of $u_{t}$ in $t$ we have $v_{0}\leq u_{1}-u_{0}=u$
and hence $-v_{0}\geq0.$ Take $t>0$ and $\omega_{0}$ such that
$\mbox{Ric}\omega_{0}>t\omega_{0}.$ Using \ref{eq:fk in terms of deriv-1}
there is a constant $C$ (depending on $t)$ such that $\mathcal{F}_{k}(u)\leq$\[
\leq\frac{n}{2}(1-t)(\frac{1}{\mbox{Vol}(-K_{X})}+Ck^{-1})\int(\frac{\omega_{0}^{n}}{n!}))(-v_{0})=-(\frac{n}{2}(1-t)+Ck^{-1})\frac{\partial\mathcal{E}_{\omega_{0}}(u_{t})}{\partial t}_{t=0}=\]
\[
=-(\frac{n}{2\mbox{Vol}(-K_{X})}(1-t)+Ck^{-1})\mathcal{E}_{\omega_{0}}(u_{1})\]
using that $\mathcal{E}_{\omega_{0}}(u_{t})$ is affine and $\mathcal{E}_{\omega_{0}}(u_{0})=0$
in the last step. Since we have assumed that $\sup u=0$ this means,
by a standard argument, that there is a constant $C$ such that \ref{eq:pf of torsion step one}
holds for \emph{any} $u$. 

\emph{Step 2:} For any $\omega_{0}\in c_{1}(-K_{X})$ there is a positive
constant $C$ such that \[
\mathcal{T}_{k}(u)\leq(-\frac{1}{2n}+Ck^{-1})(J_{\omega_{0}}(u)+C)+\mathcal{F}_{k}(u)\]

This follows immediately from Prop \ref{pro:form for torsion in general n}
combined with the basic inequality \ref{eq:ineq for i-j}.

Finally, combining the previous two steps shows that $\mathcal{T}_{k}(u)$
is bounded from above for $k$ sufficently large as long as $\frac{n}{2}(1-R(X))-\frac{1}{2n}<0,$
i.e. if \[
R(X)>1-n^{-2}.\]
This proves the boundedness with the reference metric $\omega_{0}$
depending on $t.$ The case of a general fixed reference weight then
follows immediately from the anomaly formula in \cite{b-g-s} (this
observation was already used in \cite{g-s}).

In the case when $X$ is a Fano surface (i.e. $n=2)$ it follows from
\cite{t-y} that $R(X)=1$ unless $X$ is equal to $\P^{2}$ blown
up in one or two points. In the first case $R(X)=21/28$ and in the
second case $R(X)=21/25$ \cite{li} and hence the condition on $R(X)$
above is always satisfied when $n=2.$ 

\emph{Maximizers on $\P^{2}:$ }First observe that in the case of
$\P^{2}$ the error term $O(1/k)$ and $h_{\omega_{0}}$ in the formula
in Prop \ref{pro:form for torsion in general n} vanish : \begin{equation}
(\tau_{F}(e^{-ku}h_{0},\omega_{0})-\tau_{F}(h_{0},\omega_{0}))/N_{k}=-\frac{1}{2\mbox{Vol}(-K_{X})}\left(I_{\omega_{0}}(u)-J_{\omega_{0}}(u)\right)+\mathcal{F}_{k\omega_{0}}(u),\label{eq:anal torsion for p2}\end{equation}
 Indeed, since $n=2$ the error term is of the form $\int u\beta^{2,2}+C$
for $\beta$ a form of maximal degree on $X$ depending only on $\omega_{0}.$
Now, since $\omega_{0}$ is invariant under the $SU(n+1)$ action
so is $\beta^{2,2},$ i.e. it is equal to a constant $a$ times $(\omega_{0})^{n}.$
But since the analytic torsion is invariant under $u\mapsto u+c$
it follows that $a=0.$ Finally, by the inequality \ref{eq:ineq for i-j}
the term $I_{\omega_{0}}(u)-J_{\omega_{0}}(u)$ is non-negative and
vanishes if and only if $J_{\omega_{0}}$ does, i.e. if and only if
$u$ is constant. Since, by Corollary \ref{cor:homeg} $\mathcal{F}_{k\omega_{0}}(u)\leq0$
this hence finishes the proof.

It is is also interesting to note that the following concavity property
holds:
\begin{prop}
Let $(X,\omega_{0})$ be equal to $\P^{1}$ or $\P^{2}$ equipped
with the Fubini-Study metric. Then the functional on $\mathcal{H}_{k\omega_{0}}$
defined by the analytic torsion associated to $-kK_{X}$ (with respect
$\omega_{0})$ is concave along smooth geodesics in $\mathcal{H}_{k\omega_{0}}$
for any $k\geq0.$ Moreover, if $(X,\omega_{0})$ is a Kähler-Einstein
manifold of semi-positive Ricci curvature of dimension $n\leq2$ then
the large $k$ limit of the analytic torsion functionals, i.e. minus
the relative entropy $u\mapsto\mathcal{S}(\mu_{u},\mu_{0})$ (see
remark \ref{rem:torsion entropy}) is also concave along smooth geodesics. \end{prop}
\begin{proof}
By assumption $h_{\omega_{0}}=0.$ By formula \ref{eq:anal torsion for p2}
(and its analogue in section \ref{sub:Analytic-torsion} for $n=1)$
combined with the concavity of $\mathcal{F}_{k\omega_{0}}$ established
in previous sections it is enough to observe that the functional $I_{\omega_{0}}-J_{\omega_{0}}$
is convex along smooth geodesics when $n\leq2.$ To this end first
note that when $n=1$ we have $I_{\omega_{0}}-J_{\omega_{0}}=J_{\omega_{0}}/2$
which is convex according to Prop \ref{pro:e is affine}. As for the
case $n=2$ a direct computation, using integration by parts and the
geodesic equation gives 

\[
d^{2}\left(I_{\omega_{0}}(u_{t})-J_{\omega_{0}}(u_{t})\right)/d^{2}t=\frac{1}{2\pi}\int\omega_{0}\wedge(\partial_{t}\partial_{t}u_{t}\omega_{u_{t}}-i\partial_{X}(\partial_{t}u_{t})\wedge\bar{\partial}_{X}(\partial_{t}u_{t})).\]
 Next, observe that \[
i\partial_{X}(\partial_{t}u_{t})\wedge\bar{\partial}_{X}(\partial_{t}u_{t})\leq\left|\bar{\partial}_{X}\partial_{t}u\right|_{\omega_{u_{t}}}^{2}\omega_{u_{t}}\]
using the basic fact that a positive $(1,1)$ form $\eta$ may be
point-wise estimated from above by $tr(\eta)\omega$ where $tr(\eta)$
is the trace of $\eta$ wrt the Kähler form $\omega.$ Combining this
inequality with the previous equation and the geodesic equation \ref{eq:geod eq concrete}
for $u_{t}$ finally proves the claimed convexity of $I_{\omega_{0}}(u_{t})-J_{\omega_{0}}(u_{t}).$\end{proof}
\begin{rem}
The analytic torsion funtionals are not bounded from above on the
space of\emph{ all} smooth metrics on $-K_{X}.$ Indeed, as explained
in Remark \ref{rem:j-f} there are smooth functions $u$ such that
the functionals $\mathcal{F}_{k}(tu)$ and $-(I_{\omega_{0}}(tu)-J_{\omega_{0}}(tu))$
tend to infinity when $t\rightarrow\infty.$ Hence, it follows from
formula \ref{eq:anal torsion for p2} that $\mathcal{F}_{k}(u)$ also
tends to infinity.
\end{rem}

\subsection{\label{sub:Proof-of-a conj aub}Proof of a conjecture of Aubin}

Given a polarized manifold $(X,L)$ we define an invariant $B(X,L)$
as the following infimum over all Kähler metrics $\omega\in c_{1}(L):$
\[
B(X,L):=\inf_{\omega}B(X,\omega),\,\,\, B(X,\omega_{u}):=\sup_{X}(\frac{\beta_{u}}{\omega_{u}^{n}/Vn!}).\]
Note that $B(X)\geq1$ since $\beta_{u}$ and $\omega_{u}^{n}/Vn!$
are both probability measures . When $X$ is a Fano manifold we will
write $B(X)=B(X,-K_{X}),$ i.e. \begin{equation}
B(X)=\inf_{\omega}\sup_{X}(\exp(-h_{\omega}(x))\label{eq:def of b of x}\end{equation}
 where $h_{\omega}$ is the {}``canonical'' Ricci potential of $\omega$
(i.e. $h_{\omega}$ satisfies equation \ref{eq:ricci pot eq} with
$\mu=1$ and $\int_{X}e^{h_{\omega}}\frac{\omega_{0}^{n}}{Vn!}=1).$
Repeating the argument in the proof of step one in section \ref{sub:Proof-of-Theorem fano tors},
but now using that \[
\frac{d}{dt}_{t=0+}\mathcal{F}_{\omega_{0}}(u_{t})=\int_{X}(\frac{\beta_{u}}{\omega_{u}^{n}/Vn!}-1)\omega_{u}^{n}/Vn!(-v_{0})\]
 immediately gives the following result of independent interest. 
\begin{thm}
\label{thm:general upper bound on f}For any ample line bundle $L\rightarrow X$
the following inequality holds (with the same notation as in Theorem
\ref{thm:main}): for any $\epsilon>0$ there is a constant $C_{\epsilon}$
(also depending on $\omega_{0})$ such that \[
\mathcal{F}_{\omega_{0}}(u)\leq\frac{1}{V}(B(X,L)-1+\epsilon)J_{\omega_{0}}(u)+C_{\epsilon}\]
In particular, when $X$ is Fano and $L=-K_{X}$ this means that for
any $\epsilon>0$ there is a constant $C_{\epsilon}$ such that \begin{equation}
\log\int_{X}e^{-u}\frac{\omega_{0}^{n}}{Vn!}\leq\frac{1}{V}(B(X)+\epsilon)J_{\omega_{0}}(u)+C_{\epsilon}\label{eq:aubin ineq with b}\end{equation}
 for all $u\in\mathcal{H}_{\omega_{0}}$ such that $\int_{X}u\omega_{0}^{n}=0$
\end{thm}
The second inequality \ref{eq:aubin ineq with b} above answers in
the affirmative the following conjecture of Aubin \cite{au} (in the
case when $\omega_{0}\in c_{1}(-K_{X}))$: there exists positive constants
$A$ and $B$ such that \begin{equation}
\log\int_{X}e^{-u}\frac{\omega_{0}^{n}}{Vn!}\leq\frac{A}{V}J_{\omega_{0}}(u)+B,\label{eq:conjectured ineq of aubin}\end{equation}
for all $u\in\mathcal{H}_{\omega_{0}}$ such that $\int_{X}u\omega_{0}^{n}=0$
(regardsless if $X$ admits a Kähler-Einstein metric or not). Such
inequalities were further studied by Ding \cite{din} who proved that
one may take $A=n+1$ if one restricts to all $u$ such that $\mbox{Ric }\omega_{u}>\epsilon\omega_{u}$
for some $\epsilon>0.$ He also observed that the general conjectured
inequality \ref{eq:conjectured ineq of aubin} of Aubin implies the
existence of a solution $\omega$ to Aubin's equation \[
\mbox{Ric }\omega=t\omega+(1-t)\omega_{0}\]
 for any given $t\in[0,$ $1/A^{n}[$ (see Remark 2 on p.468 in \cite{din}).
Hence, applying Ding's argument gives (since clearly $\mbox{Ric }\omega>t\omega)$
the following Corallory of the previous theorem
\begin{cor}
Let $X$ be a Fano manifold of complex dimension $n$. Then the following
relation between the two invariants $R(X)$ and $B(X),$ defined by
formulas \ref{eq:def of r} and \ref{eq:def of b of x}, respectively,
holds: \[
R(X)\geq1/B(X)^{n}\]

\end{cor}
In view of the previous inequality it would be interesting to know
whether there is a\emph{ universal upper bound }on $B(X)$ only depending
on the \emph{dimension }$n$ of the Fano manifold $X?$ As pointed
out by Tian-Yau \cite{t-y-kahler} and Ding \cite{din} such a univeral
lower bound on $R(X)$ is equivalent to a universal upper bound on
$\mbox{Vol}(-K_{X}).$ This latter bound was a well-known conjecture
in algebraic geometry at the time, subsequently proved in \cite{k-m-m}
using Mori's {}``bend-and-brake'' method of producing rational curves.

One approach to proving a universal upper bound on $B(X)$ could be
to study $\liminf_{t\rightarrow\infty}B(X,\omega_{t})$ for $\omega_{t}$
evolving according to the (normalized) Kähler-Ricci flow \cite{t-z}.
By the fundamental estimate of Perelman for $h_{\omega_{t}.}$ the
previous limit is always finite and equal to one if $X$ admits a
Kähler-Einstein metric \cite{t-z}. 

Finally, it should be pointed out that a different relation between
$R(X)$ and the properness of Mabuchi's $K-$energy has been established
by Székelyhidi \cite{sz}.

\subsection{\label{sub:Comparison-with-Donaldson's}Comparison with Donaldson's
setting and balanced metrics}

In the setting of Donaldson \cite{don1} the role of the space $H^{0}(X,L+K_{X})$
is played by the space $H^{0}(X,L).$ Any given function $u$ in $\mathcal{H}_{\omega_{0}}$
induces an Hermitian norm $Hilb(u)$ on $H^{0}(X,L)$ defined by \[
Hilb(u)[s]^{2}:=\int_{X}\left|s\right|_{\psi_{0}}^{2}e^{-u}(\omega_{u})^{n}/n!\]
 Then the functional that we will refer to as $\mathcal{L}_{D}(u),$
which plays the role of $\mathcal{L}_{\omega_{0}}(u)$ in Donaldson's
setting, is defined as in formula \ref{eq:def of l intro}, but using
the scalar product on $H^{0}(X,L)$ corresponding to $Hilb(u).$ In
other words, \begin{equation}
\mathcal{L}_{D}(u)=\mathcal{L}_{\omega_{0}+\mbox{Ric}\omega_{0}}(u-\log(\omega_{u}^{n}/\omega_{0}^{n}))\label{eq:l donaldson as l}\end{equation}
It turns out that $\mathcal{L}_{D}(u)$ is \emph{concave} along smooth
geodesics (see Theorem 3.1 in \cite{bern2} for a generalization of
this fact). However, it does not appear to be concave along a general
psh paths, which makes approximation more difficult in this setting.
Moreover, Theorem 2 in \cite{don1} says that the critical points
of $\mathcal{E}_{\omega_{0}}-\mathcal{L}_{D}$ are in fact \emph{minimizers.}
\footnote{Comparing with the notation in \cite{don1}, $\mathcal{L}_{D},$ $\mathcal{E}$
and $u$ correspond to $-\mathcal{L},$ $-I$ and $-\phi,$ respectively. %
}. On the other hand using the relation \ref{eq:l donaldson as l}
Theorem \ref{thm:main} corresponds to an\emph{ upper }bound on $\mathcal{E}_{\omega_{0}}-\mathcal{L}_{D}$
in terms of an explicit local functional (as long as $\omega_{0}+\mbox{Ric}\omega_{0}\geq0$
if $n>1).$

A major technical advantage of Donaldson's setting is that the critical
points (which are called \emph{balanced} in \cite{don1}) of the functional
$\mathcal{E}-\mathcal{L}_{D}$ acting on all of $\mathcal{C}^{\infty}(X)$
are automatically of the form \begin{equation}
\psi=\log({\textstyle \frac{1}{N}\sum_{i}}\left|S_{i}\right|^{2})\label{eq:bergm metric}\end{equation}
 fore some base $(S_{i})$ in $H^{0}(X,L).$ In particular, $u$ is
automatically in $\mathcal{H}_{\omega_{0}}$ (assuming that $L$ is
very ample). This is then used to replace the space $\mathcal{H}_{\omega_{0}}$
by the sequence of \emph{finite dimensional} symmetric spaces $GL(N,\C)/U(N)$
corresponding to the set of metrics on $L$ of the form \ref{eq:bergm metric}
(called Bergman metrics). In particular, the new geodesics, defined
wrt the Riemannian structure in the symmetric space $GL(N,\C)/U(N)$
are automatically smooth and the analysis in \cite{don1} is reduced
to this finite dimensional situation.

Note also that in this setting there is a sign difference in the expansion
\ref{eq:asym of dl}, where $\mbox{Ric \ensuremath{\omega_{u}}}$
is replaced by $-\mbox{Ric \ensuremath{\omega_{u}}. }$ As a consequence,
in Donaldson's case the functional corresponding to $\mathcal{F}_{k}$
converges to $\mathcal{M}_{\omega_{0}}$ (without the minus sign!),
which hence becomes \emph{convex} along smooth geodesics, which is
consistent with the conclusion reached above, as it must.

\subsection{\textmd{\label{sub:Conjectures}}Conjectures }

In the light of the work of Donaldson referred to above it is natural
to make the following 
\begin{conjecture}
\label{con:conv etc}Assume that $\mbox{Aut\ensuremath{_{0}(X)}}$
is trivial $c_{1}(L)$ contains a Kähler metric $\omega_{0}$ of contstant
scalar curvature. Then the following holds:
\begin{itemize}
\item for $k$ sufficently large the functional $\mathcal{F}_{k}$ admits
a critical point $u_{k}$ in $\mathcal{H}_{\omega_{0}}$
\item The convergence $\omega_{k}\rightarrow\omega_{0}$ as $k\rightarrow\infty$
holds in the $\mathcal{C}^{\infty}-$topology.
\end{itemize}
\end{conjecture}
Note that when $L$ is a degree one line bundle over a complex curve
of genus at least one, then we have already established existence
for $k=1.$ Since, loosely speaking existence also holds for $k=\infty$
(since the limiting functional is minus Mabuchi's $K-$energy whose
ciritical points are constant curvature metrics) we make the following
more precise 
\begin{conjecture}
Let $L\rightarrow X$ be an ample line bundle over a complex curve
of genus at least one. Then the functionals $\mathcal{F}_{k}$ admit
a critical point for any $k.$
\end{conjecture}
In particular, the case when $X$ has genus one would together with
the results in the present paper confirm the conjecture of Gillet-Soulé,
since $\mathcal{F}_{k}$ coincides with the corresponding analytic
torsion functional in this case.

Motivated by the case when $L=-K_{X}$ \cite{p-s+} we further make
the following
\begin{conjecture}
Assume that $\mbox{Aut\ensuremath{_{0}(X)}}$ is trivial and that
$\mathcal{F}_{\omega_{0}}$ admits a critical point in $\mathcal{H}_{\omega_{0}}.$
Then $\mathcal{F}_{\omega_{0}}$ is coercive (see remark \ref{rem:existence}).
\end{conjecture}
Note that if $\mathcal{H}_{\omega_{0}}$ were finite dimensional then
this would follow immediately from the strict concavity furnished
by Theorem \ref{thm:(Berndtsson)-Let-} (in this case the equivalence
is closely related to the Kempf-Ness principle for stability in geometric
invariant theory \cite{th}).

\subsection{\label{sub:Relation-to-random}Relation to random sections and gauge
theory}

The critical point equation and hence the conjectures above can be
given a natural interpretation in terms of random sections. Indeed,
as shown by Shiffman-Zelditch \cite{s-z} the Bergman measure $\beta_{\phi}$
of a weight $\phi$ on $L$ admits the following probabilistic interpretation
(which is different from the one discussed in the appendix): let $d\nu_{\phi}$
be the probability measure on $H^{0}(X,L+K_{X})$ defined as the Gaussian
probablity measure on the Hilbert space $(H^{0}(X,L+K_{X}),\left\langle \cdot,\cdot\right\rangle _{\phi}).$
Then \[
\beta_{\phi}=\E(s\wedge\bar{s}e^{-\phi}),\]
 where $\E$ denotes the expectation (i.e. average) over $s\in H^{0}(X,L+K_{X})$
wrt $d\nu_{\phi}$ of the measure-valued random variable $s\mapsto s\wedge\bar{s}e^{-\phi}.$
Then the critical point equation may be formulated as \begin{equation}
(dd^{c}\phi)^{n}/n!V=\E(s\wedge\bar{s}e^{-\phi})\label{eq:prob interpr critical}\end{equation}
for a positively curved metric $e^{-\phi}$ on $L.$ 

It may be illuminating to reformulate the previous equation in terms
of gauge theory \cite{d-k}. For concreteness we will only consider
the case when $X$ is a genus one curve, i.e. a one-dimensional torus
and denote by $\omega_{0}$ the standard invariant Kähler metric on
$X.$ Fix an Hermitian line bundle $E\rightarrow X$ (i.e. with structure
group $U(1))$ of positive topological degree $k.$ Consider the following
equation for a unitary connection $A$ on $E:$ 

\begin{equation}
\frac{i}{2\pi k}F_{A}=\E(|\Psi|^{2})\omega_{0}\label{eq:gauge form critic}\end{equation}
 where $F_{A}$ is the curvature two-form of $A$ and where $\E$
now denotes the expectation wrt the Gaussian probability measure on
the $k-$dimensional Hilbert space \[
\mathcal{H}_{A}=\left\{ \Psi:\,\bar{\partial}_{A}\Psi=0\right\} \]
 realized as a subspace of the space of all smooth sections $\Psi$
with values in $E$ equipped with the Hermitian product obtain by
integrating against $\omega_{0}.$ Here $\bar{\partial}_{A}$ denotes
the $(0,1)-$part of the connection $A$ (identified with a differential
operator). A standard argument (see for example \cite{d-k,br}) then
yields an isomorphism between the moduli space $\mbox{UN}_{k}(X)$
of gauge equivalence classes of solutions $A$ of equation \ref{eq:gauge form critic}
with the space of solutions $\phi$ (modulo scaling) of the critical
point equation for $\phi$ when $L$ ranges over all equivalence classes
of degree $k$ holomorphic line bundles over $X.$ In this formulation
the second point in conjecture \ref{con:conv etc} in this case is
equivalent to the convergence of $F_{A_{k}}/k$ towards $\omega_{0}$
when $A_{k}$ is in $\mbox{UN}_{k}(X).$ 

Finally, it should be pointed out that when $k=1$ equation \ref{eq:gauge form critic}
is clearly equivalent to the following equation for $(A,\Psi):$ \begin{equation}
\bar{\partial}_{A}\Psi=0,\,\,\,\,\,\frac{i}{2\pi}F_{A}=|\Psi|^{2}\omega_{0}\label{eq:jacki}\end{equation}
In the physics litterature such $(A,\Psi)$ describe the periodic
static one-vortex solutions of the Jackiw-Pi model (also called non-topological
non-relativistic self-dual \emph{Chern-Simons-Higgs solitons} or the
static \emph{gauged non-linear Schrödinger equation;} see \cite{acsv}
and references therein.) In this formulation the uniqueness result
in Corollary \ref{cor:mean field} (for $\mu=4\pi)$ amounts to saying
the solutions of \ref{eq:jacki} are uniquely determined, up to gauge
equivalence, by the zero $p$ of the corresponding {}``Higgs-field''
$\Psi.$ Note also that the equation corresponding to parameter $\mu\in]0,4\pi]$
in Corollary \ref{cor:mean field} is obtained by replacing $|\Psi|^{2}$
in equation \ref{eq:jacki} with $t|\Psi|^{2}+(1-t)$ for $\eta=4\pi t$
for $t>0.$ For $t$ \emph{negative} this equation is equivalent to
the \emph{Yang-Mills-Higgs (Ginzburg-Landau) equations }on a Riemann
surface \cite{br} (after rescaling $\Psi).$ For a general Hermitian
line bundle $L$ of degree $V>0$ and an effective divisor $D:=\sum_{i}q_{i}p_{i}$
where $p_{i}\in X$ and $q_{i}$ is a positive integer the Yang-Mills-Higgs
equations admit a unique gauge equivalence class of solutions $(A,\Psi)$
such that $D$ is the zero-divisor of $\Psi.$ Note however that for
the Chern-Simons-Higghs equation uniqueness fails when $V=2$ and
$D=2p$ since the equation is then equivalent to equation \ref{eq:intro mean field for v}
(compare the discussion below Theorem \ref{cor:mean field}).

\section{Appendix}

\subsection{\label{sub:Bergman-kernels}Bergman kernels, Toeplitz operators and
determintal point processes}

Given a function $u$ corresponding to the weight $\psi:=\psi_{0}+u$
on the line bundle $L$ we denote by $K_{u}(x,y)$ the \emph{Bergman
kernel} of the Hilbert space $(H^{0}(X,L+K_{X}),\left\langle \cdot,\cdot\right\rangle _{\psi_{0}+u}),$
i.e. \[
K_{u}(x,y):=i^{n^{2}}\sum_{i=1}^{N}s_{i}(y)\wedge\bar{s_{i}(x),}\]
represented in terms of a given orthonormal base $(s_{i})$ in $(H^{0}(X,L+K_{X}),\left\langle \cdot_{i},\cdot\right\rangle _{\psi_{0}+u}).$
This kernel may be caracterized as the integral kernel of the correponding
orthogonal projection $\Pi_{u}$ onto $(H^{0}(X,L+K_{X}),\left\langle \cdot_{i},\cdot\right\rangle _{\psi_{0}+u}),$
i.e. for any smooth section $s$ of $L+K_{X}$

\begin{equation}
(\Pi_{u}s)(x)=\int_{X_{y}}s(y)\wedge\bar{K}(x,y)e^{-(\psi(y)}\label{eq:pi in terms of k}\end{equation}
The \emph{Toeplitz operator $T[f]$ with symbol $f\in C^{0}(X),$}
acting on $(H^{0}(X,L+K_{X}),\left\langle \cdot_{i},\cdot\right\rangle _{\psi_{0}})$
(defined below formula \ref{eq:l as toeplitz det intro}) may then
be expressed as \begin{equation}
(T[f])(x)=\int_{X_{y}}f(y)s(y)\wedge\bar{K}(x,y)e^{-\psi(y)}\label{eq:toeplitz as k}\end{equation}
Applying \ref{eq:pi in terms of k} $K_{u}(x,\cdot)$ gives the following
{}``integrating out'' formula \begin{equation}
N\beta_{u}(x):=K_{u}(x,x)e^{-\psi(x)}:=\int_{X_{y}}\left|K(x,y)\right|^{2}e^{-(\psi(x)+\psi(y)}\label{eq:integr out}\end{equation}
When studying the dependence of $\beta_{u}$ on $u$ it is useful
to express $\beta_{u}(x)$ as the normalized\emph{ one-point correlation
measure }of a determinantal random point process. To this end we recall
that $\psi=\psi_{0}+u$ induces a probability measure $\gamma$ on
the $N-$fold product $X^{N}$ defined by the following local density
\cite{be2} \[
\rho(x_{1},...,x_{N})=\left|\det_{1\leq i,j\leq N}(s_{i}(x_{i}))_{i,j}\right|^{2}e^{-\psi(x_{1})}...e^{-\psi(x_{N})}/Z_{\psi}\]
 where $(s_{i})$ is an orthonormal base in the Hilbert space $(H^{0}(X,L+K_{X}),\left\langle \cdot_{i},\cdot\right\rangle _{\psi_{0}})$
and $Z_{\psi}$ is the normalizating constant (partition function).
Integration over $X^{N}$ will be denoted by respect to $\gamma$
by $\E_{\psi}$ (=expectation). As is essentially well-known \cite{be2}
$\log Z_{\psi}=-N\mathcal{L}_{\omega_{0}}(u),$ i.e. \begin{equation}
\mathcal{L}_{\omega_{0}}(u)=-\log\E_{\psi_{0}}(e^{-(u(x_{1})+...+u(x_{n})})/N\label{eq:l as expect}\end{equation}
 We also have \cite{be2} \begin{equation}
\beta_{u}(x)=\frac{1}{N}\E_{\psi}(\sum_{i=1}^{N}\delta_{x_{i}})=\int_{X^{N-1}}\left|(\det S_{0})(x,x_{2},...,x_{N}\right|^{2}e^{-\psi(x)}e^{-\psi(x_{2})}...e^{-\psi(x_{N})}/Z_{\psi}\label{eq:beta as one-pt}\end{equation}
 In particular, the map $(x,t)\mapsto(\beta_{u_{t}}(x)/\omega_{0}^{n})$
is \emph{continuous} if $u_{t}$ is a continuous path and hence there
is a positive constant $C$ such that\begin{equation}
1/C\leq(\beta_{u_{t}}(x)/\omega_{0}^{n})\leq C\label{eq:bounds on beta}\end{equation}
 on $[0,1]\times X,$ if $L+K_{X}$ is globally generated, i.e. if
$\beta_{u_{t}}(x)>0$ point-wise. Formula \ref{eq:beta as one-pt}
also shows, by the dominated convergence theorem, that $\frac{d\beta_{u_{t}}(x)}{dt}_{t=0+}$
exists under the assumptions in the following lemma. 
\begin{lem}
\label{lem:deriv of beta}Let $u_{t}$ be a family of continuous functions
on $X$ such that the right derivative $v_{t}:=\frac{du_{t}}{dt}_{+}$
exists and is uniformly bounded on $[0,1]\times X.$ Then \begin{equation}
=\frac{d\beta_{u_{t}}(x)}{dt}_{t=0+}=\int_{X_{y}}\left|K_{u}(x,y)\right|^{2}e^{-(\psi_{0}(x)+\psi_{0}(y)}v_{0}(y)-\beta_{u_{t}}(x)v_{0}(x):=R[v](x)\label{eq:r operator}\end{equation}
.\end{lem}
\begin{proof}
The proof of the formula was obtained in \cite{be} (formula $5),$
at least in the smooth case. For completeness we recall the simple
proof. By the discussion above we may differentiate formula \ref{eq:integr out}
and use Leibniz product rule to get \[
\partial_{t}(K_{t}(x,x))=2\mbox{Re }\int_{X_{y}}\partial_{t}(K_{t}(x,y)\wedge\bar{K}(x,y)e^{-\psi_{t}(y)}-\int_{X_{y}}\left|K_{t}(x,y)\right|^{2}(\partial_{t}\psi_{t}(y))e^{-(\psi_{t}(x)+\psi_{t}(y)}\]
Applying formula \ref{eq:pi in terms of k} to the holomorphic section
$s(\cdot)=\partial_{t}K_{t}(x,\cdot)$ shows that the second term
above equals $2\partial_{t}(K_{t}(x,x)).$ Hence, \[
\partial_{t}(K_{t}(x,x))=\int_{X_{y}}\left|K_{t}(x,y)\right|^{2}(\partial_{t}u_{t}(y))e^{-(\psi_{t}(x)+\psi_{t}(y)},\]
 which proves the lemma, since $N\beta_{u}(x)=K(x,x)e^{-\psi(x)}$.
\end{proof}

\subsection{A {}``Bergman kernel proof'' of the convexity statement in Theorem
\ref{thm:(Berndtsson)-Let-}}

Let $\psi_{t}:=\psi_{0}+u_{t}.$ As will be shown below, differentiating
$\mathcal{L}_{\omega_{0}}(u_{t})$ gives \begin{equation}
\partial_{t}\partial_{\bar{t}}\mathcal{L}_{\omega_{0}}(u_{t})=\frac{1}{N}\sum_{i=1}^{N}(\left\Vert (\partial_{t}\partial_{\bar{t}}u_{t})s_{i}\right\Vert _{\psi_{t}}^{2}-\left\Vert (\partial_{\bar{t}}u_{t}s_{i})-\Pi_{u_{t}}(\partial_{\bar{t}}u_{t}s_{i})\right\Vert _{\psi_{t}}^{2}),\label{eq:second der of l as proj}\end{equation}
where $(s_{i})$ is orthnormal wrt $\psi_{t}=\psi_{0}+u_{t}.$ Given
this formula the argument proceeds exactly as in \cite{bern2}; by
the definition of $\Pi_{u_{t}},$ the second term inside the sum is
the $L^{2}$- norm of the smooth solution $\alpha$ to the inhomogenous
$\bar{\partial}-$equation on $X:$ \begin{equation}
\bar{\partial_{X}}\alpha=\bar{\partial_{X}}(\partial_{t}u))s_{i},\label{eq:dbar equation appendix}\end{equation}
 which has minimal norm wrt $\left\Vert \cdot\right\Vert _{\psi_{t}}^{2}.$
Now the Hörmander-Kodaira $L^{2}-$inequality for the solution gives
\begin{equation}
i^{n^{2}}\int_{X}\alpha\wedge\bar{\alpha}e^{-\psi_{t}}\leq i^{n^{2}}\int_{X}\left|\bar{\partial_{X}}(\partial_{t}u_{t})\right|_{\omega_{u_{t}}}^{2}s\wedge\bar{s}e^{-\psi_{t}},\label{eq:horm est}\end{equation}
 using that $\omega_{u_{t}}>0.$ Hence, by formula \ref{eq:second der of l as proj},
\begin{equation}
\partial_{t}\partial_{\bar{t}}\mathcal{L}_{\omega_{0}}(u_{t})\geq\frac{1}{N}\sum_{i=1}^{N}(\left\Vert (\partial_{t}\partial_{\bar{t}}u_{t})-\left|\bar{\partial_{X}}(\partial_{t}u_{t})\right|_{\omega_{u_{t}}}^{2})s_{i}\right\Vert _{\psi_{t}}^{2}\label{eq:pf of convexity of l}\end{equation}
But since, by assumption, $(dd^{c}U+\pi_{X}^{*}\omega_{0})^{n+1}\geq0$
the rhs is non-negative (compare formula \ref{eq:expanding monge as c}),
which proves that $\mathcal{L}_{\omega_{0}}(u_{t})$ is \emph{convex
}wrt real $t.$ Note that $\mathcal{L}_{\omega_{0}}(u_{t})$ is \emph{affine}
precisely when \ref{eq:horm est} is an\emph{ equality. }By examining
the Bochner-Kodaira-Nakano-Hörmander\emph{ identity} implying the
inequality \ref{eq:horm est} one sees that the remaining term appearing
in the identity has to vanish. In turn, this is used to show that
the vector field $V_{t}$ defined by formula \ref{eq:int multip}
has to be \emph{holomorphic} on $X$ (see \cite{bern2}). Integrating
$V_{t}$ finally gives the existence of the automorphism $S_{1}$
in Theorem \ref{thm:(Berndtsson)-Let-}, as explained in section \ref{sub:Proof-of-Theorem-main}.
A slight generalization of this argument, in the one dimensional case,
will be given in the next section.

In \cite{bern2} formula \ref{eq:second der of l as proj} was derived
using the general formalism of holomorphic vector bundles and their
curvature. We will next give an alternative {}``Bergman kernel proof''.
First formula \ref{eq:r operator} and Leibniz product rule give \[
\partial_{t}\partial_{\bar{t}}\mathcal{L}_{\omega_{0}}(u_{t})=\int_{X}(\partial_{t}\partial_{\bar{t}}u_{t})\beta_{u_{t}}+\frac{d\beta_{u_{t}}}{dt}_{t=0+}(\partial_{t}u_{t})\]
Next, by formula \ref{eq:r operator} the second term may be expressed
in terms of the Bergman kernel $K_{t}(x,y)$ associated to the weight
$\psi_{t}$ as \[
\frac{1}{N}\int_{X\times X}\left|K_{t}(x,y)\right|^{2}e^{-(\psi_{t}(x)+\psi_{t}(y)}((\partial_{t}u_{t})(x)(\partial_{t}u)(y)-\int_{X}\beta(\partial_{t}u_{t})^{2},\]
 By simple and  well-known identities for Toeplitz operators this
last expression, for $t=0,$ is precisely the trace of the operator
$T[\partial_{t}u_{t}])^{2}-T[(\partial_{t}u_{t})^{2}].$ All in all
we obtain \[
\partial_{t}\partial_{\bar{t}}\mathcal{L}_{\omega_{0}}(u_{t})=\frac{1}{N}\mbox{Tr}(T[\partial_{t}\partial_{\bar{t}}u_{t}]+(T[\partial_{t}u_{t}])^{2}-T[(\partial_{t}u_{t})^{2}]),\]
for $t=0.$ Expanding in terms of an orthonormal base $s_{i}$ hence
gives \[
\partial_{t}\partial_{\bar{t}}\mathcal{L}_{\omega_{0}}(u_{t})=\frac{1}{N}\sum_{i=1}^{N}(\left\Vert (\partial_{t}\partial_{\bar{t}}u_{t})s_{i}\right\Vert _{\psi_{0}+u_{t}}^{2}+\left\Vert \Pi_{u_{t}}(\partial_{t}u_{t}s_{i})\right\Vert _{\psi_{0}+u_{t}}^{2}-\left\Vert \partial_{t}u_{t}s_{i}\right\Vert _{\psi_{0}+u_{t}}^{2}),\]
for $t=0$ (and hence for all $t$ by symmetry) which finally proves
\ref{eq:second der of l as proj}, using {}``Pythagora's theorem''.

\subsection{\label{sub:Conditions-for-equality}Conditions for equality in the
Hörmander-Kodaira estimate when $n=1$}

Denote by $\Omega^{p,q}(X,L)=\{f^{p,q}\}$ the space of all smooth
$(p,q)-$forms with values in the holomorphic line bundle $L\rightarrow X.$
We will next assume that $\dim_{:}X=1$ and use the natural identification
$\Omega^{1,q}(X,L)=\Omega^{0,q}(X,L+K_{X}).$ Given a metric on $L$
that we will write, abusing notation slightly, as $e^{-\psi}$ we
get a natural Hermitian product on $\Omega^{1,0}(X,L).$ The next
proposition is a generalization of Proposition 2.2. in \cite{bern2}
to possible degenerate curvature forms $i\partial\bar{\partial}\psi\geq0,$
in the one dimensional case.
\begin{prop}
\label{pro:horm append}Let $L$ be a line bundle over a compact complex
manifold $X$ of dimension one and let $e^{-\psi}$ be a smooth metric
on $L$ such that \textup{$i\partial\bar{\partial}\psi\geq0$ with
strict inequality almost everywhere. For any given $g^{1,1}\in\Omega^{1,1}(X,L)$
the $L^{2}(e^{-\psi})-$minimal solution of }\[
\overline{\partial}\alpha^{1,0}=g^{1,1}\]
satisfies \[
\left\Vert \alpha^{1,0}\right\Vert ^{2}:=i\int_{X}\alpha^{1,0}\wedge\overline{\alpha^{1,0}}e^{-\psi}\leq\left\Vert g^{1,\text{1}}\right\Vert ^{2}:=i\int_{X}g^{1,1}\frac{\overline{g^{1,1}}}{\partial\bar{\partial}\psi}e^{-\psi}<\infty\]
 with equality if and only if $g^{1,1}/\partial\bar{\partial}\psi$
defines a holomorphic section of $L$ over $X.$ 
\end{prop}
Next, we turn to the proof of the proposition, which is rather close
to that in \cite{bern2}, but taking advantage of the specially simple
structure of one complex dimension. As a courtesy to the reader it
will be essentially self-contained. Let $D$ denote the Chern connection
on the Hermitian holomorphic line bundle $(L,e^{-\psi})$ and decompose
$D=D^{1,0}+\overline{\partial}$ according to bidegree. In other words,
the operator $D^{1,0}$ on $\Omega^{0,0}(X,L)$ is uniquely determined
by\label{eq:comp with metric} the relation \[
\left\langle D^{1,0}f^{0,0},g^{1,0}\right\rangle :=i\int_{X}D^{1,0}f^{0,0}\wedge\overline{g^{0,1}}e^{-\psi}=-i\int_{X}f^{0,0}\overline{\bar{\partial}g^{0,1}}e^{-\psi}\]
for all $f^{0,0}$ and $g^{0,1}.$ Note that in the usual local representations
where the metric in $L$ is reprented by $e^{-\psi}$ this means that
$D^{1,0}=\partial-\partial\psi\wedge$ locally, which gives the following
basic commutation relation: \[
[\overline{\partial},D^{1,0}]=\partial\bar{\partial}\psi,\]
where $\partial\bar{\partial}\psi$ is the (non-normalized) curvature
form. Now, for any given smooth $g^{1,1}$ we consider the following
equation for $f^{0,0}$ smooth

\begin{equation}
\overline{\partial}D^{1,0}f^{0,0}=g^{1,1},\label{eq:second order equation}\end{equation}
Note that, if $f^{0,0}$ is a solution then, by \ref{eq:comp with metric},
$D^{1,0}f^{0,0}$ is orthogonal to the null space of $\bar{\partial}$
on $\Omega^{1,0}(X,L)$ and hence $u^{1,0}:=D^{1,0}f^{0,0}$ is the
$L^{2}-$minimal solution of the inhomogenous $\bar{\partial}$- equation
in Prop \ref{pro:horm append}. Next, we observe that the equuation
\ref{eq:second order equation} admits a smooth solution, which is
unique. Indeed, since $\overline{\partial}D^{1,0}$ is elliptic this
is equivalent to the null space of $\overline{\partial}D^{1,0}$ being
trivial, i.e. that the null space of $D^{1,0}$ is trivial (using
\ref{eq:comp with metric}). This latter fact follows in turn from
the following Hörmander-Kodaira identity obtained by integrating the
commutation relation above against $|f^{0,0}|^{2}$ and using \ref{eq:comp with metric}
again: for any $f^{0,0}\in\Omega^{0,0}(X,L):$ \begin{equation}
\left\Vert D^{1,0}f^{0,0}\right\Vert ^{2}+\left\Vert \overline{\partial}f^{0,0}\right\Vert ^{2}=\left\Vert f^{0,0}\right\Vert ^{2}:=i\int_{X}|f^{0,0}|^{2}\partial\bar{\partial}\psi,\label{eq:hormand-kod ident}\end{equation}
 where the last equality defines an Hermitian product on $\Omega^{0,0}(X,L).$
In particular, if $D^{1,0}f^{0,0}=0$ on $X$ then $f^{0,0}$ vanishes
on an open set, which by the identity principle for homogenous elliptic
equations means that it vanishes identically, which finishes the proof
of the existence and uniqueness. 

Next, we define the following map, introducing the (possibly degenerate)
Hermitian product on the image that makes the map into a unitary one:
\[
*:\,\,\Omega^{0,0}(X,L)\rightarrow\Omega^{1,1}(X,L),\,\,\,*f^{0,0}:=if^{0,0}\partial\bar{\partial}\psi\]
 whose inverse, defined on the image $*(\Omega^{0,0}(X,L))$ will
also be denoted by $*.$ Then \begin{equation}
\left\Vert D^{1,0}f^{0,0}\right\Vert ^{2}=i\int_{X}(\overline{\partial}D^{1,0}f^{0,0})\overline{f^{0,0}}\leq\left\Vert \overline{\partial}D^{1,0}f^{0,0}\right\Vert \left\Vert *f^{0,0}\right\Vert \leq\left\Vert \overline{\partial}D^{1,0}f^{0,0}\right\Vert \left\Vert D^{1,0}f^{0,0}\right\Vert \label{eq:pf of hormander append}\end{equation}
 (using the Cauchy-Schwartz inequality in the first inequality and
\ref{eq:hormand-kod ident} in the second one), which proves the estimate
in the proposition above upon division by $\left\Vert D^{1,0}f\right\Vert .$
Next, if equality holds in the estimate in the proposition above then
this forces equalities in \ref{eq:pf of hormander append}. Hence,
the first equality gives equality in the Cauchy-Schwartz inequality
and hence $*f^{0,0}=\overline{\partial}D^{1,0}f^{0,0}(:=g^{1,1})$
(possibly up to a multiplicative constant of unit modulus). Concretely,
this means that $g^{1,1}/i\partial\bar{\partial}\psi=f^{0,0},$ which
is smooth on all of $X$ by ellipticity. But then the second equality
above gives $\left\Vert *f\right\Vert ^{2}=\left\Vert D^{1,0}f\right\Vert ^{2}$
so that the identity \ref{eq:hormand-kod ident} forces $\overline{\partial}f^{0,0}=0,$
i.e. $\overline{\partial(}*g)=0.$ Conversly, if $\overline{\partial(}*g)=0$
then the equality in the estimate follows immediately from the identity
\ref{eq:hormand-kod ident}.
\begin{cor}
\label{cor:strict convex in append}Let $L$ be a line bundle over
a compact complex manifold $X$ of dimension one and let $e^{-\psi_{t}}$
be a one-parameter family of metrics on $L$ such that
\begin{itemize}
\item for each $t$ $\partial_{t}\partial_{t}\psi_{t}$ is in $L^{1}(X,\beta_{\psi_{t}}),$
where $\beta_{\psi_{t}}$ is the Bergman measure
\item \textup{for each $t$ $\psi_{t}$ is smooth} on $X$ with \textup{$i\partial_{X}\bar{\partial_{X}}\psi_{t}$
strictly positive almost everywhere and }
\item \textup{for all $h^{1,0}\in H^{1,0}(X,L)$ we have \[
ih^{1,0}\wedge\overline{h^{1,0}}e^{-\psi_{t}}\leq C_{t}i\partial\bar{\partial}\psi_{t}\]
 on $X$ for some constant $C_{t}$ depending on $t.$ }
\end{itemize}
\textup{Then, $d^{2}\mathcal{L}(\psi_{t})/d^{2}t=0$ at $t$ (see
definition }\ref{eq:l as toeplitz det intro}\textup{) if and only
if there is $h^{1,0}\in H^{1,0}(X,L)$} such that $h^{1,0}\wedge\bar{\partial}(\partial_{t}\psi_{t})/\partial_{X}\bar{\partial_{X}}\psi_{t}$
defines a holomorphic section of $L$ over $X.$\end{cor}
\begin{proof}
Under the assumptions above the proof in the previous section carries
over essentially verbatim to give equality in \ref{eq:pf of convexity of l}.
In particular, we can take $h^{1,0}=s_{i}$ for some $i$ and apply
the previous proposition to $g^{1,1}:=h^{1,0}\wedge\bar{\partial}(\partial_{t}\psi_{t})$
which, by assumption, satisfies $\left\Vert g^{1,1}\right\Vert ^{2}:=i\int_{X}g^{1,1}\overline{g^{1,1}}e^{-\psi}/\partial\bar{\partial}\psi<\infty.$ \end{proof}

\end{document}